\documentclass[a4paper]{amsart}
\usepackage{amsmath,amsthm,amssymb}
\usepackage{graphicx}
\usepackage{hyperref} 

\newtheorem{theorem}{Theorem}[section]
\newtheorem{proposition}[theorem]{Proposition}
\newtheorem{lemma}[theorem]{Lemma}
\newtheorem{corollary}[theorem]{Corollary}

\theoremstyle{definition}
\newtheorem{definition}[theorem]{Definition}
\newtheorem{remark}[theorem]{Remark}

\newtheorem{example}[theorem]{Example}
\newtheorem{question}[theorem]{Question}

\newtheorem{convention}[theorem]{Convention}

\theoremstyle{theorem}

\title[Corks, exotic 4-manifolds and genus functions]{Corks, exotic 4-manifolds and genus functions}
 
\author[Kouichi Yasui]{Kouichi Yasui}
\date{January 31, 2025. \textit{Revised}: December 25, 2025.}

\address{Department of Pure and Applied Mathematics, Graduate School of Information Science and Technology, The University of Osaka, 
1-5 Yamadaoka, Suita, Osaka 565-0871, Japan}
\email{kyasui@ist.osaka-u.ac.jp}
\begin{document}

\begin{abstract} 
			We prove that every 4-dimensional oriented handlebody without 3- and 4-handles can be modified to admit infinitely many exotic smooth structures, and moreover prove that their  genus functions are pairwise equivalent. We furthermore show that for any 4-manifold admitting an embedding into a symplectic 4-manifold with weakly convex boundary, its genus function is algebraically realized as those of infinitely many pairwise exotic 4-manifolds. In addition, we prove that algebraic inequivalences of genus functions are stable under connected sums and boundary sums with a certain type of 4-manifolds having arbitrarily large second Betti numbers. Besides, we introduce a notion of genus function type for diffeomorphism invariants, and show that any such invariant shares properties similar to all the preceding results and yields lower bounds for the values of genus functions. As an application of our exotic 4-manifolds, we also prove that for any (possibly non-orientable) 4-manifold, every submanifold of codimension at most one satisfying a mild condition can be modified to admit infinitely many exotically knotted copies. 
\end{abstract}

\maketitle

\section{Introduction}\label{sec:intro}
The genus function of a 4-manifold is a powerful invariant of smooth structures and has been extensively studied dating back to the work  of Kervaire and Milnor \cite{KM61} in 1961. In general, it is quite hard to determine values of genus functions, and they are completely determined only for few 4-manifolds (see \cite{Law92},  \cite{Law97}, \cite{DL}). Nevertheless, genus functions can detect (infinitely many) exotic smooth structures for numerous 4-manifolds (e.g.\ \cite{Law97}, \cite{AY13}, \cite{Y11}, \cite{G17_genus}, \cite{AY_JSG14}, \cite{Y14}, \cite{Y19Tr}) and have various applications to low dimensional topology (e.g.\ \cite{AM97}, \cite{AM98}, \cite{LMY}, \cite{Y15}, \cite{Y19Tr}, \cite{Y19GT}, \cite{KNY}, \cite{ACMPS}, \cite{Tor}, \cite{DHHRS}). 

Lawson conjectured in his 1997 survey \cite{Law97} that, for simply connected compact 4-manifolds, genus functions determine smooth structures. Indeed, for rational surfaces and ruled surfaces, genus functions uniquely determine their smooth structures among symplectic 4-manifolds (Proposition~\ref{prop:genus_symplectic}). 
It later turned out that his conjecture is false (\cite{AR}, \cite{G17GT}), but the counterexamples are not particularly interesting, because their second homology groups are trivial, and hence the underlying topological structures uniquely determine genus functions. Thus, the conjecture remains intriguing especially for 4-manifolds with large second Betti numbers, and naturally raises the problem whether a 4-manifold can simultaneously admit exotic structures with the same genus function and those with distinct genus functions.  

Let us recall that the genus function of a connected oriented smooth 4-manifold $X$ is defined to be the map $g_X:H_2(X)\to \mathbb{Z}$ which sends a class $\alpha$ to the minimum genus of a smoothly embedded surface representing $\alpha$.  Genus functions of two 4-manifolds will be called \textit{equivalent} if there exists a homeomorphism between them whose induced isomorphism preserves the genus functions. They will be called \textit{algebraically equivalent} if there exists an isomorphism between their second homology groups that preserves their genus functions and intersection forms. 

In this paper, we prove that every 4-manifold with boundary satisfying mild conditions can be modified to admit infinitely many exotic smooth structures. Applying this result, we answer the aforementioned problems on genus functions and also modify submanifolds of 4-manifolds into exotically knotted ones. Furthermore, we reveal several interesting properties of genus functions. Besides, we introduce a notion of genus function type for diffeomorphism invariants of 4-manifolds, and show that any such invariant shares many properties of genus functions.

 \begin{convention}\label{sec:intro:convention}
 Throughout this paper, unless otherwise stated, all manifolds, submanifolds, boundaries, handlebodies, and surfaces are assumed to be compact, connected, oriented, and smooth, and all embeddings of manifolds are assumed to be smooth but not necessarily proper. A pair of 4-manifolds is called \textit{exotic} (resp.\ \textit{orientedly exotic})  if they are orientedly homeomorphic but not unorientedly diffeomorphic (resp.\ not orientedly diffeomorphic). Here, ``orientedly'' and ``unorientedly'' respectively mean ``orientation-preserving'' and ``not necessarily orientation-preserving''.
\end{convention}
  We require exotic 4-manifolds to be orientedly homeomorphic, since this assumption is necessary to ensure that their connected sums and boundary sums with a common 4-manifold always yield homeomorphic 4-manifolds. We note that boundary 3-manifolds of two exotic 4-manifolds are always orientedly diffeomorphic. 
\subsection{Exotic 4-manifolds and genus functions}\label{sec:intro:subsec:exotic}
For convenience, we introduce a notion of topological similarity for 4-manifolds, which we call \textit{HIHC-equivalence}.    

\begin{definition}Two 4-manifolds with (possibly disconnected) boundary are called \textit{HIHC-equivalent} if they are homotopy equivalent, their intersection forms are isomorphic, and their boundaries are homology cobordant. 
\end{definition}
We note that two orientedly homeomorphic 4-manifolds are necessarily HIHC-equivalent, and the converse also holds for simply connected closed (smooth) 4-manifolds due to Freedman's celebrated theorem.  

Let us recall that a \textit{$2$-handlebody} is a handlebody consisting of finitely many handles whose indices are at most two. $2$-handlebodies form a very large class of 4-manifolds with boundary, realizing arbitrarily large integers as their second Betti numbers $b_2$. (Note that in some papers, 2-handlebodies mean handlebodies consisting only of  0- and 2-handles, which form a much smaller class.) 

We prove that every 4-dimensional 2-handlebody can be modified to admit infinitely many exotic smooth structures by using corks,  and moreover show that these exotic structures have pairwise equivalent genus functions. 
\begin{theorem}\label{sec:intro:thm:2-handlebody:infinite}	For any $4$-dimensional $2$-handlebody $X$, there exist infinitely many pairwise exotic $4$-manifolds which are HIHC-equivalent to $X$. Moreover, they can be chosen so that their genus functions are pairwise equivalent. 
\end{theorem}

Even without the conclusion on genus functions, this result is new and shows the existence of infinitely many exotic smooth structures for a very large class of 4-manifolds, substantially strengthening the following earlier results: every 2-handlebody with non-trivial second homology group can be modified to admit \textit{arbitrarily many}  exotic structures (Akbulut and the author \cite{AY13});  there exist contractible 4-manifolds admitting infinitely many exotic structures (Gompf \cite{G17GT}). 

The corollary below is straightforward from the above theorem. 
 
\begin{corollary}\label{sec:intro:cor:invariant:infinite}$(1)$ For any integral symmetric bilinear form $Q$, there exist infinitely many pairwise exotic simply connected $4$-manifolds with pairwise equivalent genus functions such that their intersection forms are isomorphic to $Q$. 

$(2)$ For any finitely presented group $G$, there exist infinitely many pairwise exotic  $4$-manifolds with pairwise equivalent genus functions such that their fundamental groups are isomorphic to $G$. 

$(3)$ For any closed $3$-manifold $M$, there exist infinitely many pairwise exotic simply connected $4$-manifolds with pairwise equivalent genus functions such that their boundaries are homology cobordant to $M$.  
\end{corollary}

We furthermore show that a given 4-manifold embedded in a symplectic $4$-manifold with weakly convex boundary can be modified to admit infinitely many exotic smooth structures. We moreover prove that the genus function of the given 4-manifold is algebraically realized as their genus functions, showing that there are various genus functions which cannot determine smooth structures. 

\begin{theorem}\label{sec:intro:thm:symplectic:infinitely}
	For any $4$-manifold $X$ admitting an embedding into a symplectic $4$-manifold with weakly convex boundary, there exist infinitely many pairwise exotic $4$-manifolds which are HIHC-equivalent to  $X$. Moreover, they can be chosen so that their genus functions are pairwise equivalent and are algebraically equivalent to that of $X$. 
\end{theorem}

The assumption of this theorem holds for a large class of smooth 4-manifolds. 
For example, blow ups of 4-manifolds admitting Stein structures and their codimension-$0$ submanifolds satisfy the assumption. Furthermore, these 4-manifolds provide varieties of genus functions even within fixed homeomorphism types  (\cite{AY13}, \cite{Y14}, \cite{Y19Tr}).  
We here state only simple examples. 

\begin{corollary}\label{sec:intro:cor:nCP2} $(1)$ For each positive integer $n$, there exist infinitely many pairwise exotic simply connected $4$-manifolds $($with boundary$)$ such that their genus functions are pairwise equivalent and are algebraically equivalent to that of $n\mathbb{C}P^2$.  
	
	$(2)$ For any $D^2$-bundle over a closed surface, there exist infinitely many pairwise exotic $4$-manifolds HIHC-equivalent to  the total space $X$ of the bundle such that their genus functions are pairwise equivalent and are algebraically equivalent to that of $X$. 
\end{corollary}

The claim (1) contrasts with the aforementioned fact that the genus function of  $\mathbb{C}P^2$ determines its smooth structure among closed symplectic 4-manifolds. 

We next show that a fixed  4-manifold can simultaneously admit exotic smooth structures with pairwise equivalent genus functions and  those with pairwise inequivalent genus functions. In fact, we prove the much stronger  theorem below. 

\begin{theorem}\label{sec:intro:thm:Stein nuclei:infinitely:infinitely}
There exist infinitely many pairwise exotic simply connected $4$-manifolds with pairwise algebraically inequivalent genus functions such that, for each $4$-manifold $X$ among them, there exist infinitely many pairwise exotic $4$-manifolds homeomorphic to $X$ whose genus functions are all equivalent to that of $X$. 
\end{theorem}

We can construct similar infinite families for varieties of homeomorphism types. For example, any finitely presented group is realized as the fundamental group of infinite families of  such exotic 4-manifolds. 
Here we give similar examples for much more homeomorphism types, by relaxing the condition ``infinitely many'' to ``arbitrarily many''. 
\begin{theorem}\label{sec:intro:thm:2-handlebody:infinitely:arbitrarily}
		For any $4$-dimensional $2$-handlebody $X$ with non-trivial second homology group, there exist arbitrarily many pairwise exotic $4$-manifolds HIHC-equivalent to $X$ and having pairwise algebraically inequivalent genus functions such that, for each $4$-manifold $X'$ among them, there exist infinitely many pairwise exotic $4$-manifolds homeomorphic to $X'$ whose genus functions are all  equivalent to that of $X'$. 
\end{theorem}

We remark that exotic 4-manifolds in Theorems~\ref{sec:intro:thm:2-handlebody:infinite}, \ref{sec:intro:thm:symplectic:infinitely} and  \ref{sec:intro:thm:2-handlebody:infinitely:arbitrarily} have more properties inherited from a given 4-manifold $X$, see Theorems~\ref{sec:exotic:thm:2-handlebody:infinite},  \ref{sec:exotic:thm:symplectic:infinitely} and \ref{sec:exotic genus:thm:2-handlebody:infinitely:arbitrarily}. 

It remains interesting to see whether exotic closed 4-manifolds with equivalent genus functions exist. For closed 4-manifolds, we instead show that algebraic equivalent classes do not determine homeomorphism types even within 4-manifolds having isomorphic cohomology rings and that there are many such 4-manifolds, realizing arbitrarily large $b_2$. See Theorem~\ref{sec:genus:sum:thm:closed:algebraically equivalent} for details. During the preparation of this paper, Stipsicz and Szab\'{o} \cite{SS24} proved that infinitely many pairwise non-diffeomorphic closed 4-manifolds with the cohomology ring of $S^2\times S^2$ obtained in \cite{FPS} have pairwise algebraically equivalent genus functions. It is not known whether their 4-manifolds are pairwise  homeomorphic. They showed algebraic equivalence by explicitly determining all the values of the genus functions, but our proof is completely different from their approach. 

We remark that, after the first version of this paper appeared on arXiv on January 2025, Nakamura submitted the paper \cite{Nak} to arXiv which gives an example of Corollary~\ref{sec:intro:cor:invariant:infinite}.(1) in the case where the intersection form is represented by the $1\times 1$ matrix $0$, by using a different method.  

\subsection{Stabilities of algebraic inequivalences of genus functions}\label{sec:intro:subsec:stability}
To show that our exotic 4-manifolds have equivalent genus functions, we prove that attachments of a class of cobordisms called  \textit{quasi-invertible cobordisms} always preserve the genus functions, which is of independent interest (see Sections~\ref{sec:quasi-invertible} and \ref{sec:Genus functions and quasi-invertible cobordisms}).  As an application, we discuss behavior of genus functions under connected sums and boundary sums with a certain type of 4-manifolds. 

For an $S^2$-link $L$ in $S^4$, let $S^4(L)$ denote the closed 4-manifold obtained from $S^4$ by performing the surgery on $L$, that is, by removing the interior of the tubular neighborhood $S^2\times D^2$ for each component of $L$ and gluing $S^1\times D^3$ along the resulting boundary. A simple example of $S^4(L)$ is $n(S^1\times S^3)$. 
As we will show in Subsection~\ref{sec:genus:sum:subsec:example}, each $S^4(L)$ has rich examples of codimension-$0$ submanifolds, realizing arbitrarily large $b_2$. For instance, $\Sigma_g\times D^2$ $(g\geq 0)$ with 2-handles attached along a strongly slice link (see Definition~\ref{sec:quasi-invertible:def:strongly_slice}) with 0-framings  are such examples, where $\Sigma_g$ denotes the closed surface of genus $g$.  

In the rest of this subsection, for each $i=1,2$, let $X_i$ be a $($possibly closed$)$ 4-manifold,  $L_i$ be an $S^2$-link in $S^4$, and $Z_i$ be a connected summand or a codimension-$0$ submanifold of $S^4(L_i)$. 

We will prove that the values of the genus functions of 4-manifolds are preserved under connected sums and boundary sums with $Z_1$ and $Z_2$ (Theorem~\ref{sec:genus:sum:thm:genus preserving}). Applying this result, we show that algebraic inequivalences of genus functions are stable under these operations. In the case where $H_2(Z_i)=0$, we prove the following theorem. 

\begin{theorem}\label{sec:intro:thm:sum:iff}Suppose $H_2(Z_1)=H_2(Z_2)=0$. Then the following hold. 
	
	$(1)$ The genus functions of $X_1\#Z_1$ and  $X_2\#Z_2$ are algebraically equivalent  if and only if those of $X_1$ and $X_2$ are algebraically equivalent. 
	
	$(2)$ In the case where each $\partial X_i$ and $\partial Z_i$ are non-empty, the genus functions of $X_1\natural Z_1$ and  $X_2\natural Z_2$ are algebraically equivalent if and only if those of $X_1$ and  $X_2$ are algebraically equivalent. 
\end{theorem}

This result contrasts with the fact that the Seiberg--Witten invariants of closed symplectic 4-manifolds vanish after taking a connected sum with $S^1\times S^3$ and hence cannot detect them after the connected sum. 

In the general case, we prove similar stabilities under additional assumptions.    
\begin{theorem}\label{sec:intro:thm:equivalent:non-degenerate:sum}
	Suppose that the intersection form of each $X_i$ is non-degenerate and that each $H_2(Z_i)$ is torsion-free. Then the following hold. 
	
	$(1)$  If the genus functions of $X_1$ and $X_2$ are algebraically inequivalent, then those of $X_1\#Z_1$ and $X_2\#Z_2$ are algebraically inequivalent. 
	
	$(2)$ Assume each $\partial X_i$ and $\partial Z_i$ are non-empty. If the genus functions of $X_1$ and $X_2$ are algebraically inequivalent, then those of $X_1\natural Z_1$ and $X_2\natural Z_2$ are algebraically inequivalent. 
\end{theorem}

Furthermore, we give similar results and refinements without assuming the torsion-free conditions (Theorems~\ref{sec:genus:sum:thm:equivalent:non-degenerate} and \ref{sec:genus:sum:thm:equivalent:non-degenerate:torsion-free}). However, the above theorem does not hold if the non-degeneracy assumption is removed (Remark~\ref{sec:genus:remark:sum}). 

These stabilities for algebraic inequivalences show that for any family of exotic 4-manifolds with algebraically inequivalent genus functions, the family remain (orientedly) exotic under connected sums and boundary sums with homeomorphic 4-manifolds satisfying the assumptions.  One can apply these stabilities to various families of exotic 4-manifolds (e.g.\ \cite{AY13}, \cite{Y11}, \cite{Y14}, \cite{Y19Tr}) and produce further examples of exotic families by connected sums and boundary sums. Furthermore, we give similar stabilities for attachments of quasi-invertible cobordisms, see Subsection~\ref{subsec:Genus functions and quasi-invertible cobordisms:stabilities}. We remark that the proofs of these stabilities use no gauge theory.

As an application of our arguments, we give exotic 4-manifolds which are stably exotic under boundary sums with arbitrary Stein manifolds having non-degenerate intersection forms.  
\begin{theorem}\label{sec:intro:thm:sum:Stein}
	For any positive integer $n$, there exist pairwise exotic simply connected $4$-manifolds $X_1,X_2,\dots, X_n$ admitting Stein structures such that, for any pairwise orientedly homeomorphic $4$-manifolds $Y_1,Y_2,\dots, Y_n$ having non-degenerate intersection forms and admitting Stein structures, the boundary sums $X_1\natural Y_1$, $X_2\natural Y_2$, $\dots$, $X_n\natural Y_n$ remain pairwise exotic. 
\end{theorem}
In fact, we realize various homeomorphism types including not simply connected ones as those of $X_1,\dots, X_n$, see Theorem~\ref{sec:exotic:thm:sum:Stein}. Note that the summands $Y_1,\dots, Y_n$ can be taken as pairwise exotic 4-manifolds due to a result of \cite{AY13} (see also Theorem~\ref{sec:cork:exotic Stein}).  Although the above theorem assumes that each $Y_i$ admits a Stein structure, the above theorem holds under  the weaker condition that each $Y_i$ admits an embedding into a symplectic 4-manifold with weakly convex boundary. This will be discussed in the paper \cite{Y_sta} announced in \cite{Y19GT}. Moreover, in \cite{Y_sta}, we will produce various families of exotic 4-manifolds which are stably exotic under connected sums with arbitrary closed definite 4-manifolds.  
\subsection{Diffeomorphism invariants of genus function type}
In a different direction, we introduce a notion of \textit{genus function type} for diffeomorphism invariants of 4-manifolds by extracting a simple property of genus functions, see Definition~\ref{sec:genus type:def:invariant}. We then give their examples and show that any such invariant shares many properties of genus functions. 

Specifically, we give examples of this type which can detect families of exotic 4-manifolds, by utilizing the Seiberg--Witten invariants of closed 4-manifolds. The lasagna $s$-invariant recently introduced by Ren and Willis \cite{RW} is also an example of this type according to their result. See Example~\ref{sec:genus type:example:invariant} for these examples. 

Furthermore, we show that every invariant of genus function type yields lower bounds for values of genus functions, which can be sharp depending on the invariant, see Proposition~\ref{sec:genus type:prop:lower_bound}. We also prove stabilities of algebraic inequivalences similar to Theorems~\ref{sec:intro:thm:sum:iff} and \ref{sec:intro:thm:equivalent:non-degenerate:sum}, see Theorems~\ref{sec:genus type:thm:quasi-invertible:many_results} and \ref{sec:genus type:thm:sum:many_results}. 
 These stabilities imply that if (algebraic equivalent classes of) an invariant of genus function type can detect a family of orientedly exotic 4-manifolds, then the invariant can detect infinitely many families of orientedly exotic 4-manifolds, see Corollary~\ref{sec:genus type:cor:detect_infinite}. By contrast, we also prove that there exist various families of exotic 4-manifolds which can not be detected by any invariant of genus function type, see Theorems~\ref{sec:genus type:thm:2-handlebody:undetectable} and \ref{sec:genus type:thm:symplectic:undetectable}. We will give further related results in Section~\ref{sec:genus type}.  

\subsection{Exotically knotted  submanifolds}
Applying our exotic 4-manifolds, we also discuss exotically  knotted and exotically embedded submanifolds. 
Here a pair of submanifolds of a manifold will be called \textit{exotically knotted} if  they are diffeomorphic, topologically ambient isotopic, but not smoothly ambient isotopic. Also, a pair of submanifolds will be called \textit{exotically embedded} if they are diffeomorphic, the ambient manifold admits a self-homeomorphism that maps one to the other, but the ambient manifold does not admit such a self-diffeomorphism.  
We furthermore require natural conditions for orientations and boundaries as well (see Definition~\ref{sec:knotted:def:exotically embedded}), but our results still hold under the above simpler conditions. 
 
 Many examples of exotically knotted submanifolds are also exotically embedded, but interestingly, exotically knotted submanifolds are not necessarily exotically embedded (e.g.\ \cite{Sch19}, \cite{KMT_ar24}).  In this paper, we focus on submanifolds which are both exotically knotted and exotically embedded. Such examples indicate strong differences between smooth and topological categories. We remark that our results in this subsection are new as exotically knotted submanifolds and also as exotically embedded submanifolds. 

It is a major problem in 4-dimensional topology to see whether a given submanifold of a 4-manifold admits exotically knotted and/or exotically embedded copies. Many examples of submanifolds admitting such copies  have been constructed (see, for example, \cite{FKV88}, \cite{FS97}, \cite{AY13}, \cite{AKMRS19}, \cite{Sch19}, \cite{IKMT25}, \cite{Miy}, \cite{KMT_ar24} and the references therein), but this problem remains open in general. Utilizing our exotic 4-manifolds, here we show that every codimension-$0$ submanifold satisfying a mild condition on the exterior (i.e., the closure of the complement) can be modified to admit infinitely many exotically knotted and exotically embedded copies. 
\begin{theorem}\label{sec:intro:thm:exotic:embedded}
		Let $X$ be a $4$-manifold embedded in a  $4$-manifold $Z$. These manifolds may be non-orientable and may have $($possibly disconnected$)$ boundaries. Suppose the exterior of $X$ either is a  $($orientable$)$ $2$-handlebody or admits an embedding into a symplectic $4$-manifold with weakly convex boundary, and assume the boundary of the exterior is connected. Then, $Z$ admits infinitely many pairwise exotically knotted and exotically embedded codimension-$0$ submanifolds which are HIHC-equivalent to  $X$. 
\end{theorem}
 Furthermore, we show that each of these exotically knotted submanifolds contains the given submanifold $X$ as a submanifold and admits a deformation retraction onto $X$, indicating that the modification to $X$ is small, see Theorem~\ref{sec:knotted:thm:exotic:embedded} and also  Remark~\ref{sec:knotted:rem:exotic:embedded}.  Note that the modification is necessary, since the 4-ball admits neither exotically knotted nor exotically embedded copies due to the disk theorem. The above theorem improves an earlier result of Akbulut and the author \cite{AY13}, which modified $X$ into \textit{arbitrarily many} exotically embedded submanifolds additionally assuming that the exterior is a 2-handlebody with non-trivial second homology group. Their examples were not shown to be exotically knotted in general. 

Since every 4-manifold has numerous 2-handlebodies as their submanifolds, the corollary below is straightforward.  
\begin{corollary}\label{sec:intro:cor:exotic:embedded:all}Every $($possibly non-orientable$)$ $4$-manifold admits infinitely many $($possibly non-orientable$)$ codimension-$0$ submanifolds which are pairwise exotically knotted and exotically embedded. 
\end{corollary}

We also discuss codimension-1 submanifolds by applying Theorem~\ref{sec:intro:thm:exotic:embedded}. Prior results in the literature are realizations of many 3-manifolds as exotically knotted codimension-1 submanifolds of some orientable 4-manifolds (see e.g.\ \cite{IKMT25} and the references therein). By contrast, we show that for a given (possibly non-orientable) 4-manifold, every codimension-$1$ separating submanifold satisfying a mild condition can be modified to admit infinitely many exotically knotted (and exotically embedded) copies, using homology cobordisms embedded in the given 4-manifold.    
\begin{theorem}\label{sec:intro:thm:exotic:embedded:3-manifold}Let $M$ be a closed $3$-manifold embedded in the interior of   a $($possibly non-orientable$)$ $4$-manifold $Z$ which may have $($possibly disconnected$)$ boundary. Suppose the exterior of the tubular neighborhood of $M$ has a connected component with connected boundary which either is a $4$-dimensional $($orientable$)$ $2$-handlebody or admits an embedding into a symplectic $4$-manifold with weakly convex boundary. Then, $Z$ admits infinitely many pairwise exotically knotted and exotically embedded codimension-$1$ submanifolds which are homology cobordant to $M$.  
\end{theorem}

As corollaries,  we show that  there exist infinitely many closed $3$-manifolds such that every $($possibly non-orientable$)$ $4$-manifold realizes each of them as infinitely many exotically knotted submanifolds. 
\begin{corollary}\label{sec:intro:cor:exotic:embedded:3-manifold:S4}
For any  closed $3$-manifold $M$ embedded in $S^4$,  every $($possibly non-orientable$)$ $4$-manifold admits infinitely many pairwise exotically knotted and exotically embedded codimension-$1$ submanifolds which are homology cobordant to $M$. 
\end{corollary}

\begin{corollary}\label{sec:intro:cor:exotic:embedded:3-manifold:any}
		There exist infinitely many pairwise not homology cobordant closed $3$-manifolds such that, for each $3$-manifold $M$ among them, every $($possibly non-orientable$)$ $4$-manifold admits infinitely many pairwise exotically knotted and exotically embedded submanifolds diffeomorphic to $M$. 
\end{corollary}
Note that many closed 3-manifolds do not have the above property. Indeed, there are many closed 3-manifolds not admitting any embedding into $S^4$ (see e.g.\ \cite{IM20} and the references therein).

For a given closed 3-manifold $M$, we consider the question of what is the minimum second Betti number of a simply connected closed 4-manifold which contains a pair of exotically knotted (or exotically embedded) submanifolds diffeomorphic to $M$ (cf.\ \cite[Definition~2.2]{AGL}, \cite[Question~7.7]{IKMT25}). 
This question naturally arises from the well-known fact that  every closed 3-manifold admits an embedding into a simply connected closed 4-manifold (see e.g.\ \cite[Theorem~2.1]{AGL}).  
Here we show that every non-negative even integer is realized as such a minimum second Betti number by using  Theorem~\ref{sec:intro:thm:exotic:embedded:3-manifold}. 

\begin{corollary}\label{sec:intro:cor:exotic:embedded:3-manifold:some_4-manifold}
	For each even integer $n\geq 0$, there exist infinitely many pairwise not homology cobordant closed $3$-manifolds such that, for each $3$-manifold $M$ among them, a simply connected closed $4$-manifold with $b_2=n$ contains infinitely many pairwise exotically knotted and exotically embedded submanifolds diffeomorphic to $M$, but no simply connected closed $4$-manifold with $b_2 < n$ contains a submanifold diffeomorphic to $M$. 
\end{corollary}
In fact, we show that the above $4$-manifold with $b_2=n$ can be taken as $\frac{n}{2}(S^2\times S^2)$, independently of the choice of the 3-manifold $M$. 

We close this subsection by posing a question. As is well-known, every 4-manifold admits no pair of exotically knotted codimension-3 submanifolds. By contrast,  Corollary~\ref{sec:intro:cor:exotic:embedded:3-manifold:any} gives closed 3-manifolds such that every 4-manifold admits pairwise exotically knotted and exotically embedded submanifolds diffeomorphic to them. Thus, the following question naturally arises for the remaining codimensions.  

\begin{question}
	For $n=2,4$,  does there exist an $n$-manifold such that every 4-manifold admits a pair of exotically knotted (or exotically embedded) submanifolds diffeomorphic to the $n$-manifold?
\end{question}

This question is open even if ambient 4-manifolds are restricted to be simply connected and closed.  Corollary~\ref{sec:intro:cor:exotic:embedded:all} gives a partial answer to the $n=4$ case, but to the best of the author's knowledge, it is not known whether the corresponding result holds in the $n=2$ case. 

\subsection*{Acknowledgements}Parts of this work were presented at conferences MSJ Autumn Meeting 2023, Four Dimensional Topology, and Differential Topology '24. The author would like to thank the organizers and participants for useful discussions. The author was partially supported by JSPS KAKENHI Grant Numbers 19K03491, 23K03090 and 19H01788.
\section{Genus functions}\label{sec:genus}
In this section, we introduce notions of equivalences for genus functions, as well as a notion of the torsion-free genus function. We then prove that the values of the genus function do depend on the torsion part of the second homology group, and hence the genus function is not necessarily determined by the torsion-free genus function. Finally, we show that genus functions can determine smooth structures among symplectic 4-manifolds. 

We here fix notations of this paper. Every homology group is integral coefficient. For a continuous map $f:X\to Y$, the induced homomorphism $H_i(X)\to H_i(Y)$ is denoted by $f_*$. For a manifold $Z$, the manifold equipped with the reverse orientation is denoted by $\overline{Z}$. 

Throughout this section, let $X,Y$ be 4-manifolds which may have possibly disconnected boundaries. We recall the definition of the genus function. 

\begin{definition}For a second homology class $\alpha$ of $X$, the integer $g_X(\alpha)$ is defined to be the minimum genus of a closed surface smoothly embedded in $X$ that represents $\alpha$. The map $g_X:H_2(X)\to \mathbb{Z}$ is called the \textit{genus function} of $X$.
\end{definition}
This function is well-defined, because every second homology class of $X$ is represented by a closed surface (see \cite[Exercise 4.5.12(b)]{GS}). We occasionally treat disconnected 4-manifolds as well, and for such 4-manifolds, we allow a representative surface to be disconnected and define the genus of such a surface as the sum of the genera of its connected components. Then, we define the genus function of a disconnected 4-manifold in the same way as above. Note that, for a (connected) 4-manifold, these two definitions give the same function as easily seen.   
 
We will say that a continuous map $f:X\to Y$ \textit{preserves the genus functions} if the induced homomorphism $f_*:H_2(X)\to H_2(Y)$ preserves the genus functions, that is, $f_*$ satisfies $g_X(\alpha)=g_Y(f_*(\alpha))$ for any $\alpha\in H_2(X)$. 
We here introduce a notion of equivalence for genus functions of 4-manifolds.

\begin{definition} The genus functions of $X$ and $Y$ are called \textit{equivalent} if there exists a homeomorphism $f:X\to Y$ that  preserves the genus functions. If $f$ furthermore preserves the orientations of $X$ and $Y$, then the genus functions are called \textit{orientedly equivalent}. 
\end{definition}

It is easy to see that if two 4-manifolds are diffeomorphic  (resp.\ orientedly diffeomorphic), then their genus functions are equivalent (resp.\ orientedly equivalent). Hence, the genus function is a diffeomorphism invariant of 4-manifolds. 

We also introduce a notion of algebraic equivalence for genus functions of 4-manifolds. We say that a homomorphism $\varphi:H_2(X)\to H_2(Y)$ preserves the intersection forms of $X$ and $Y$ if $\varphi(\alpha)\cdot \varphi(\beta)=\alpha\cdot \beta$ for any $\alpha,\beta\in H_2(X)$, where $\alpha\cdot \beta$ denotes the intersection number of $\alpha,\beta$. 

\begin{definition}(1) The genus functions of $X$ and $Y$ are called \textit{algebraically equivalent} if there exists an isomorphism $H_2(X)\to H_2(Y)$ that preserves the genus functions and the intersection forms. 

(2) The genus functions of $X$ and $Y$ are called \textit{algebraically equivalent for some orientations} if $g_X$ is algebraically equivalent to $g_Y$ or to $g_{\overline{Y}}$. 

(3) For subgroups $F_X\subset H_2(X)$ and $F_Y\subset H_2(Y)$, the restrictions of the genus functions of $X$ and $Y$ respectively to $F_X$ and $F_Y$ are called  \textit{algebraically equivalent} if there exists an isomorphism $F_X\to F_Y$ that preserves the restrictions of the genus functions and the intersection forms. 
\end{definition}

It is clear that if the genus functions of two 4-manifolds are algebraically inequivalent, then they are orientedly inequivalent. Hence, if they are algebraically inequivalent for any orientations, then they are inequivalent. 
We remark that the definition of equivalence in the paper \cite{SS24} of Stipsicz and Szab\'{o} is the same as our algebraic equivalence, not our (ordinary) equivalence.

As we will show in Sections~\ref{sec:Genus functions and quasi-invertible cobordisms} and \ref{sec:genus:sum}, it is useful to consider the genus function modulo torsion part. We thus introduce the following function.  
\begin{definition}For a class $A$ of $H_2(X)/\mathrm{Tor}$, we define an integer $g^*_X(A)$ by
\begin{equation*}
	g^*_X(A)=\min\{g_X(\alpha)\mid \textnormal{$\alpha\in H_2(X)$ represents $A$.}\}.
\end{equation*}	
The map $g^*_X:H_2(X)/\mathrm{Tor}\to \mathbb{Z}$ is called the \textit{torsion-free genus function} of $X$. 
\end{definition}
Clearly, the torsion-free genus function is completely determined by the genus function. The torsion-free genus function is often easier to deal with, since the adjunction inequality for the Seiberg--Witten invariant gives a lower bound for $g_X(\alpha)$ up to torsion. By contrast, here we show that the values of the genus functions do depend on torsion parts even for closed 4-manifolds having non-vanishing Seiberg--Witten invariants, which is of independent interest. Hence, the genus functions are not necessarily determined by the torsion-free genus functions.  

\begin{theorem}\label{sec:genus:thm:torsion} There exists a closed $4$-manifold $X$ with $b_2^+>1$ having non-vanishing Seiberg-Witten invariant that admits a class $\alpha$ and a torsion class $\tau$ of $H_2(X)$  satisfying $g_{X}(\alpha)\neq g_{X}(\alpha+\tau)$.  Furthermore, there exist such closed $4$-manifolds having arbitrarily large $b_2^+$ and $b_2^-$.  
\end{theorem}
\begin{proof}Let $Y$ be an arbitrary closed spin 4-manifold, and let $R$ be the closed spin 4-manifold that is the total space of an $S^2$-bundle over $\mathbb{R}P^2$ (see \cite[Exercises 5.7.7.(a)]{GS}). We note that $R$ is a rational homology 4-sphere satisfying $H_2(R)\cong\mathbb{Z}/2\mathbb{Z}$.  
We consider the closed 4-manifold $X=Y\#2\overline{\mathbb{C}P^2}\#R$. In the case where $Y$ satisfies $b_2^+>1$ and has a non-vanishing Seiberg--Witten invariant, the 4-manifold $X$ also satisfies these conditions due to the blow-up formula \cite[Proposition~2]{KMW}. Since each elliptic surface $E(2n)$ ($n\in \mathbb{N}$) is an example of such a 4-manifold $Y$, the 4-manifold $X$ realizes arbitrarily large $b_2^+$ and $b_2^-$. 

Let $e_1,e_2$ be a standard orthogonal basis of $H_2(2\overline{\mathbb{C}P^2})$. Then, the mod 2 reduction of the class $\alpha=3e_1+e_2$ of $H_2(X)$ is the Poincar\'{e} dual of $w_2(X)\in H^2(X;\mathbb{Z}/2\mathbb{Z})$. Due to Rochlin's signature theorem, the signature $\sigma(X)$ of $X$ satisfies $\alpha\cdot \alpha-\sigma(X)\equiv -8\not\equiv 0\pmod{16}$. We thus obtain $g_X(\alpha)>0$ by \cite[Exercises 5.7.7.(b)]{GS}). (In fact, $g_X(\alpha)=1$ as seen from the rest of this proof.)  Let $V$ be the 4-manifold obtained from the 4-ball by attaching a 2-handle along the $(-10)$-framed figure-eight knot. This manifold ${V}$ admits an embedding into $2\overline{\mathbb{C}P^2}$ that sends a generator of $H_2(V)$ to the class $\alpha$ due to a result of Aceto et.al.\ \cite{ACMPS}. Furthermore, by a result of Levine \cite[Theorem~1.1 and Remark~1.2]{Le}, the figure-eight knot bounds a disk in the rational homology 4-ball $R-\mathrm{int}\, D^4$. The above inequality thus implies $g_{X}(\alpha+\tau)=g_{2\overline{\mathbb{C}P^2}\#R}(\alpha+\tau)=0< g_{X}(\alpha)$, where $\tau$ denotes the torsion class of $H_2(R)$. Therefore, the desired claim follows.   
\end{proof}

We define equivalences for torsion-free genus functions in the same way as before. We state two of them for accuracy. 
\begin{definition}(1) The torsion-free genus functions of $X$ and $Y$ are called \textit{equivalent} if there exists a homeomorphism $f:X\to Y$ such that the isomorphism $H_2(X)/\mathrm{Tor}\to H_2(Y)/\mathrm{Tor}$ induced by $f$ preserves the torsion-free genus functions.  
	
	(2) The torsion-free genus functions of $X$ and $Y$ are called \textit{algebraically equivalent} if there exists an isomorphism $H_2(X)/\mathrm{Tor}\to H_2(Y)/\mathrm{Tor}$ that preserves the torsion-free genus functions and the intersection forms. 
\end{definition}

Let us recall that a rational surface is orientedly diffeomorphic to $\mathbb{C}P^2\#n \overline{\mathbb{C}P^2}$, and that a ruled surface  is orientedly diffeomorphic to the total space of an $S^2$-bundle over a closed surface or to its blow-ups. 
We observe that genus functions uniquely determine smooth structures of rational surfaces and ruled surfaces among symplectic 4-manifolds. This fact easily follows from \cite[Proposition~4.3]{LMY}. 

\begin{proposition}[cf.\  {\cite{LMY}}]\label{prop:genus_symplectic} 
If the genus function of a given closed symplectic $4$-manifold is equivalent to that of a rational surface or a ruled surface, then the given  $4$-manifold is diffeomorphic to the surface. 
\end{proposition}
\begin{proof}Let $N_g$ denote the total space of a non-trivial $S^2$-bundle over  the closed surface $\Sigma_g$ of genus $g\geq 0$, and  let  $R$ be a rational surface or a ruled surface. Then $R$ is orientedly diffeomorphic to one of the $4$-manifolds $\mathbb{C}P^2\#n \overline{\mathbb{C}P^2}$, $(S^2\times \Sigma_g)\#n \overline{\mathbb{C}P^2}$ and $N_g\#n \overline{\mathbb{C}P^2}$ for some  non-negative integers $n,g$. Note that $R$ satisfies $b_2^+=1$. We easily see that $R$ contains a closed surface representing a non-torsion second homology class whose self-intersection number is larger than $2g'-2$, where  $g'$ denotes the genus of the surface. In the case where $b_2^+(\overline{R})\geq 1$, the orientation reversal $\overline{R}$ also contains such a surface as seen from its diffeomorphism type.  

Now suppose that the genus function of a closed symplectic 4-manifold $X$ is equivalent to that of $R$. Then we have $b_2^+(X)\geq 1$, and it thus follows from the above paragraph that $X$ contains a closed  surface of genus $g'$ representing a non-torsion second homology class whose self-intersection number is larger than $2g'-2$. Hence, by \cite[Proposition~4.3]{LMY}, $X$ is orientedly diffeomorphic to a rational surface or to a ruled surface. We see that two 4-manifolds in these classes of 4-manifolds are diffeomorphic if and only if they are homeomorphic, by comparing their homology groups and intersection forms. Since $R$ is also diffeomorphic to one of these 4-manifolds and is homeomorphic to $X$, it follows that $X$ is diffeomorphic to $R$.   
\end{proof}

The above proposition holds for torsion-free genus functions as well, because the second homology groups of rational surfaces and ruled surfaces are torsion-free.

\section{Corks, invertible cobordisms and exotic 4-manifolds}
In this section, we review (generalized) corks, invertible cobordisms and constructions of exotic 4-manifolds that utilize corks and invertible homology cobordisms. Furthermore, we improve the method of Akbulut and Ruberman~\cite{AR} that modifies  generalized corks into exotic 4-manifolds. 

\subsection{Generalized corks}
In this paper, we define generalized corks as follows (cf.\ \cite{AY08}, \cite{T15}). 

\begin{definition}(1) Let $C$ be a (not necessarily contractible) 4-manifold and let $f:\partial C\to \partial C$ be a self-diffeomorphism of the boundary $\partial C$ that preserves the orientation. The pair $(C,f)$ is called a \textit{generalized cork} if $f$ extends to a self-homeomorphism of $C$, but does not extend to any self-diffeomorphism of $C$. The \textit{order} of a generalized cork $(C,f)$ is defined to be the minimum positive integer $k$ for which the power $f^k$ extends to a self-diffeomorphism of $C$. If  $f^k$ does not extend to any self-diffeomorphism of $C$ for any positive integer $k$, then $(C,f)$ is called an \textit{infinite order generalized cork}. A generalized cork $(C,f)$ is called a \textit{cork} if $C$ is contractible. 

(2) Let $(C,f)$ be a generalized cork, and assume $C$ is embedded in a 4-manifold $X$. We remove $C$ from $X$ and glue it back along the resulting boundary via the gluing map $f$. This operation is called a \textit{generalized cork twist} along $(C,f)$. 

(3) Two generalized corks $(C,f)$ and $(C,g)$ are called \textit{equivalent} if the self-diffeomorphism $g^{-1}\circ f$ of $\partial C$ extends to a self-diffeomorphism of $C$. 
\end{definition}

Clearly, a generalized cork twist always preserves the homeomorphism type of an ambient 4-manifold. We note that in the case where $C$ is contractible, every self-diffeomorphism of $\partial C$ extends to a self-homeomorphism of $C$ (\cite{Fr}, \cite{Bo86}). 
\subsection{Invertible cobordisms}
We recall the definitions of an invertible cobordism and related terminologies and fix conventions. Throughout this paper, we denote the interval $[0,1]$ by $I$.  Let $M,N$ be closed (possibly disconnected) 3-manifolds, and let $X$ be a 4-manifold having $M$ as an oriented boundary component. For convenience, we allow $N$ to be empty. 

\begin{definition}
(1) A possibly disconnected 4-manifold $P$ is called a \textit{$($$4$-dimensional$)$ cobordism} from $M$ to $N$ if $\partial P$ is orientedly diffeomorphic to the disjoint union $\overline{M}\sqcup N$, and every connected component of $P$ has a non-empty boundary component that belongs to $\overline{M}$. In this case,  $M$ and $N$ are respectively called a negative boundary and a positive boundary of $P$.  We often identify $\partial P$ with the disjoint union $\overline{M}\sqcup N$ by fixing a diffeomorphism. For an orientation-preserving self-diffeomorphism $f$ of $M$, let $X\cup_f P$ denote the 4-manifold constructed from $X$ by attaching $P$ along $M$, where we glue each $x\in M\subset \partial X$ with $f(x)\in M\subset \overline{\partial P}$. When we do not specify the gluing map $f$, we often denote a glued 4-manifold by $X\cup_M P$. Unless otherwise stated, we always attach a cobordism to a 4-manifold along a negative boundary of the cobordism.   

(2) For a cobordism $P$ from $M$ to $N$, the cobordism $\overline{P}$ from ${N}$ to ${M}$ is called the \textit{upside down cobordism} of $P$.  

(3) For a cobordism $P$ from $M$ to $N$, if the homomorphisms $H_*({M})\to H_*(P)$ and $H_*(N)\to H_*(P)$ induced by the inclusions are isomorphisms, then $P$ is called a  \textit{homology cobordism}. If there exists a homology cobordism from $M$ to $N$, then $M$ and $N$ are called \textit{homology cobordant}.  

(4) A cobordism $P$ from $M$ to $N$ is called \textit{invertible} if there exists a cobordism $Q$ from $N$ to $M$ such that a glued 4-manifold $P\cup_N Q$  is orientedly diffeomorphic to the product $I\times M$. The cobordism $Q$ is called an \textit{inverse cobordism} of $P$. 
\end{definition}

We use the following terminology. 
\begin{definition}For a positive integer $n$, attach $n$ 1-handles $h_1^1, h_2^1, \dots, h_n^1$ to a 4-manifold $X$, and then attach $n$ 2-handles to the resulting 4-manifold $X'$. The set of these 1- and 2-handles is called \textit{homotopically canceling pairs of $1$- and $2$-handles} if the link in $\partial X'$ consisting of the attaching circles of 2-handles is homotopic to a link consisting of knots $K_1, K_2, \dots, K_n$ in $\partial X'$, such that each $K_i$ intersects the belt sphere of the 1-handle $h_i^1$ geometrically exactly once.   
\end{definition}
Here we give a simple construction of an invertible cobordism as well as its properties. 
We will later give much more flexible constructions, see Lemma~\ref{sec:quasi-invertible:lem:quasi-invertible:generate}. 

\begin{lemma}\label{sec:cork:thm:inverse cobordism}Suppose a  cobordism $P$ from $M$ to $N$ is obtained from $I\times M$ by attaching homotopically canceling pairs of $1$- and $2$-handles, none of which are attached to $0\times M$. Then $P$ is invertible, and the upside down cobordism $\overline{P}$ is its inverse cobordism. Moreover, $P\cup_{\mathrm{id}_N} \overline{P}$ is diffeomorphic to $I\times M$ fixing the boundary $\overline{M}\sqcup M$ pointwise, where $\overline{M}\sqcup M$ is identified with $\{0,1\}\times M$ via the natural identification.  
\end{lemma}
\begin{proof}We consider the 5-manifold $P\times I$. This manifold is obtained from $(I\times M)\times I$ by attaching homotopically canceling pairs of $1$- and $2$-handles in $1\times M\times I$. Since a homotopy of  attaching circles of these 2-handles can be approximated by a smooth isotopy, we can eliminate these pairs.  We thus have a diffeomorphism between $P\times I$ and $(I\times M)\times I$ that fixes $(0\times M)\times I$ pointwise. 
	By restricting this diffeomorphism to the exteriors of $(0\times M)\times I$ in the boundaries of these two 5-manifolds, we see that $P\cup_{\mathrm{id}_N} \overline{P}$ is diffeomorhic to $I\times M$ fixing the boundary $\overline{M}\sqcup M$ pointwise.  
\end{proof}
\begin{remark}\label{sec:cork:rem:inverse cobordism}Under the assumption of this lemma, $P$ is a homology cobordism. Furthermore, the homotopy extension property (see \cite[Proposition~0.18]{Hat}) implies that $P$ admits a deformation retraction onto the boundary component $M$. 
\end{remark}

\subsection{Modifying 2-handlebodies into exotic Stein manifolds}\label{subsec:exotic Stein}
In this subsection, we review the method of Akbulut and the author~\cite{AY13} that modifies a given 4-dimensional 2-handlebody with  non-trivial second homology group into arbitrarily many pairwise exotic  4-manifolds admitting Stein Structures. 

We say that a 4-dimensional \textit{$2$-handlebody} is a \textit{Legendrian handlebody} if the attaching circles of 2-handles form a Legendrian link in the standard tight contact structure on the boundary $k(S^1\times S^2)$ of the sub 1-handlebody $\natural_k S^1\times D^3$. We note that every 4-dimensional 2-handlebody becomes a Legendrian handlebody by isotoping the attaching circles of 2-handles. We say that a 4-dimensional 2-handlebody is a \textit{Stein handlebody} if it is Legendrian handlebody such that the framing of the attaching circle of each 2-handle is  one less than the Thurston--Bennequin number. By results of Eliashberg~\cite{E90} and Gompf~\cite{G98}, a 4-manifold admits a Stein structure if and only if it is orientedly diffeomorphic to a Stein handlebody. 

\begin{figure}[ht!]
\begin{center}
\includegraphics[width=0.8\columnwidth]{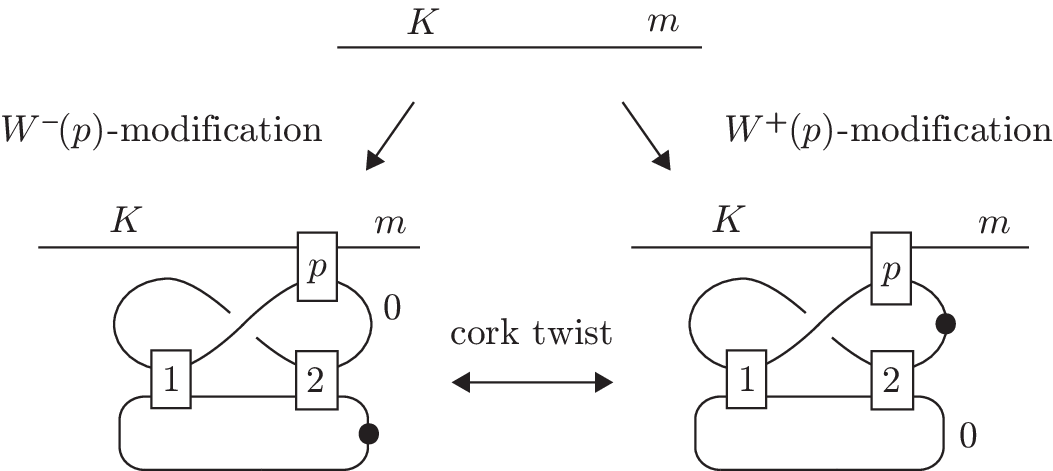}
\caption{$W^-(p)$- and $W^+(p)$-modifications}
\label{fig:W_modifications}
\end{center}
\end{figure}

Let $X$ be a 4-dimensional 2-handlebody, and let $K$ be the attaching circle of a 2-handle of $X$. We consider a handlebody diagram of $X$. Akbulut and the author~\cite{AY13} introduced two local diagrammatic operations for a subarc of $K$ shown in Figure~\ref{fig:W_modifications}. For a positive integer $p$, the left side operation is called a \textit{$W^-(p)$-modification} to $K$, and the right side operation is called a \textit{$W^+(p)$-modification} to $K$. We note that these operations do not change the framing coefficient $m$ of $K$. When we do not specify the coefficient $p$, these operations are called a $W^-$-modification and a $W^+$-modification, respectively. 

Each operation creates an embedded copy of the Akbulut cork $W$ (\cite{A91a}) in the resulting handlebody, and these two operations are related to each other by the cork twist along the embedded $W$ as shown in the figure. So, the 4-manifold obtained by the $W^+(p)$-modification to $K$ is orientedly homeomorphic to the 4-manifold obtained by the $W^-(p)$-modification to $K$.

Every $W^-$-modification attaches a homotopically canceling pair of 1- and 2-handles to $X$. Thus, this operation is an attachment of an invertible homology cobordism to $X$, and its upside down cobordism is its inverse cobordism (see Lemma~\ref{sec:cork:thm:inverse cobordism}). Furthermore, this invertible cobordism has a deformation retraction onto $\partial X$ (see Remark~\ref{sec:cork:rem:inverse cobordism}), and hence the result of any $W^-$-modification admits a deformation retraction onto $X$. Consequently, the result of any $W^-$-modification and hence that of any $W^+$-modification is HIHC-equivalent to  $X$. In addition, the 4-manifold obtained by any $W^+$-modification becomes diffeomorphic to the original 4-manifold $X$ by attaching the upside down cobordism of the invertible cobordism given by the corresponding $W^-$-modification (\cite[Proposition 4.5]{AY13}). 

By drawing Legendrian pictures, we obtain Legendrian versions of $W^-$- and $W^+$-modifications for Legendrian handlebodies, and these have the following two effects (\cite[Proposition 4.7]{AY13}). (1) Each of these two operations modifies a Stein handlebody into another Stein handlebody. (2) If we apply a $W^+$-modification to a Legendrian knot $K$ of a Legendrian handlebody, then the Thurston--Bennequin number of the result of $K$ increases by $p$, and the Thurston--Bennequin number for the attaching circle of the 2-handle of $W$ is $+2$. 

Thus, by applying $W^+$-modifications to a Legendrian handlebody and adding zig-zags to each Legendrian knot (i.e., stabilizations), we obtain a Stein handlebody, implying the following.  

\begin{proposition}[\cite{AY13}]\label{sec:cork:prop:2-handlebody:Stein}
	For any $4$-dimensional $2$-handlebody $X$, there exists a $4$-dimensional $2$-handlebody $X'$ having the following properties. 
\begin{itemize}
 \item $X'$ admits a Stein structure. 
 \item $X'$ is HIHC-equivalent to  $X$ and is homeomorphic to a $4$-manifold that deformation retracts onto an embedded copy of $X$. 
 \item $X'$ becomes diffeomorphic to $X$ after attaching the upside down cobordism of an invertible homology cobordism along $\partial X'
 $.  
\end{itemize}
\end{proposition}

Applying $W^-$-modifications several times and replacing some of them with $W^+$-operations by cork twists, we obtain arbitrarily many pairwise homeomorphic 2-handlebodies. In the case where $H_2(X)\neq 0$, by carefully choosing coefficients of $W^+$- and $W^-$-modifications and applying adjunction inequalities, Akbulut and the author showed that the resulting 2-handlebodies have pairwise  algebraically inequivalent genus functions for any orientations (\cite[Lemma 5.13 and Proposition 6.2]{AY13}), giving the following result. 

\begin{theorem}[\cite{AY13}]\label{sec:cork:exotic Stein}For any $4$-dimensional $2$-handlebody $X$ with $H_2(X)\neq 0$ and any positive integer $n$, there exist pairwise exotic $4$-manifolds $X_0, X_1,\dots,X_n$ having the following properties. 
\begin{itemize}
 \item Each $X_i$ is HIHC-equivalent to $X$ and admits an embedding into $X$.  
 \item $X_0$ is obtained from $X$ by attaching an invertible homology cobordism, and admits a deformation retraction onto $X$.  
 \item $X_i$ admits a Stein structure for any $i\neq 0$.  
 \item The genus functions of $X_0, X_1,\dots,X_n$ are pairwise algebraically inequivalent for any orientations. 
\end{itemize}
\end{theorem}

The proof of this theorem yields a concrete handlebody diagram of each $X_i$ from that of $X$ and hence produces various explicit examples of exotic 4-manifolds.   
\subsection{Modifying generalized corks into exotic 4-manifolds}In this subsection, we review the method of Akbulut and Ruberman~\cite{AR} that modifies generalized corks into exotic 4-manifolds. We also improve their results. 

Throughout this subsection, we fix an arbitrary closed 3-manifold $M$. Note that $M$ is connected due to our convention. 
Akbulut and Ruberman~\cite{AR} gave a particularly useful invertible homology cobordism.  
\begin{theorem}[\cite{AR}]\label{thm:AR:cobordism}
There exist a closed $3$-manifod $N$ and an invertible homology cobordism $\mathbb{P}$ from $M$ to $N$ that satisfy the following. 
\begin{itemize}
 \item Every self-diffeomorphism of $N$ extends to a self-diffeomorphism of $\mathbb{P}$ whose restriction to $M$ is isotopic to the identity map.  
 \item The homomorphism $\pi_1(M)\to \pi_1(\mathbb{P})$ induced by the inclusion is an isomorphism. 
 \end{itemize}
\end{theorem}

Utilizing this cobordism, Akbulut and Ruberman~\cite{AR}  constructed an exotic ``pair'' of 4-manifolds from a given generalized cork by attaching the above cobordism and then twisting the embedded generalized cork. In order to produce infinitely many pairwise exotic 4-manifolds, they constructed another invertible homology cobordism having properties different from the above cobordism. They  constructed infinitely many 4-manifolds from a given infinite order cork by attaching the latter cobordism and then twisting the embedded cork, and showed that at least infinitely many of them are pairwise exotic using the pigeonhole principle. However, for a non-simply connected generalized cork of infinite order, it is unclear whether exotic 4-manifolds given by the latter construction have the same fundamental groups as that of the generalized cork, since the latter cobordism has a weaker property on the fundamental groups.

Here we observe that, in the former construction, all 4-manifolds constructed by twisting a given generalized cork are pairwise exotic, by   modifying their proof. We prove this claim in a more general form.  

\begin{theorem}\label{thm:AR:arbitrary}
Let $f_1,f_2$ be self-diffeomorphisms of $M$, and let $X_1,X_2$ be $4$-manifolds with $\partial X_1=\partial X_2=M$. If $f_2^{-1}\circ f_1$ does not extend to any diffeomorphism $X_1\to X_2$, then the following hold. 
\begin{itemize}
 \item  The $4$-manifolds $X_1\cup_{f_1} \mathbb{P}$ and $X_2\cup_{f_2} \mathbb{P}$ are not diffeomorphic to each other. 
 \item Each $X_i\cup_{f_i} \mathbb{P}$  is  HIHC-equivalent to  $X_i$ and admits a deformation retraction onto $X_i$.   
\end{itemize}
\end{theorem}
\begin{remark}\label{rem:AR-cobordism used in our theorem}
	Akbulut and Ruberman constructed a cobordism $\mathbb{P}$ in Theorem~\ref{thm:AR:cobordism} using an arbitrary doubly slice knot in $S^3$ having certain properties, and they showed that the Kinoshita--Terasaka knot $J$ (see \cite[Figure~2]{AR}) can be used for the construction. In Theorem~\ref{thm:AR:arbitrary}, the cobordism $\mathbb{P}$ denotes the one corresponding to the knot $J$ in order to guarantee the homotopy equivalence. The rest of the properties hold for any  $\mathbb{P}$ in Theorem~\ref{thm:AR:cobordism} as seen from our proof. 
\end{remark}

\begin{proof}[Proof of Theorem~$\ref{thm:AR:arbitrary}$]
Let $\mathbb{Q}$ be an inverse cobordism of $\mathbb{P}$. Then the 4-manifold $\mathbb{P}\cup_{q} \mathbb{Q}$ is diffeomorphic to the product $I\times M$ for some self-diffeomorphism $q$ of $N$. 
For convenience, we set $\partial \mathbb{P}=\overline{M_{\mathbb{P}}}\sqcup N_{\mathbb{P}}$ and $\partial \mathbb{Q}=\overline{N_{\mathbb{Q}}}\sqcup M_{\mathbb{Q}}$, where $M_{\mathbb{P}}, M_{\mathbb{Q}}$ are identified with $M$, and $N_{\mathbb{P}}, N_{\mathbb{Q}}$ are identified with $N$.  
We may assume that the restriction $M_{\mathbb{P}}\to M(=0\times M)$ of the diffeomorphism $\mathbb{P}\cup_q \mathbb{Q}\to I\times M$ is the identity map if necessary by composing a self-diffeomorphism of $I\times M$. Let $\varphi:M_{\mathbb{Q}}\to M(=1\times M)$ denote the restriction of this diffeomorphism $\mathbb{P}\cup_q \mathbb{Q}\to I\times M$.  

We set $X'_i=X\cup_{f_i} \mathbb{P}$ for each $i$. Then each $X'_i$ deformation retracts  to $X$ by the forthcoming Lemma~\ref{lem:AR:homotopy}. Since $\mathbb{P}$ is a homology cobordism, it follows that $X'_i$ is HIHC-equivalent to  $X_i$. 

We will show that $X'_1$ is not diffeomorphic to $X'_2$. Suppose to the contrary that there exists a diffeomorphism $X'_1\to X'_2$. Since every self-diffeomorphism of $N_{\mathbb{P}}$ extends to a self-diffeomorphism of $\mathbb{P}$ whose restriction to $M_{\mathbb{P}}$ is isotopic to the identity, we may assume that the diffeomorphism $X'_1\to X'_2$ is the identity on $N_{\mathbb{P}}$ if necessary by composing a self-diffeomorphism of $X'_2$. We thus have a diffeomorphism $\Phi:X'_1\cup_q \mathbb{Q}\to X_2'\cup_q \mathbb{Q}$ whose restriction to $M_{\mathbb{Q}}$ is the identity map. 

For each $i$, we also use a diffeomorphism $X'_i\cup_q \mathbb{Q}\to X_i\cup_{\mathrm{id}_M}(I\times M)$ whose restriction to the boundary is $f_i^{-1}\circ \varphi:M_{\mathbb{Q}}\to M$. 
This diffeomorphism is obtained as follows. We have a diffeomorphism $X'_i\cup_{q} \mathbb{Q}\to X_i\cup_{f_i}(I\times M)$ whose restriction to the boundary is the diffeomorphism $\varphi$. Since $f_i^{-1}$ extends to a self-diffeomorphism of $I\times M$, we have a diffeomorphism $X_i\cup_{f_i}(I\times M)\to X_i\cup_{\mathrm{id}_M}(I\times M)$ whose restriction to the boundary $1\times M$ is $f_i^{-1}$. By composing these two diffeomorphisms, we obtain a desired diffeomorphism $X'_i\cup_{q} \mathbb{Q}\to X_i\cup_{\mathrm{id}_M}(I\times M)$. 

Composing the last diffeomorphisms for $i=1,2$ and $\Phi$, we obtain a diffeomorphism $X_1\cup_{\mathrm{id}_M}(I\times M)\to X_2\cup_{\mathrm{id}_M}(I\times M)$ whose restriction to the boundary is $f_2^{-1}\circ f_1$. 
Since each  $X_i\cup_{\mathrm{id}_M}(I\times M)$ is diffeomorphic to $X_i$ fixing the boundary $M$ pointwise, it follows that $f_2^{-1}\circ f_1$ extends to a diffeomorphism $X_1\to X_2$, giving a contradiction. Therefore, $X'_1$ is not diffeomorphic to $X'_2$. 
\end{proof}

The homotopy equivalence of  $X_i\cup_{f_i} \mathbb{P}$ and $X_i$ is claimed in \cite[Theorem~A]{AR}, but the proof is given only in the case where $X_i$ is contractible, and the argument does not seem to extend to the general case. For completeness, we prove the general case  for $\mathbb{P}$ used in Theorem~\ref{thm:AR:arbitrary} (see Remark~\ref{rem:AR-cobordism used in our theorem}) by showing the lemma below. 
\begin{lemma}\label{lem:AR:homotopy}
	The cobordism $\mathbb{P}$  used in Theorem~$\ref{thm:AR:arbitrary}$ admits a deformation retraction onto the boundary component $M$.  
\end{lemma}
\begin{proof}
	In order to construct $\mathbb{P}$,  Akbulut and Ruberman~\cite{AR} used a certain link $L$ in $M$ and a certain concordance $C(\subset I\times S^3)$ from the unknot $U$ to the Kinoshita--Terasaka knot  $J$. Let $\nu(L)$ and $\nu(C)$ respectively denote closed tubular neighborhoods of $L$ and $C$, and let $\mathring{\nu}(L)$ and $\mathring{\nu}(C)$ denote the open tubular neighborhoods. 
	Let $n$ be the number of the components of $L$. 
	
	They constructed $\mathbb{P}$ from $I\times (M-\mathring{\nu}(L))$ by attaching  $n$ copies of $(I\times S^3)-\mathring{\nu}(C)$ along $I\times \partial \nu(L)$ so that the resulting negative boundary  remains $M$. Note that there is a diffeomorphism from $[-1,0]\times M \cup_{0\times  (M-\mathring{\nu}(L))} I\times (M-\mathring{\nu}(L))$  to $I\times M$ that maps $(0\times \nu(L))\cup (I\times \partial \nu(L))$ to $1\times \nu(L)$. Thus,  $\mathbb{P}$ is obtained from $I\times M$ by attaching  $n$ copies of $(I\times S^3)-\mathring{\nu}(C)$ using a diffeomorphism from each component of $1\times \nu(L)$ to $\nu(m)$, where $m$ denotes the meridian of the unknot  $U$ in $0\times S^3$, and $\nu(m)$ denotes a closed tubular neighborhood of $m$ in $\partial (I\times S^3-\mathring{\nu}(C))$. 
	
	Akbulut and Ruberman gave a handlebody diagram  of the exterior of a slice disk for $J$ in $D^4$ (\cite[Figure~4.(a)]{AR}), and it is easy to see that this 4-manifold is diffeomorphic to  $(I\times S^3)-\mathring{\nu}(C)$. Note that $m$ is isotopic to the meridian of a dotted circle in the diagram. Since the unique 2-handle of this handlebody and the 1-handle corresponding to the other dotted circle  form a homotopically canceling pair,  Remark~\ref{sec:cork:rem:inverse cobordism} implies that $(I\times S^3)-\mathring{\nu}(C)$ deformation retracts to the tubular neighborhood $\nu(m)$. Therefore, $\mathbb{P}$ admits a deformation retraction onto the boundary component $0\times M$. 
\end{proof}
\begin{remark}\label{rem:AR:cobordism}(1) An inverse cobordism of the cobordism $\mathbb{P}$ used in Theorem~\ref{thm:AR:arbitrary} (see Remark~\ref{rem:AR-cobordism used in our theorem}) is its upside down cobordism $\overline{\mathbb{P}}$, though we do not use this fact in this section. This claim is not stated in \cite{AR},  but it is implicit and follows from their claim that the Kinoshita--Terasaka knot $J$ is superslice, that is,  the double of a slice disk for $J$ is an unknotted 2-sphere in $S^4$ (see the proof of \cite[Proposition~2.6]{AR}). The supersliceness can be alternatively proved using \cite[Proposition~4.3]{LiMe15}, and it is clear from the latter proof that the slice disk obtained in \cite{AR} gives a concodance in $I\times S^3$ from the unknot to $J$ such that the double is isotopic to the product concordance fixing $\{0,1\}\times S^3$ pointwise. Therefore, it follows from their construction of  $\mathbb{P}$ that the double $\mathbb{P}\cup_{\mathrm{id}_M} \overline{\mathbb{P}}$ admits a diffeomorphism to the product $I\times M$ whose restriction to the boundary $\overline{M}\sqcup M$ is the identity map.   
	
(2) The cobordism $\mathbb{P}$ is obtained from $I\times M$ by attaching homotopically canceling pairs of 1- and 2-handles, and Lemma~\ref{lem:AR:homotopy} and the claims in the above (1) follow from this fact and Lemma~\ref{sec:cork:thm:inverse cobordism} as well. This fact is implicit from the handlebody diagram of an example of a contractible 4-manifold in \cite[Example~4.1]{AR} and can be proved using  Akbulut's cylinder method \cite[Section~3.2]{A_book} as follows. The proof of Lemma~\ref{lem:AR:homotopy} implies that the cobordism $\mathbb{P}$ is obtained from the disjoint union of $I\times M$ and copies of the exterior of the slice disk for $J$ by attaching copies of the cylinder $I\times S^1\times D^2$ along $\partial I\times S^1\times D^2$.  Since this cylinder is a union of a 1-handle and a 2-handle, the cobordism $\mathbb{P}$ is obtained from the boundary sum of $I\times M$ with the copies of  the exterior by attaching 2-handles. Due to the handlebody diagram of the exterior (\cite[Figure~4.(a)]{AR}), we can prove the aforementioned fact by canceling a 1-handle of each copy of the exterior using the 2-handle of the cylinder. We remark that this argument gives a method for drawing handlebody diagrams of $\mathbb{P}$ and each glued 4-manifold  $X_i\cup_{f_i} \mathbb{P}$, once the link $L$ in $M$ is specified.  
\end{remark}

We obtain the following corollary, improving a result of Akbulut and Ruberman. 

\begin{corollary}\label{cor:AR:arbitrary}
	For any family $\{(C,f)\}_{f\in A}$ of pairwise inequivalent generalized corks,  there exists a family $\{C_{f}\}_{f\in A\sqcup \{\mathrm{id}_{\partial C}\}}$ of pairwise exotic $4$-manifolds which are HIHC-equivalent to  $C$. 
\end{corollary}
\begin{proof}We set $M=\partial C$, and let $\mathbb{P}$ be the cobordism used in Theorem~\ref{thm:AR:arbitrary}. For each $f\in A\sqcup \{\mathrm{id}_{\partial C}\}$, let $C_f$ denote the 4-manifold $C\cup_f \mathbb{P}$. Then, it follows from Theorem~\ref{thm:AR:arbitrary} that $\{C_{f}\}_{f\in A\sqcup \{\mathrm{id}_{\partial C}\}}$ is a desired family. 
\end{proof}	
We note that for any generalized cork $(C,f)$ of infinite order, $\{(C,f^k)\}_{k\in \mathbb{Z}-\{0\}}$ satisfies the assumption of this corollary. The assumption holds for families given by finite order generalized corks as well.  
\subsection{Infinite order corks}\label{subsec:Gompf's corks}
\begin{figure}[ht!]
	\begin{center}
		\includegraphics[width=0.32\columnwidth]{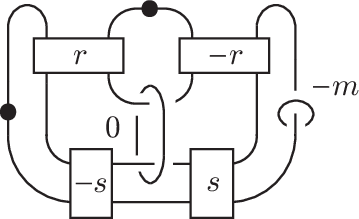}
		\caption{$C(r,s;m)$}
		\label{fig:infinite_order_cork}
	\end{center}
\end{figure}

In this subsection,  we briefly review Gompf's infinite order corks (\cite{G17GT}). For integers $r,s,m>0$, let $C(r,s;m)$ be the 4-manifold in Figure~\ref{fig:infinite_order_cork}. (In this notation, we changed the sign of $m$ from the original one in \cite{G17GT} so that we may assume $m>0$.) It is easy to see that $C(r,s;m)$ is a contractible 4-manifold of \textit{Mazur type}, that is, it admits a handle decomposition consisting of a 0-handle and an algebraically canceling pair of 1- and 2-handles (i.e., a 1-handle and a 2-handle whose attaching circle intersects the belt sphere of the 1-handle algebraically once). Note that such a pair attached to the 0-handle is necessarily a homotopically canceling pair. The boundary $\partial C(r,s;m)$ gives infinitely many pairwise non-homeomorphic homology 3-spheres by varying $r,s,m$ (\cite[p.2477]{G17GT}). 

Gompf proved that each pair $(C(r,s;m),\tau)$ is an infinite order cork, where the self-diffeomorphism $\tau$ on the boundary is the torus twist on a certain embedded torus. Roughly speaking, a torus twist is $(\textnormal{Dehn twist on an annulus})\times \mathrm{id}_{S^1}$, and hence $\tau$ is the identity map outside the tubular neighborhood of the torus. For the precise definition of a torus twist, see \cite[Definition~2.1]{G17AGT}. 

\begin{figure}[ht!]
		\begin{center}
	\begin{minipage}[b]{0.48\linewidth}
		\begin{center}
			\includegraphics[width=0.52\columnwidth]{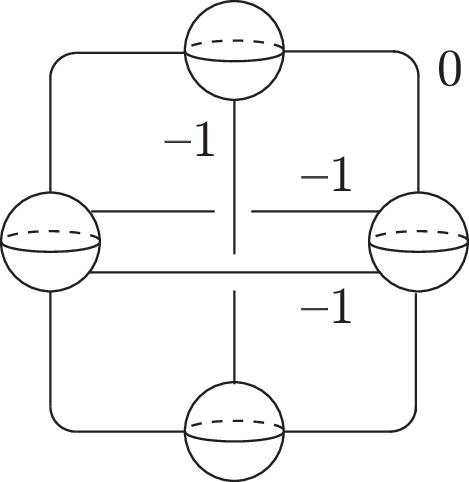}
			\caption{$E$}
			\label{fig:cork_nbd}
		\end{center}
	\end{minipage}%
	\begin{minipage}[b]{0.48\linewidth}
		\begin{center}
			\includegraphics[width=0.38\columnwidth]{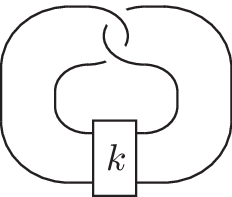}
			\caption{Twist knot $T_k$}
			\label{fig:twist_knot}
		\end{center}
		\end{minipage}
			\end{center}
\end{figure}

We recall an important property of Gompf's corks. Let $E$ be the 4-manifold given in Figure~\ref{fig:cork_nbd}. We denote the obvious torus $T^2\times 0 \subset T^2\times D^2$ in the figure by $T$. Also, for an integer $k$, we denote the twist knot in Figure~\ref{fig:twist_knot} by $T_k$. Let $E_{k}$ be the 4-manifold obtained from $E$ by applying the Fintushel--Stern knot surgery (\cite{FS98}) along the torus $T$ using the knot $T_k$. 
Gompf~\cite{G17GT} proved the following theorem to show that each $(C(r,s;m),\tau)$ is an infinite order cork. For details, see \cite[the proof of Theorem~3.5]{G17AGT}. 

\begin{theorem}[\cite{G17GT}, \cite{G17AGT}]\label{thm:cork_embedding}Let $r,s,m$ be positive integers, and set $N=r+s+m-3$. Then there exists an embedding of $C(r,s;m)$ into $E\#N\overline{\mathbb{C}P^2}$ such that for each integer $k$, the knot surgered $4$-manifold $E_{k}\#N\overline{\mathbb{C}P^2}$ is diffeomorphic to the $4$-manifold obtained from $E\#N\overline{\mathbb{C}P^2}$ by removing the embedded $C(r,s;m)$ and gluing it back via $\tau^k$. Furthermore, this diffeomorphism is isotopic to the identity map on the boundary $\partial E$ of these two $4$-manifolds. 
\end{theorem}
We note that, in the last claim of this theorem, we identified the boundaries of these two 4-manifolds with $\partial E$, since the knot surgery and the cork twist are applied in the interior of $E\#N\overline{\mathbb{C}P^2}$. 

\section{Producing infinite order generalized corks}
In this section, we show that infinite order generalized corks exist for a large class of homeomorphism types, which is of independent interest. 

 The first example of an infinite order generalized cork was given by Tange~\cite[Theorem~1.6]{T15}, and his example is simply connected and satisfies $b_2=2$. Gompf~\cite{G17GT} later gave the first examples of infinite order contractible corks. We show that every 2-handlebody can be modified into an infinite order generalized cork, showing that infinite order generalized corks exist in abundance. 

\begin{theorem}\label{sec:cork:thm:2-handlebody:cork}For any $4$-dimensional $2$-handlebody $X$, there exists a generalized cork $(X',f)$ of infinite order such that $X'$ is HIHC-equivalent to  $X$. Furthermore, in the case where $H_2(X)\neq 0$, there exist arbitrarily many  generalized corks of infinite order such that the $4$-manifolds constituting these corks are pairwise exotic and are HIHC-equivalent to $X$. 
\end{theorem}

We will prove this theorem together with further properties in Section~\ref{sec:exotic}, see Theorem~\ref{sec:exotic:thm:2-handlebody:cork}. 

Let us recall that a symplectic 4-manifold $(Z,\omega)$ is said to have weakly convex boundary if  $\partial Z$ admits a positive contact structure $\xi$ satisfying $\omega|_\xi>0$, that is, $\omega(v,w)>0$ for any oriented basis $v,w$ for $\xi$ (see \cite{Et09}). Here we prove the following proposition to prove the above theorem and main results of this paper. 
\begin{proposition}\label{sec:cork:prop:construction}
	Suppose  that a $4$-manifold $X$ admits an embedding into a symplectic $4$-manifold with weakly convex boundary. 
	Then, there exists an infinite order generalized cork $(X',f)$ such that $X'$ is obtained from $X$ by attaching an invertible homology cobordism and admits a deformation retraction onto $X$. Furthermore, $f$ extends to a self-homeomorphism of $X'$ that induces the identity map on $H_2(X')$. 
\end{proposition}

\begin{figure}[ht!]
\begin{center}
\includegraphics[width=0.32\columnwidth]{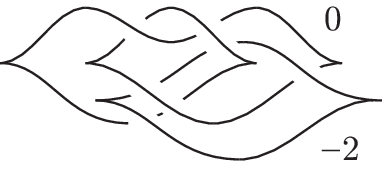}
\caption{Stein handlebody diagram of $E$}
\label{fig:Stein_cork_nbd}
\end{center}
\end{figure}

Let $E$ be the 4-manifold given by Figure~\ref{fig:cork_nbd}. We show the lemma below. 
\begin{lemma}$E$ admits a Stein structure. 
\end{lemma}
\begin{proof}By canceling the 1-handles and isotoping the resulting 2-handles, we obtain the Stein handlebody diagram in Figure~\ref{fig:Stein_cork_nbd}, where the coefficients are framings relative to the Seifert framings. Therefore, $E$ admits a Stein structure. 
\end{proof}

\begin{proof}[Proof of Proposition~$\ref{sec:cork:prop:construction}$]We fix positive integers $r,s,m$ and abbreviate the contractible 4-manifold $C(r,s;m)$ in Subsection~\ref{subsec:Gompf's corks} by $C$, and we set $N=r+s+m-3$. Since $C$ is of Mazur type, its double 
 $C\cup \overline{C}$ is diffeomorphic to $S^4$. 
 
Let $H$ be a homology cobordism obtained from $I\times \partial X$ by taking a boundary sum with $C$, where  the sum is taken outside the tubular neighborhood of the torus on which Gompf's torus twist $\tau$ is applied. 
We see that the double $H\cup \overline{H}$ is diffeomorphic to $(I\times \partial X)\# (C\cup \overline{C})$ and hence to $I\times \partial X$.  Thus,  $H$ is an invertible homology cobordism from $\partial X$ to $\partial X\# \partial C$. 

Let $X'$ be the 4-manifold obtained from $X$ by attaching the invertible homology cobordism $H$. Then, $X'$ deformation retracts to $X$ and is diffeomorphic to $X\natural C$. Let $f$ be the self-diffeomorphism of $\partial X'$ that extends the torus twist $\tau$ of $\partial C$ as the identity on the rest of $\partial X'$. Since $\tau$ extends to a self-homeomorphism of $C$, the self-diffeomorphism $f$ extends to a self-homeomorphism of $X'$ that is identity on $X$. We thus see that this self-homeomorphism induces the identity map on $H_2(X')$. 

We will show that $(X',f)$ is an infinite order generalized cork. By the assumption, $X$ admits an embedding into a symplectic 4-manifold  $Y$ with weakly convex boundary. 
Let $E'$ be the 4-manifold obtained from $E$ by attaching a 2-handle along a meridian of the 0-framed knot with sufficiently large negative framing. Then $Y\natural E'$  is obtained from $Y$ by attaching 2-handles along a Legendrian link with the contact $(-1)$-framings, and thus $Y\natural E'$ and its blow-ups $Y\natural E'\#N\overline{\mathbb{C}P^2}$ are symplectic 4-manifolds with weakly convex boundary (\cite[Theorem 2.5]{EH02}, \cite{Wei}, see also \cite[Theorem~5.8]{Et09}). We embed $X'$ into $Y\natural E'\#N\overline{\mathbb{C}P^2}$ by taking a boundary sum of $X\subset Y$ and $C\subset E'$ using a thickened path from $\partial X$ to $\partial C$. 

By a result of Eliashberg~\cite{El04} and Etnyre~\cite{Et04}, we can embed $Y\natural E'$ into a closed symplectic 4-manifold with $b_2^+>1$. We note that the torus $T$ in $E$ represents a second homology class of $Z$ having infinite order, since $T$ algebraically intersects the 2-sphere coming from the above 2-handle of $E'$. Let $Z_k$ denote the 4-manifold obtained from $Z$ by applying the knot surgery along $T$ using the twist knot $T_k$. It follows from the blow-up formula (\cite{FS95}) and the knot surgery formula (\cite{FS98}, \cite{F06}, see also \cite{Su}) that $Z_{k}\#N\overline{\mathbb{C}P^2}$ yields pairwise non-diffeomorphic 4-manifolds by varying $k$. On the other hand, Theorem~\ref{thm:cork_embedding} and the definition of $f$ show that each $Z_{k}\#N\overline{\mathbb{C}P^2}$ is diffeomorphic to the 4-manifold obtained from $Z\#N\overline{\mathbb{C}P^2}$ by twisting $X'$ using the diffeomorphism $f^k$. Therefore, $f^k$ $(k\neq 0)$ does not extend to any self-diffeomorphism of $X'$, and hence $(X',f)$ is an infinite order generalized cork. This completes the proof. 
\end{proof}

\begin{remark}For each $X$, the above proof gives infinitely many pairwise non-homeomorphic 4-manifolds having the properties of $X'$ in Proposition~\ref{sec:cork:prop:construction} by varying $r,s,m$ in the proof. Also, we can relax the assumption of Proposition~\ref{sec:cork:prop:construction} that $X$ admits an embedding into a symplectic 4-manifold with weakly convex boundary. Indeed, as seen from the proof, if $X\natural E$ admits an embedding into a closed  4-manifold $Z$ with $b_2^+>1$ and non-vanishing Seiberg--Witten invariant such that the torus $T$ in $E$ represents an infinite order second homology class of $Z$, then the conclusion still holds. 
\end{remark}

\section{Isomorphisms between symmetric bilinear forms}\label{sec:bilinear}
In this section, under certain conditions, we construct isomorphisms between symmetric bilinear forms that preserve certain maps modeled on genus functions and related invariants. 

Throughout this section, we use the following setting. For each $i=1,2$, let $A_i, B_i$ be finitely generated $\mathbb{Z}$-modules, and set $A'_i=A_i\oplus B_i$. Let $p_{A_i}:A'_i\to A_i$ and $p_{B_i}:A'_i\to B_i$ denote the projections. 
Take an integral symmetric bilinear form $Q_{A_i}: A_i\times A_i\to \mathbb{Z}$, and define  symmetric bilinear forms $Q_{B_i}:B_i\times B_i\to \mathbb{Z}$ and $Q_i:A'_i\times A'_i\to \mathbb{Z}$ by $Q_{B_i}(b,b')=0$ and $Q_i(a+b, a'+b')=Q_{A_i}(a,a')$, where $a, a'\in A_i$ and $b,b'\in B_i$. We set $\widetilde{A}_i=A_i/\mathrm{Tor}$ and assume that the symmetric bilinear form $Q_{\widetilde{A}_i}:\widetilde{A}_i\times \widetilde{A}_i\to \mathbb{Z}$ induced by each $Q_{A_i}$ is non-degenerate. 
These maps are models of intersection forms of 4-manifolds. 

We show the following lemma.  

\begin{lemma}\label{sec:bilinear:lem:projection} Suppose the modules $A_1$ and $A_2$ or the modules  $B_1$ and $B_2$ are torsion-free. 
	If there exists an isomorphism $\Phi:A'_1\to A'_2$ that preserves $Q_1$ and $Q_2$, then $p_{A_2}\circ\Phi|_{A_1}:A_1\to A_2$ is an isomorphism that preserves $Q_{A_1}$ and $Q_{A_2}$, and $p_{B_2}\circ\Phi|_{B_1}:B_1\to B_2$ is an isomorphism that preserves $Q_{B_1}$ and $Q_{B_2}$. 
\end{lemma}
\begin{proof}
	Due to the definition of $Q_i$, we can check that $\mathrm{rank}\, A_i$ is equal to the maximum rank of a submodule $C_i$ of $A'_i$ such that the symmetric bilinear form  $(C_i/\mathrm{Tor})\times (C_i/\mathrm{Tor})\to \mathbb{Z}$ induced by $Q_i$ is non-degenerate. This fact implies  $\mathrm{rank}\, A_1=\mathrm{rank}\, A_2$ and also $\mathrm{rank}\, B_1=\mathrm{rank}\, B_2$, since the bilinear forms on $A'_1/\mathrm{Tor}$ and $A'_2/\mathrm{Tor}$ are isomorphic. Due to the assumption that $A_1, A_2$ or $B_1, B_2$ are torsion-free, it follows that $A_1, B_1$ are respectively isomorphic to $A_2, B_2$. 
	
	Let each $T_{i}$ be the torsion submodule of $A_i$. Note that $T_1=T_2=0$ in the case where $A_1, A_2$ are torsion-free. Let each $F_i$ be a free submodule of $A_i$ such that $A_i=F_i\oplus T_i$. We note that $F_1, T_1$ are isomorphic to $F_2, T_2$, respectively. 
	Let $M_i$ be a matrix representing the bilinear form $F_i\times F_i\to \mathbb{Z}$ induced by $Q_{A_i}$. Then each $\det M_i$ is independent of the choices of $F_i$ and its basis, and $\det M_i$ is non-zero due to the non-degeneracy of $Q_{\widetilde{A}_i}$. 
	
	Let each $p_{F_i}:A'_i=F_i\oplus T_i\oplus B_i\to F_i$ be the projection. We here show that  $p_{F_2}\circ \Phi|_{F_1}:F_1\to F_2$ is an isomorphism. Let $u_1,\dots, u_k$ be a basis of $F_1$, and define $u_1',\dots, u_k'\in F_2$  by $u_j'=p_{F_2}(\Phi(u_j))$. We define the matrix $M_1$ in the last paragraph by  $M_1=(Q_{A_1}(u_i,u_j))$.  
	Since the projection $p_{F_2}$ preserves $Q_2,Q_{A_2}$, the homomorphism $p_{F_2}\circ \Phi|_{F_1}:F_1\to F_2$ preserves $Q_{A_1},Q_{A_2}$, showing that the matrix $(Q_2(u_i',u_j'))$ is equal to $M_1$. Due to the condition $\det M_1\neq 0$, we easily see that $u_1',\dots,u_k'$ is linearly independent (cf.\ \cite[Corollary~1.2.13.]{GS}). We thus have a basis $v_1,\dots,v_k$ of $F_2$ and positive integers $a_1,\dots, a_k$ such that $a_1v_1,\dots, a_kv_k$ generates the submodule generated by $u_1',\dots,u_k'$. It follows that $\det M_1=m\cdot \det M_2$ for some positive integer $m$. By reversing the roles of $A_1$ and $A_2$ and applying the same argument, we see that $\det M_1=\det M_2$. Consequently, $u_1',\dots,u_k'$ is a basis of $F_2$, and thus $p_{F_2}\circ \Phi|_{F_1}:F_1\to F_2$ is an isomorphism. 
	
	We now show that $p_{A_2}\circ \Phi|_{A_1}:A_1\to A_2$ is an isomorphism that preserves $Q_{A_1}, Q_{A_2}$. The projection $p_{A_2}$ clearly preserves $Q_2, Q_{A_2}$, and hence $p_{A_2}\circ \Phi|_{A_1}$ preserves $Q_{A_1}, Q_{A_2}$. 
	Due to the assumption that $A_1, A_2$ or $B_1, B_2$ are torsion-free, the isomorphism $\Phi$ bijectively maps $T_{1}$ to $T_{2}$, implying that $p_{A_2}\circ \Phi|_{T_1}:T_1\to T_2$ is an isomorphism. 
	We define a submodule $F_2'$ of $A_2$ by $F_2'=p_{A_2}(\Phi(F_1))$. 
	Then the composition of $p_{A_2}\circ \Phi|_{F_1}:F_1\to F_2'$ and the projection $A_2=F_2\oplus T_2\to F_2$ is equal to the isomorphism $p_{F_2}\circ \Phi|_{F_1}:F_1\to F_2$. This fact implies that the surjection $p_{A_2}\circ \Phi|_{F_1}:F_1\to F_2'$ is an isomorphism, and thus 
	$p_{F_2}|_{F_2'}:F_2'\to F_2$ is also an isomorphism. It follows $A_2=F_2'\oplus T_2$, implying $p_{A_2}\circ \Phi|_{A_1}(A_1)=A_2$. Therefore, $p_{A_2}\circ \Phi|_{A_1}:A_1\to A_2$ is an isomorphism. 
	
	Finally, we show that $p_{B_2}\circ \Phi|_{B_1}:B_1\to B_2$ is an isomorphism that preserves $Q_{B_1}, Q_{B_2}$. Due to the assumption that  $A_1, A_2$, or $B_1, B_2$, are torsion-free, $\Phi$ bijectively maps the torsion submodule of $B_1$ to that of $B_2$. We can thus easily  check  $p_{B_2}\circ \Phi|_{B_1}(B_1)$ is an isomorphism that preserves  $Q_{B_1}, Q_{B_2}$, by using the non-degeneracy assumption of each $Q_{\widetilde{A}_i}$.
\end{proof}

In the rest of this section, let $R$ be an ordered set, and let  $G_i:A'_i\to R$ be an arbitrary map for each $i=1,2$. These maps are not necessarily homomorphisms and are models of genus functions and related invariants. 
We will construct isomorphisms that preserve the restrictions of $Q_1, Q_2$ and $G_1, G_2$ by using  Lemma~\ref{sec:bilinear:lem:projection}.

\begin{proposition}\label{sec:bilinear:prop:genus:zero-part}  
	Suppose there exists an isomorphism $\Phi:A'_1\to A'_2$ that preserves $Q_1,Q_2$ and $G_1,G_2$. If $A_1$ and $A_2$ are torsion-free, then there exists an isomorphism between $B_1$ and $B_2$ that preserves  $Q_{B_1},Q_{B_2}$ and $G_1|_{B_1},G_2|_{B_2}$. 
\end{proposition}
\begin{proof}By Lemma~\ref{sec:bilinear:lem:projection}, the map $p_{B_2}\circ \Phi|_{B_1}:B_1\to B_2$ is an isomorphism that preserves $Q_{B_1}$, $Q_{B_2}$. 
Since  $A_1, A_2$ are torsion-free, and $Q_{A_1}, Q_{A_2}$ are non-degenerate, we see $\Phi(B_1)=B_2$. The isomorphism $p_{B_2}\circ \Phi|_{B_1}$ thus preserves $G_1|_{B_1},G_2|_{B_2}$. 	
\end{proof}

\begin{proposition}\label{sec:bilinear:prop:genus} Assume that each $G_i$ satisfies $G_i(a)\leq G_i(a+b)$ for any $a \in A_i$ and $b\in B_i$, and further that there exists an isomorphism $\Phi:A'_1\to A'_2$  that preserves $Q_1,Q_2$ and $G_1,G_2$. If  the modules $A_1$ and $A_2$ or the modules  $B_1$ and $B_2$ are torsion-free, then there exists an isomorphism between $A_1$ and $A_2$ that preserves $Q_{A_1},Q_{A_2}$ and $G_1|_{A_1},G_2|_{A_2}$.  
\end{proposition}

\begin{proof}
By Lemma~\ref{sec:bilinear:lem:projection}, we see that $p_{A_2}\circ \Phi|_{A_1}:A_1\to A_2$ is an isomorphism that preserves $Q_{A_1}$, $Q_{A_2}$. 
	 Let $a$ be an element of $A_1$ and set $\widehat{a}=p_{A_2}\circ \Phi(a)$. We will prove $G_1(a)=G_2(\widehat{a})$, showing that $p_{A_2}\circ \Phi|_{A_1}:A_1\to A_2$ is a desired isomorphism. We divide the proof into two cases. Note that  $\Phi(a)=\widehat{a}+\widehat{a}_B$ for some $\widehat{a}_B\in B_2$.

	(1) The case where $A_1, A_2$ are torsion-free.     
	In this case, we see $\widehat{a}_B=\Phi(b)$ for some $b\in B_1$ due to the non-degeneracy of $Q_{\widetilde{A}_1}=Q_{A_1}$. By the assumptions on $G_1,G_2,\Phi$, we see     
	\begin{equation*}
		G_1(a) = G_2(\widehat{a}+\widehat{a}_B)\geq G_2(\widehat{a})=G_1(a-b) \geq G_1(a),  
	\end{equation*}  
	giving $G_1(a)=G_2(\widehat{a})$. 
	
	(2) The case where $B_1, B_2$ are torsion-free. 
	We see  $\widehat{a}_B=\Phi(b+t)$ for some element $b\in B_1$ and some torsion element $t\in A_1$ due to the non-degeneracy of $Q_{\widetilde{A}_1}$. We denote the torsion element $\Phi(t)\in A_2$ by $\widehat{t}$. For each integer $n$, we have   $\Phi(a+nt)=\widehat{a}+n\widehat{t}+\widehat{a}_B$ and  $\Phi(a+(n-1)t-b)=\widehat{a}+n\widehat{t}$. 
	By the assumptions on $G_1,G_2,\Phi$, we thus see
	\begin{alignat*}{2}
		G_1(a+nt) = G_2(\widehat{a}+n\widehat{t}+\widehat{a}_B) &\geq G_2(\widehat{a}+n\widehat{t})\\
		&=G_1(a+(n-1)t-b)\geq G_1(a+(n-1)t),      
	\end{alignat*}
	and hence we have 
	\begin{equation*}
		G_1(a)\geq G_2(\widehat{a})\geq G_1(a-t)\geq G_1(a-2t)\cdots\geq G_1(a-nt)\geq \cdots.
	\end{equation*}
	Since $t$ is a torsion, this implies $G_1(a)=G_2(\widehat{a})$. This completes the proof. 
\end{proof}


\section{Quasi-invertible cobordisms}\label{sec:quasi-invertible}
In this section, we introduce the notion of quasi-invertible cobordisms to study genus functions. We also provide their examples and study their algebraic properties. Throughout this section, \textit{we allow manifolds, submanifolds, and their boundaries to be disconnected}. For such manifolds, the uniqueness of connected sums and boundary sums do not necessarily hold, but our results hold for any choices of sums unless otherwise stated. So, we will use $\#$ and $\natural$ for these operations. Let $M,N$ be  closed 3-manifolds. 

\subsection{Definitions and examples}\label{sec:quasi-invertible:subsec:example} 
Let us recall that a 4-dimensional 1-handlebody of genus $k$ $(k\geq 0)$ is the $k$-fold boundary sum $\natural_k (S^1\times D^3)$. For convenience, we use the following terminologies. 

\begin{definition}
Suppose a 4-dimensional 1-handlebody of genus $k$ is embedded in the interior of a 4-manifold. An operation of removing the 1-handlebody from the 4-manifold and gluing $\natural_k S^2\times D^2$ back along the resulting boundary component is called a \textit{surgery along a $4$-dimensional $1$-handlebody}. The reverse operation is called a \textit{surgery along an embedded copy of $\natural_k (S^2\times D^2)$}.   
\end{definition}

Note that the way to glue $\natural_k S^2\times D^2$ is not unique in general, but we use a convenient way in each context. By contrast, the reverse operation 
is uniquely determined by an embedded $\natural_k (S^2\times D^2)$ (\cite{LP72}). 

Here we introduce notions of a quasi-product cobordism, a quasi-invertible cobordism and an $H_2$-surjective cobordism. 

\begin{definition}(1) A cobordism $R$ from $M$ to $N$ is called \textit{quasi-product} if it is orientedly diffeomorphic to the product $I\times M$ possibly after performing surgeries along finitely many embedded $4$-dimensional $1$-handlebodies. If furthermore the homomorphism $H_2(M)\to H_2(R)$ induced by the inclusion is an isomorphism, then $R$ is called \textit{strongly quasi-product}. 

	(2) A cobordism $P$ from $M$ to $N$ is called \textit{quasi-invertible} (resp.\ \textit{strongly quasi-invertible}) if it admits a cobordism $Q$ from $N$ to $M$ such that a glued 4-manifold $P\cup_N Q$ is a quasi-product cobordism (resp.\ strongly quasi-product cobordism) from the boundary component $M$ of $P$ to that of $Q$. The cobordism $Q$ is called a \textit{quasi-inverse cobordism} (resp.\ \textit{strongly quasi-inverse cobordism}) of $P$. 
	
	(3) A cobordism $P$ from $M$ to $N$ is called \textit{$H_2$-surjective} if the homomorphism $H_2(M)\to H_2(P)$ induced by the inclusion is surjective.  
\end{definition}

As we will later show, there are rich examples of quasi-invertible cobordisms. To produce examples, let us recall a few terminologies. A \textit{surgery on a circle} embedded in a 4-manifold is an operation of removing the tubular neighborhood $S^1\times D^3$ of the circle and gluing $S^2\times D^2$ back along the resulting boundary component. In general, there are two ways to glue $S^2\times D^2$ along the boundary component (\cite{Gl62}), but we use a convenient way in each context. Clearly, any surgery on a circle is realized by a surgery along an embedded $4$-dimensional $1$-handlebody of genus one. 
The \textit{surgery on a $2$-sphere} embedded in a 4-manifold with trivial normal bundle is an operation of removing the tubular neighborhood $S^2\times D^2$ of the 2-sphere from the 4-manifold and gluing $S^1\times D^3$ back along the resulting boundary. We note that the resulting 4-manifold is uniquely determined by the embedded 2-sphere (\cite{LP72}). An \textit{$S^2$-link} in a 4-manifold is a submanifold of the interior of the 4-manifold whose each connected component is a 2-sphere. 

We use the following notation. 

\begin{definition}\label{sec:quasi-invertible:def:S^4(L)}
For an $S^2$-link $L$ in $S^4$, let $S^4(L)$ be the closed 4-manifold obtained from $S^4$ by performing the surgery on $L$, that is, by performing the surgery on each component of $L$. Note that the normal bundle of any 2-sphere embedded in $S^4$ is trivial. 
\end{definition}

It is easy to see that $n(S^1\times S^3)$ is an example of $S^4(L)$ for each $n\geq 0$. 
We here study quasi-product cobordisms. Clearly, $I\times M$ is strongly quasi-product. Furthermore, we have simple examples below. 
\begin{lemma}\label{sec:quasi-invertible:lem:quasi-product}
	$(1)$ A $4$-manifold is a quasi-product cobordism from $M$ to $M$ if and only if the $4$-manifold is diffeomorphic to a $4$-manifold obtained from $I\times M$ by performing surgeries along finitely many submanifolds, each of which is diffeomorphic to the boundary sum of copies of $S^2\times D^2$. 
	
	$(2)$ For any  $S^2$-link $L$ in $S^4$,  the connected sum $(I\times M)\#S^4(L)$ is a strongly quasi-invertible cobordism from $M$ to $M$. 
\end{lemma}
\begin{proof}The claim (1) is straightforward. We have $H_2(S^4(L))=0$ by Lemma~\ref{sec:genus:sum:lem:torsion-free},  so the claim (2) follows from  (1). 
\end{proof}
\begin{remark}\label{sec:quasi-invertible:rem:cobordism:non-strong}
	Quasi-product cobordisms are not necessarily strongly quasi-product. For example, let $P$ be the 4-manifold obtained from $I\times S^1\times S^2$ by performing the surgery on the 2-sphere $a \times b\times S^2$ ($a\in \mathrm{int}\, I$, $b\in S^1$). 
	Then $P$ is a quasi-product cobordism from $S^1\times S^2$ to $S^1\times S^2$, but the homomorphism $H_2(0\times S^1\times S^2)\to H_2(R)$ induced by the inclusion is the zero-map, and thus $P$ is not strongly quasi-product. 
\end{remark}

Now we produce examples of quasi-invertible cobordisms. Clearly, invertible cobordisms are strongly quasi-invertible, and quasi-product cobordisms are quasi-invertible. To produce varieties of examples, we recall the notion of a strongly slice link. 
\begin{definition}\label{sec:quasi-invertible:def:strongly_slice}
A (1-dimensional) link $\ell$ in $N$ is called \textit{strongly slice} if $\ell$ in $1\times N$ bounds a disjoint union of smoothly and properly embedded disks in $I\times N$. We note that each component of a strongly slice link is null-homologous, and thus its 0-framing is defined as the Seifert framing. 
\end{definition}

In Subsection~\ref{sec:genus:sum:subsec:example}, we will produce various examples of codimension-$0$ submanifolds of $S^4(L)$. Using those examples, Lemma~\ref{sec:quasi-invertible:lem:quasi-product} and the lemma below, we can generate various examples of (strongly) quasi-invertible cobordisms from product cobordisms. 

\begin{lemma}\label{sec:quasi-invertible:lem:quasi-invertible:generate} Let $L$ be an $S^2$-link in $S^4$, and let $V$ be a codimension-$0$ submanifold of $S^4(L)$. Then, for any quasi-invertible $($resp.\ strongly quasi-invertible$)$ cobordism $P$ from $M$ to $N$,  the following hold.  
	
	$(1)$ Every  codimension-$0$ submanifold of $P$ containing the boundary component $M$ of $P$ is a quasi-invertible $($resp.\ strongly quasi-invertible$)$ cobordism  from $M$ to the other boundary component. 
	
	$(2)$ $P\#V$ is a quasi-invertible $($resp.\ strongly quasi-invertible$)$ cobordism  from $M$ to the disjoint union $N\sqcup \partial V$. 
	
	$(3)$ If $\partial V$ is non-empty, then $P\natural V$ is a quasi-invertible $($resp.\ strongly quasi-invertible$)$ cobordism from $M$ to $N\#\partial V$, where the boundary sum is taken for the boundary component $N$ of $P$.

	$(4)$ Every $4$-manifold obtained by attaching $1$-handles to $P$ along a connected component of $N$ is a quasi-invertible $($resp.\ strongly quasi-invertible$)$ cobordism from $M$ to the new boundary component.

	$(5)$ Every $4$-manifold obtained by attaching homotopically canceling pairs of $1$- and $2$-handles to $P$ along $N$ is a quasi-invertible $($resp.\ strongly quasi-invertible$)$ cobordism from $M$ to the new boundary component.

	$(6)$ Every $4$-manifold obtained by attaching  $2$-handles to $P$ along a strongly slice link  $\ell$ in $N$ with $0$-framings is a quasi-invertible $($resp.\ strongly quasi-invertible$)$ cobordism from $M$ to the $0$-surgery on $\ell$ in $N$.

	Furthermore,  the quasi-invertible cobordisms given by these claims are  invertible in the case where $P, L$ are respectively an invertible cobordism and an empty set.  
\end{lemma}

\begin{proof}For the claim (1), attach the quasi-inverse cobordism of $P$ to the exterior of the submanifold. The resulting manifold gives a quasi-inverse cobordism of the submanifold, implying the claim (1). The claims (2) and (3) are straightforward from Lemma~\ref{sec:quasi-invertible:lem:quasi-product} and the claim (1).  For (4), by further attaching a canceling 2-handle to each 1-handle, we obtain a 4-manifold diffeomorphic to $P$. The claim (4) thus follows from (1).   For (5), by further attaching the upside down cobordism of the cobordism consisting of the handle pairs, we obtain a 4-manifold diffeomorphic to $P$ due to Lemma~\ref{sec:cork:thm:inverse cobordism}. The claim (5) thus follows from (1).  The claim (6) is easily checked as follows. The link $\ell$ in $0\times N$ bounds a disjoint union of smoothly embedded disks in $I\times N$, and the tubular neighborhood of these disks give the 2-handles attached to $P$. Hence, the 4-manifold in (6) can be embedded into $P\cup _{\mathrm{id}_N}(I\times N)$, and thus the claim (6) follows from (1).  The last claim easily follows from these constructions. 
\end{proof}

\subsection{Algebraic properties}\label{sec:quasi-invertible:subsec:algebraic}
We begin with quasi-product cobordisms.

\begin{lemma}\label{sec:quasi-invertible:subsec:algebraic:lem:product}
	Let $R$ be a quasi-product cobordism from $M$ to $N$. Then $R$ is $H_2$-surjective, and the homomorphism $H_1(M)\to H_1(R)$  induced by the inclusion is injective. Consequently, $R$ is strongly quasi-product if and only if the homomorphism $H_2(M)\to H_2(R)$ induced by the inclusion is injective.  
\end{lemma}
\begin{proof}Let $\widehat{R}$ be the 4-manifold obtained from $R$ by performing surgeries along finitely many embedded 1-handlebodies $\natural_{k_i} S^1\times D^3$ so that $\widehat{R}$ is diffeomorphic to $I\times M$. 
	We remove the interior of these 1-handlebodies from $R$ and denote the resulting 4-manifold by $R_c$. Note that $R=\widehat{R}=R_c$ in the case where $R$ is a product cobordism. Let ${j}:M\hookrightarrow {R}$, $j_{c}:M\hookrightarrow R_c$, $h:R_c\hookrightarrow R$, $\widehat{h}:R_c\hookrightarrow \widehat{R}$ and  $\widehat{j}:M\hookrightarrow \widehat{R}$  be the inclusions.

	By the excision theorem and the Poincar\'{e} duality, we see that $H_3(\widehat{R},R_c)\cong \oplus H_3(\natural_{k_i}(S^2\times D^2), k_i(S^1 \times S^2))\cong 0$. 
	The homology exact sequence for the pair $(\widehat{R}, R_c)$ thus shows that $\widehat{h}_*:H_2(R_c)\to H_2(\widehat{R})$ is injective. Since the inclusion $\widehat{j}=\widehat{h}\circ j_{c}:M\hookrightarrow \widehat{R}$ is a homotopy equivalence, the injection $\widehat{h}_*$ and also ${j_c}_*:H_2(M)\to H_2(R_c)$ are isomorphisms.

	By the excision theorem and the Poincar\'{e} duality, we see that $H_2(R,R_c)\cong \oplus H_2(\natural_{k_i}(S^1 \times D^3), k_i(S^1 \times S^2))\cong 0$. The homology exact sequence for the pair $({R}, R_c)$ thus shows that ${h}_*:H_2(R_c)\to H_2(R)$ is surjective and that ${h}_*:H_1(R_c)\to H_1(R)$ is injective.  
	
	Now we can easily prove the claims. 	
	Since ${h}_*:H_2(R_c)\to H_2(R)$ is surjective, and ${j_c}_*:H_2(M)\to H_2(R_c)$ is bijective, the composition $j_*=h_*\circ {j_c}_*:H_2(M)\to H_2(R)$ is also surjective. Consequently, $R$ is $H_2$-surjective. 
	Since the inclusion $\widehat{j}=\widehat{h}\circ j_{c}:M\hookrightarrow \widehat{R}$ is a homotopy equivalence, we see that ${j_{c}}_*:H_1(M)\hookrightarrow H_1(R_c)$ is injective, and thus the composition $j_*=h_*\circ {j_{c}}_*:H_1(M)\to H_1(R)$ is also injective. 
\end{proof}

In a special case, we can easily check the strongness.   

\begin{corollary}\label{sec:quasi-invertible:subsec:algebraic:cor:strongly:raional}
	If $M$ is a disjoint union of rational homology $3$-spheres, then any quasi-product $($resp.\ quasi-invertible$)$ cobordism from $M$ to $N$ is strongly quasi-product $($resp.\ strongly quasi-invertible$)$. 
\end{corollary}
\begin{proof}Since the rank of $H_1(M)$ is zero, we see  $H_2(M)=0$ by the Poincar\'{e} duality and the universal coefficient theorem. The claim thus follows from Lemma~\ref{sec:quasi-invertible:subsec:algebraic:lem:product}. 
\end{proof}
We next discuss quasi-invertible cobordisms. In the rest of this subsection, let $P$ be a quasi-invertible cobordism from $M$ to $N$, and let $Q$ be a quasi-inverse cobordism of $P$. We denote a quasi-product cobordism $P\cup_N Q$ by $R$. In the case where $P$ is strongly quasi-invertible, we assume $R$ is strongly quasi-product. 
Let $I^M_P$ and $K^R_P$ respectively denote the image of the homomorphism $H_2(M)\to H_2(P)$ and the kernel of the homomorphism $H_2(P)\to H_2(R)$, where these are the homomorphisms induced by the inclusions. 
Then the following lemma holds. 

\begin{lemma}\label{sec:quasi-invertible:subsec:algebraic:lem:quasi-invertible}

	$(1)$ The homomorphism $H_1(M)\to H_1(P)$ induced by the inclusion is injective, and the decomposition $H_2(P)=I^M_P+K^R_P$ holds.
	
	$(2)$ If $P$ is strongly quasi-invertible, then the homomorphism $H_2(M)\to H_2(P)$ induced by the inclusion is injective, and the decomposition $H_2(P)=I^M_P\oplus K^R_P$ holds.
\end{lemma}
\begin{proof}
	Since the homomorphism $H_*(M)\to H_*(R)$ is the composition of homomorphisms $H_*(M)\to H_*(P)$ and $H_*(P)\to H_*(R)$, the claims follow from Lemma~\ref{sec:quasi-invertible:subsec:algebraic:lem:product}. 
\end{proof}

Finally, we discuss algebraic effects of an attachment of a quasi-invertible cobordisms. Let $X$ be a 4-manifold having a boundary component orientedly diffeomorphic to $M$, and let $X'$ denote a 4-manifold $X\cup_M P$. We attach $Q$ to $X'$ so that the submanifold $R=P\cup_N Q$ is a quasi-product cobordism. In the case where $P$ is strongly quasi-invertible (resp.\ invertible), we assume $R$ is strongly quasi-product (resp.\ product). Let $X''$  denote the resulting 4-manifold $X\cup_M P\cup_N Q$, and let $\iota:X\hookrightarrow X'$ and $\iota_P:P\hookrightarrow X'$ be the inclusions. Then the following properties hold. 

\begin{lemma}\label{sec:quasi-invertible:subsec:algebraic:lem:X'}
$(1)$ The homomorphism $H_2(X)\to H_2(X'')$ induced by the inclusion is surjective. 

$(2)$ $H_2(X')=\iota_*(H_2(X))+{\iota_P}_*(K^R_P)$.

$(3)$ The restriction of the intersection form of $X'$ to ${\iota_P}_*(K^R_P)$ is the zero-map. 

$(4)$  $\iota_*:H_2(X)\to H_2(X')$ is surjective if and only if $P$ is $H_2$-surjective. 
\end{lemma}
\begin{proof}$(1)$ By  the homology exact sequence for the pair $(R, M)$ and Lemma~\ref{sec:quasi-invertible:subsec:algebraic:lem:product}, we see $H_2(R,M)=0$.  
	The excision theorem shows $H_2(X'',X)\cong H_2(R,M)=0$. The homology exact sequence for the pair $(X'', X)$ thus shows the claim.  
	
	$(2)$ The claim (1) and the sequence $X\subset X'\subset X''$ implies $H_2(X')=\iota_*(H_2(X))+K^{X''}_{X'}$, where $K^{X''}_{X'}$ denotes the kernel of the homomorphism $H_2(X')\to H_2(X'')$ induced by the inclusion. We note ${\iota_P}_*(K^R_P)\subset K^{X''}_{X'}$. Using the Mayer--Vietoris exact sequence for the decomposition $X'=X\cup_M P$ and  Lemma~\ref{sec:quasi-invertible:subsec:algebraic:lem:quasi-invertible}, we can check $\iota_*(H_2(X))+K^{X''}_{X'}=\iota_*(H_2(X))+{\iota_P}_*(K^R_P)$, showing the claim. 

$(3)$ By Lemma~\ref{sec:quasi-invertible:subsec:algebraic:lem:product}, the intersection form of $R$ is the zero-map, and thus the claim follows.  

	$(4)$ Lemma~\ref{sec:quasi-invertible:subsec:algebraic:lem:quasi-invertible} and the Mayer--Vietoris exact sequence for  $X'=X\cup_M P$ imply that $\iota_*:H_1(X)\to H_1(X')$ is injective. It thus follows from the homology exact sequence for the pair $(X',X)$ that  $\iota_*:H_2(X)\to H_2(X')$ is surjective if and only if $H_2(X',X)=0$.  
	Similarly, by Lemma~\ref{sec:quasi-invertible:subsec:algebraic:lem:quasi-invertible}, $P$ is $H_2$-surjective if and only if $H_2(P,M)=0$. The claim thus follows from the fact $H_2(X',X)\cong H_2(P,M)$. 
\end{proof}

The useful corollary below is straightforward from the above (4). 
\begin{corollary}\label{sec:quasi-invertible:subsec:algebraic:cor:X':isom}
	Suppose $\iota_*:H_2(X)\to H_2(X')$ is injective. Then,  $\iota_*:H_2(X)\to H_2(X')$ is an isomorphism if and only if  $P$ is $H_2$-surjective. 
\end{corollary}

In the case where the homomorphism $H_2(X)\to H_2(X'')$ is injective, nice properties hold.  

\begin{proposition}\label{sec:quasi-invertible:subsec:algebraic:prop:decomposition:injective}
	Suppose the homomorphism $H_2(X)\to H_2(X'')$ induced by the inclusion is injective.  Then, $\iota_*:H_2(X)\to H_2(X')$ is injective, and the decomposition $H_2(X')=\iota_*(H_2(X))\oplus {\iota_P}_*(K^R_P)$ holds. 
\end{proposition}
\begin{proof}The injectivity is clear from the sequence $X\subset X'\subset X''$, and the decomposition easily follows from Lemma~\ref{sec:quasi-invertible:subsec:algebraic:lem:X'}. 
\end{proof}

We give sufficient conditions that the homomorphism $H_2(X)\to H_2(X'')$ is injective. 
\begin{proposition}\label{sec:quasi-invertible:subsec:algebraic:prop:injective:sufficient}
	Suppose at least one of the following conditions \textnormal{(i)}--\textnormal{(iii)} holds. 
	\begin{enumerate}
		\item [\textnormal{(i)}] $P$ is strongly quasi-invertible. 
		\item [\textnormal{(ii)}] The homomorphism $H_2(M)\to H_2(X)$  induced by the inclusion is the zero-map. 
		\item [\textnormal{(iii)}] The intersection form of $X$ is non-degenerate, and $H_2(X)$ has no torsion.
	\end{enumerate}	
	Then the homomorphism $H_2(X)\to H_2(X'')$ induced by the inclusion is injective. 
\end{proposition}
\begin{proof}In the case where (i) or (ii) holds, the Mayer--Vietoris exact sequence for $X''=X\cup_M R$ implies that the homomorphism $H_2(X)\to H_2(X'')$ is injective.  In the case where (iii) holds, every class of $H_2(X)$ algebraically intersects with some class of $H_2(X)$, implying that the condition (ii) also holds. Hence, the claim follows.  
\end{proof}

The lemma below is useful for applying the results in Section~\ref{sec:Genus functions and quasi-invertible cobordisms}. 
\begin{lemma}\label{sec:quasi-invertible:subsec:algebraic:lem:torsion-free}
	Suppose at least one of the following conditions \textnormal{(i)}--\textnormal{(iii)} holds. 
\begin{enumerate}
	\item [\textnormal{(i)}] $P$ is strongly quasi-product. 
	\item [\textnormal{(ii)}] $P$ is invertible. 
	\item [\textnormal{(iii)}] $P$ is strongly quasi-invertible, and $N$  is a disjoint union of rational homology $3$-spheres.
\end{enumerate}	
Then, $H_2(P)$ and ${\iota_P}_*(K^R_P)$ are torsion-free. 
\end{lemma}
\begin{proof}The Poincar\'{e} duality and the universal coefficient theorem imply that $H_2(N)$ and $H_2(M)$ are torsion-free. Since $P$ is strongly quasi-inverible, we see that $H_2(R)$ is also torsion-free. 
	
	We consider the Mayer--Vietoris exact sequence for $R=P\cup_N Q$, and we will show that the homomorphism $H_2(N)\to H_2(P)\oplus H_2(Q)$ in the sequence is injective. In the case where (i) holds, we may assume that $Q=I\times N$, and hence this homomorphism is injective.   In the case where (iii) holds, we see $H_2(N)=0$ as shown in the proof of Corollary~\ref{sec:quasi-invertible:subsec:algebraic:cor:strongly:raional}, and hence the homomorphism is injective.  So, we prove the injectivity in the case where (ii) holds. In this case, $R$ is diffeomorphic to $I\times M$, and hence the inclusion $M\hookrightarrow R$ is a homotopy equivalence. Due to the sequence $M\subset P\subset R$, this fact implies that the homomorphism $H_3(P)\to H_3(R)$ induced by the inclusion is surjective.  The Mayer--Vietoris exact sequence for $R=P\cup_N Q$ thus shows that the homomorphism $H_2(N)\to H_2(P)\oplus H_2(Q)$ is injective. 
	
	Now we are ready to prove the claim. Since $H_2(N)$ and $H_2(R)$ are torsion-free,  and the homomorphism $H_2(N)\to H_2(P)\oplus H_2(Q)$ is injective, the Mayer--Vietoris exact sequence for $R=P\cup_N Q$ implies that  $H_2(P)$ is torsion-free, and hence $K^R_P$ is also  torsion-free. 
	By Lemma~\ref{sec:quasi-invertible:subsec:algebraic:lem:quasi-invertible}, we have $H_2(P)=I^M_P\oplus K^R_P$. 
 It thus follows from the Mayer--Vietoris exact sequence for $X'=X\cup_M P$ that the restriction ${\iota_P}_*|_{K^R_P }:K^R_P\to {\iota_P}_*(K^R_P)$ is an isomorphism. Therefore, ${\iota_P}_*(K^R_P)$ is torsion-free. 
\end{proof}

In the case where the homomorphism $H_2(X)/\mathrm{Tor}\to H_2(X'')/\mathrm{Tor}$ is injective, we still have nice properties.  
\begin{proposition}\label{sec:quasi-invertible:subsec:algebraic:prop:decomposition:injective:torsion}
	Suppose the homomorphism $H_2(X)/\mathrm{Tor}\to H_2(X'')/\mathrm{Tor}$ induced by the inclusion is injective. Then,  $\iota_*:H_2(X)/\mathrm{Tor}\to H_2(X')/\mathrm{Tor}$ is injective, and the decomposition  $H_2(X')/\mathrm{Tor}=\iota_*(H_2(X))/\mathrm{Tor}\oplus {\iota_P}_*(K^R_P)/\mathrm{Tor}$ holds. 
\end{proposition}
\begin{proof}The proof is similar to that of Proposition~\ref{sec:quasi-invertible:subsec:algebraic:prop:decomposition:injective}.  
\end{proof}

We give sufficient conditions that the homomorphism $H_2(X)/\mathrm{Tor}\to H_2(X'')/\mathrm{Tor}$ is injective. 
\begin{proposition}\label{sec:quasi-invertible:subsec:algebraic:prop:injective:sufficient:torsion}
Suppose at least one of the following conditions \textnormal{(i)}--\textnormal{(iii)} holds. 
	\begin{enumerate}
		\item [\textnormal{(i)}] $P$ is strongly quasi-invertible. 
		\item [\textnormal{(ii)}] The homomorphism $H_2(M)\to H_2(X)/\mathrm{Tor}$ induced by the inclusion is the zero-map. 
		\item [\textnormal{(iii)}] The intersection form of $X$ is non-degenerate.
	\end{enumerate}	
	Then the homomorphism $H_2(X)/\mathrm{Tor}\to H_2(X'')/\mathrm{Tor}$ induced by the inclusion is injective. 
\end{proposition}
\begin{proof}The proof is similar to that of  Proposition~\ref{sec:quasi-invertible:subsec:algebraic:prop:injective:sufficient}. 
\end{proof}
\section{Genus functions and quasi-invertible cobordisms}\label{sec:Genus functions and quasi-invertible cobordisms}
In this section, we study behavior of genus functions of 4-manifolds under attachment of quasi-invertible cobordisms. We also show stabilities of algebraic inequivalences of genus functions under the attachments. Throughout this section, \textit{we allow manifolds, submanifolds, and their boundaries to be disconnected unless otherwise stated, but we still assume surfaces are connected.} 

\subsection{Behavior of genus functions}\label{subsec:Genus functions and quasi-invertible cobordisms:behavior}
Throughout this subsection, we use the following setting. Let $M,N$ be closed 3-manifolds, and let $P$ be a quasi-invertible cobordism from $M$ to $N$. Let $X$ be a connected 4-manifold having a boundary component orientedly diffeomorphic to $M$, and we denote a glued  4-manifold $X\cup_M P$ by $X'$. Let $\iota:X\hookrightarrow X'$ be the inclusion. 

Under these assumptions, we prove the following theorem. 

\begin{theorem}\label{sec:genus:quasi-invertible:thm:genus preserving}
	The inclusion $\iota:X\hookrightarrow X'$ preserves the genus functions. 
\end{theorem}

To prove this theorem, we show the following lemma, which reveals an important property of genus functions. 
\begin{lemma}\label{sec:genus:quasi-invertible:lem:genus:3-handle}
	$(1)$ If a  $4$-manifold $Z'$ is obtained from a connected $4$-manifold $Z$ by attaching a $3$- or a $4$-handle, then the inclusion $Z\hookrightarrow Z'$ preserves the genus functions. 
	
	$(2)$ If a $4$-manifold $Z'$ is obtained from a connected $4$-manifold $Z$ by attaching $\natural_n(S^1\times D^3)$ along a boundary component diffeomorphic to $n(S^1\times S^2)$ for some $n\geq 0$, then the inclusion $Z\hookrightarrow Z'$ preserves the genus functions. 	
\end{lemma}
\begin{proof}(1) We prove the claim in the case where the attached handle is a 3-handle $D^3\times D^1$. In the case of a 4-handle, 
	the claim follows by a much simpler argument. Let $j:Z\hookrightarrow Z'$ denote the inclusion. We note that the kernel of $j_*:H_2(Z)\to H_2(Z')$ is the image of the homomorphism $H_2(S^2\times D^1)\to H_2(Z)$ induced by the inclusion as seen from the Mayer--Vietoris exact sequence, where $S^2\times D^1$ denotes the attaching region of the 3-handle.  
	
	Let $\alpha$ be a class of $H_2(Z)$ and define the class $\alpha'$ of $H_2(Z')$ by $\alpha'={j}_*(\alpha)$. We set $g=g_Z(\alpha)$ and $g'=g_{Z'}(\alpha')$. Then $\alpha$ (resp.\ $\alpha'$) is represented by a closed surface $\Sigma$ (resp.\ $\Sigma'$) of genus $g$ (resp.\ $g'$) embedded in the interior of $Z$ (resp.\ $Z'$). Clearly, $\alpha'$ is also represented by the surface $\Sigma$, showing $g\geq g'$. 
	
	We next show $g\leq g'$. Since the cocore of the 3-handle is $0\times D^1$, we may assume that $\Sigma'$ is disjoint from a tubular neighborhood of the cocore by taking a general position. Therefore, we can isotope $\Sigma'$ into the interior of $Z$. We thus see that $\alpha'=j_*([\Sigma'])$ and hence $j_*(\alpha)=j_*([\Sigma'])$. It follows that $\alpha=[\Sigma']+\beta$ in $H_2(Z)$ for some $\beta\in \mathrm{Ker}\, j_*$. Here note that any class of $\mathrm{Ker}\, j_*$ is represented by a disjoint union of 2-spheres in the attaching region of the 3-handle. By pushing these 2-spheres into the interior of $Z$ and taking a connected sum of $\Sigma'$ with these 2-spheres, we can construct a closed surface of genus $g'$ embedded in  $Z$ which represents $\alpha$. Hence, we see $g\leq g'$. Therefore, we have  $g=g'$, showing that $j$ preserves the genus functions. 
	
	(2) Since $Z'$ is obtained from $Z$ by attaching 3-and 4-handles, the claim follows from (1). 
\end{proof}

In the rest of this subsection, we use the following setting to prove Theorem~\ref{sec:genus:quasi-invertible:thm:genus preserving} and its corollaries. We attach a quasi-inverse cobordism $Q$ of $P$ to $X'$ so that the submanifold $R=P\cup_N Q$ is a quasi-product cobordism, and we  denote the resulting 4-manifold $X\cup_M P\cup_N Q$ by  $X''$. 

\begin{proof}[Proof of Theorem~$\ref{sec:genus:quasi-invertible:thm:genus preserving}$]
	The cobordism $R$ becomes diffeomorphic to $I\times M$ possibly after performing surgeries along finitely many embedded 4-dimensional 1-handlebodies. Thus, the 4-manifold $\widehat{X}''$ obtained from $X''$ by these surgeries is diffeomorphic to $X\cup_{\mathrm{id}_{M}} (I\times M)$. We denote by $\widehat{j}:X\hookrightarrow \widehat{X}''$ the inclusion. By shrinking $I\times M$ together with the collar neighborhood of $\partial X$, we can easily construct a diffeomorphism $\varphi:\widehat{X}''\to X$ such that $\varphi\circ \widehat{j}:X\to X$ is homotopic to $\mathrm{id}_X$.  Since the diffeomorphisms $\varphi$ and $\mathrm{id}_X$ preserves the genus functions, it follows that $\widehat{j}:X\hookrightarrow \widehat{X}''$ also preserves the genus functions. 	
	
	We remove the interior of the aforementioned 1-handlebodies in $R$ from $X''$ and denote the resulting 4-manifold by $X''_c$.  Note that $X''=\widehat{X}''=X''_c$ in the case where $R$ is a product cobordism. Let  $j:X\hookrightarrow X''$, $j':X'\hookrightarrow X''$, $j_{c}:X\hookrightarrow X''_c$, $h:X''_c\hookrightarrow X''$ and $\widehat{h}:X''_c\hookrightarrow \widehat{X}''$ be the inclusions. Since $j_c$ and $\widehat{h}$ are inclusions, and $\widehat{j}=\widehat{h}\circ j_{c}$ preserves the genus functions, every $\alpha\in H_2(X)$ satisfies 
	\begin{equation*}
		g_X(\alpha)\geq g_{X''_c}({j_{c}}_*(\alpha))\geq g_{\widehat{X}''}(\widehat{h}_*({j_{c}}_*(\alpha)))=g_X(\alpha), 
	\end{equation*}
	showing that  $j_{c}:X\hookrightarrow X''_c$ also preserves the genus functions. By Lemma~\ref{sec:genus:quasi-invertible:lem:genus:3-handle}, the inclusion $h:X''_c\hookrightarrow X''$ preserves the genus functions, and thus $j=h\circ j_{c}:X\hookrightarrow  X''$ also preserves the genus functions. 
	Due to this fact and the relation $j=j'\circ \iota$, we see that every $\alpha\in H_2(X)$ satisfies 
		\begin{equation*}
		g_X(\alpha)\geq g_{X'}({\iota}_*(\alpha))\geq g_{{X}''}(j'_*({\iota}_*(\alpha)))=g_X(\alpha), 
	\end{equation*}
	showing that $\iota:X\hookrightarrow X'$ preserves the genus functions. 
\end{proof}

\begin{remark}\label{sec:genus:quasi-invertible:remark:genus property}
	(1) By Theorem~$\ref{sec:genus:quasi-invertible:thm:genus preserving}$ and Lemma~\ref{sec:quasi-invertible:lem:quasi-invertible:generate}, the genus functions are preserved under attachment of a 1-handle whose attaching region belongs to a connected boundary component and also under attachment of 2-handles along a strongly slice link with 0-framings. Applying this result, we can show non-sliceness for various knots in $S^3$.  
	
(2)	To prove a result similar to Theorem~\ref{sec:genus:quasi-invertible:thm:genus preserving} for other invariants in Section~\ref{sec:genus type}, we note that the properties of genus functions used in the above proof are the following  (i)--(iii). (i) Every (orientation-preserving) diffeomorphism preserves the genus functions. (ii) For any (orientation-preserving) embedding, the homomorphism induced by the embedding does not increase the values of genus functions. (iii) The property stated in   Lemma~\ref{sec:genus:quasi-invertible:lem:genus:3-handle}.(1).  We note that, in the case where $P$ is an invertible cobordism, the above proof of Theorem~\ref{sec:genus:quasi-invertible:thm:genus preserving} does not require  Lemma~\ref{sec:genus:quasi-invertible:lem:genus:3-handle}.(1). 

(3) There are many examples of cobordisms that do not have the property of Theorem~\ref{sec:genus:quasi-invertible:thm:genus preserving}. For instance, many cobordisms do not have this property after taking a connected sum with $S^2\times S^2$.  We can also construct examples of such homology cobordisms as follows. Take a (contractible) cork $(C,f)$ admitting a 4-manifold $Z$ such that its cork twist along $(C,f)$ changes the algebraic equivalent class of the genus function of $Z$. Note that there are various examples of such corks (e.g.\ \cite{AM97}, \cite{AY08}, \cite{AY12},  \cite{G17GT}). Then, $\overline{C}-\textnormal{int}\, D^4$ is a homology cobordism from $\partial C$ to $S^3$, and the attachment of this cobordism to  $\overline{Z-\textnormal{int}\, C}$ using either $f$ or $\mathrm{id}_{\partial C}$ does not preserve the genus function of $\overline{Z-\textnormal{int}\, C}$ due to the assumption on $Z$.    
\end{remark}

The corollary below is straightforward. 

\begin{corollary}\label{sec:genus:quasi-invertible:cor:inclusion:isomorphism}
	If ${\iota}_*:H_2(X)\to H_2(X')$ is an isomorphism, then the genus function of ${X'}$ is algebraically equivalent to that of ${X}$. 
\end{corollary}

The corollaries below provide sufficient conditions that $g_{X'}$ is algebraically equivalent to $g_X$. 

\begin{corollary}\label{sec:genus:quasi-invertible:cor:inclusion:invertible homology cobordism}
	If $P$ is an invertible homology cobordism from $M$ to $N$, then the genus function of ${X'}$ is algebraically equivalent to that of ${X}$. 
\end{corollary}
\begin{proof} This is straightforward from Corollary~\ref{sec:genus:quasi-invertible:cor:inclusion:isomorphism}. 
\end{proof}
\begin{corollary}\label{sec:genus:quasi-invertible:cor:algebraic equivalence:sufficient}
Suppose that $P$ is $H_2$-surjective and that at least one of the following conditions \textnormal{(i)}--\textnormal{(iii)} holds. 
\begin{enumerate}
\item [\textnormal{(i)}] $P$ is strongly quasi-invertible. 
\item [\textnormal{(ii)}] The homomorphism $H_2(M)\to H_2(X)$  induced by the inclusion is the zero-map. 
\item [\textnormal{(iii)}] The intersection form of $X$ is non-degenerate, and $H_2(X)$ has no torsion.
\end{enumerate}	
Then   the genus function of ${X'}$ is algebraically equivalent to that of ${X}$. 
\end{corollary}
\begin{proof} This is straightforward from Corollary~\ref{sec:genus:quasi-invertible:cor:inclusion:isomorphism}, Lemma~\ref{sec:quasi-invertible:subsec:algebraic:lem:X'}, and  Propositions~\ref{sec:quasi-invertible:subsec:algebraic:prop:injective:sufficient} and \ref{sec:quasi-invertible:subsec:algebraic:prop:decomposition:injective}. 
\end{proof}

In the case where $P$ is $H_2$-surjective, the genus function $g_{X'}$ is completely determined by $g_X$ and $\iota_*$ due to Theorem~\ref{sec:genus:quasi-invertible:thm:genus preserving} and Lemma~\ref{sec:quasi-invertible:subsec:algebraic:lem:X'}. In the case where $P$ is not $H_2$-surjective, we give estimates on the values of $g_{X'}$ for all classes. As in Subsection~\ref{sec:quasi-invertible:subsec:algebraic}, let $K^R_P$ be the kernel of the homomorphism $H_2(P)\to H_2(R)$ induced by the inclusion, and  let $\iota_P:P\hookrightarrow X'$ be the inclusion. Due to Lemma~\ref{sec:quasi-invertible:subsec:algebraic:lem:X'}, we have $H_2(X')=\iota_*(H_2(X))+{\iota_P}_*(K^R_P)$. 
We prove the following estimates for $g_{X'}$.  

\begin{theorem}\label{sec:genus:quasi-invertible:thm:estimate}
	Every classes $\alpha\in H_2(X)$ and $\beta\in K^R_P$ satisfy      
	\begin{equation*}
		g_{X}(\alpha)\leq g_{X'}(\iota_*(\alpha)+{\iota_P}_*(\beta))\leq g_X(\alpha)+g_{P}(\beta). 
	\end{equation*}
\end{theorem}
\begin{proof}The right inequality is obvious, so we prove the left inequality. 
	Let  $j:X\hookrightarrow  X''$ and $j':X'\hookrightarrow  X''$ be the inclusions. By the definition of $K^R_P$, we see $j'_*(\iota_*(\alpha)+{\iota_P}_*(\beta))=j_*(\alpha)$, showing $g_{X'}(\iota_*(\alpha)+{\iota_P}_*(\beta))\geq g_{X''}(j_*(\alpha))$. Since $j_*$ preserves the genus functions due to Theorem~\ref{sec:genus:quasi-invertible:thm:genus preserving}, we obtain $g_{X'}(\iota_*(\alpha)+{\iota_P}_*(\beta))\geq g_X(\alpha)$. 
\end{proof}
As seen from Lemma~\ref{sec:quasi-invertible:lem:quasi-invertible:generate}, there are many examples of quasi-invertible cobordisms $P$ such that $g_P$ is the zero-map, and for such examples, Theorem~\ref{sec:genus:quasi-invertible:thm:estimate} shows that $g_{X'}$ is completely determined by $g_X$ even when $P$ is not $H_2$-surjective.  

Under an additional assumption, the inclusion furthermore preserves the torsion-free genus functions.

\begin{corollary}\label{sec:genus:quasi-invertible:cor:genus preserving:torsion-free}
Suppose the homomorphism $H_2(X)/\mathrm{Tor}\to H_2(X'')/\mathrm{Tor}$ induced by the inclusion is injective. Then, the inclusion $\iota:X\hookrightarrow X'$ preserves the torsion-free genus functions. 
\end{corollary}
\begin{proof}Let $A$ be a class of $H_2(X)/\mathrm{Tor}$,  and  let $\alpha$ be a class of $H_2(X)$ that  represents $A$ and satisfies $g_X(\alpha)=g^*_X(A)$.  By Theorem~\ref{sec:genus:quasi-invertible:thm:genus preserving}, we see 
	\begin{equation*}
		g^*_X(A)=g_X(\alpha)=g_{X'}({\iota}_*(\alpha))\geq  g^*_{X'}({\iota}_*(A)). 
	\end{equation*}
Let $\alpha'$ be a class of $H_2(X')$ that represents ${\iota}_*(A)\in H_2(X')/\mathrm{Tor}$ and satisfies $g_{X'}(\alpha')=g^*_{X'}({\iota}_*(A))$. 
Proposition~\ref{sec:quasi-invertible:subsec:algebraic:prop:decomposition:injective:torsion} implies that there exist torsion classes $\tau\in H_2(X)$ and $\beta\in K^R_P$ satisfying $\alpha'={\iota}_*(\alpha)+{\iota}_*(\tau)+{\iota_P}_*(\beta)$.  Theorem~\ref{sec:genus:quasi-invertible:thm:estimate} and the above inequality thus imply   
    \begin{equation*}
    g^*_{X'}({\iota}_*(A))=g_{X'}(\alpha')\geq g_{X}(\alpha+\tau)\geq  g^*_{X}(A)\geq  g^*_{X'}({\iota}_*(A)), 
    \end{equation*}
	showing $g^*_{X'}({\iota}_*(A))=g^*_{X}(A)$. Therefore, the inclusion $\iota:X\hookrightarrow X'$ preserves the torsion-free genus functions.
\end{proof}
We remark that the assumption of the above corollary holds if $P$ is strongly quasi-invertible. For other sufficient conditions, see Proposition~\ref{sec:quasi-invertible:subsec:algebraic:prop:injective:sufficient:torsion}.

\subsection{Stabilities for algebraic inequivalences of genus functions}\label{subsec:Genus functions and quasi-invertible cobordisms:stabilities}
Throughout this subsection, we use the following setting. For $i=1,2$, let $M_i,N_i$ be closed 3-manifolds, let $P_i$ be a quasi-invertible cobordism from $M_i$ to $N_i$, and let $Q_i$ be its quasi-inverse cobordism. Also, let $X_i$ be a connected 4-manifold having a boundary component orientedly diffeomorphic to $M_i$, and we denote a glued 4-manifold $X_i\cup_{M_i} P_i$ by $X'_i$. We attach $Q_i$ to $X_i'$ so that the submanifold $R_i=P_i\cup_{N_i} Q_i$ is a quasi-product cobordism, and we denote the resulting 4-manifold $X_i\cup_{M_i} P_i\cup_{N_i} Q_i$ by $X''$. Let $\iota_i:X_i\hookrightarrow X'_i$ and $j_i:X_i\hookrightarrow X''_i$ be the inclusions. Furthermore, we define the subgroup ${{\iota_i}_{P_i}}_*(K^{R_i}_{P_i})$ of $H_2(X_i')$ as in Subsection~\ref{sec:quasi-invertible:subsec:algebraic}. 

We will show that algebraic inequivalences of (torsion-free) genus functions are stable under attaching a quasi-invertible cobordism. Our results thus show that exotic 4-manifolds with algebraically inequivalent genus functions always yield many more  (orientedly) exotic 4-manifolds by attaching quasi-invertible cobordisms. 

We first prove the following theorem, which shows stabilities by taking the contrapositions. 
\begin{theorem}\label{sec:genus:quasi-invertible:thm:stablity:H2zero}
	$(1)$ Suppose each ${\iota_i}_*:H_2(X_i)\to H_2(X'_i)$ is an isomorphism. Then, the genus functions of $X_1'$ and $X_2'$ are algebraically equivalent  if and only if those of $X_1$ and $X_2$ are algebraically equivalent.  
	
	$(2)$ Suppose that each ${\iota_i}_*:H_2(X_i)/\mathrm{Tor}\to H_2(X'_i)/\mathrm{Tor}$ is an isomorphism and that each ${j_i}_*:H_2(X_i)/\mathrm{Tor}\to H_2(X''_i)/\mathrm{Tor}$ is injective. Then, the torsion-free genus functions of $X_1'$ and $X_2'$ are algebraically equivalent if and only if those of $X_1$ and $X_2$ are algebraically equivalent.  
\end{theorem}
\begin{proof}The claim (1) is straightforward from Corollary~\ref{sec:genus:quasi-invertible:cor:inclusion:isomorphism}. For (2), the torsion-free genus functions of each $X_i'$ is algebraically equivalent to that of $X_i$  due to  Corollary~\ref{sec:genus:quasi-invertible:cor:genus preserving:torsion-free}. The claim (2) thus follows. 
\end{proof}
We remark that the assumptions of this theorem hold in the case where each $P_i$ is  strongly quasi-invertible and $H_2$-surjective. See Propositions~\ref{sec:quasi-invertible:subsec:algebraic:prop:injective:sufficient} and \ref{sec:quasi-invertible:subsec:algebraic:prop:injective:sufficient:torsion} and Lemma~\ref{sec:quasi-invertible:subsec:algebraic:lem:X'} for other sufficient conditions. 

For general quasi-invertible cobordisms, we give a stability of algebraic inequivalences 
under additional assumptions. See Propositions~\ref{sec:quasi-invertible:subsec:algebraic:prop:decomposition:injective} and \ref{sec:quasi-invertible:subsec:algebraic:prop:injective:sufficient} and Lemma~\ref{sec:quasi-invertible:subsec:algebraic:lem:torsion-free} for sufficient conditions that the assumption below holds. 

\begin{theorem}\label{sec:genus:quasi-invertible:thm:stablity:non-degenerate}
	Suppose the intersection form of each $X_i$ is non-degenerate. Suppose further that the groups $H_2(X_1)$ and $H_2(X_2)$ or the groups ${\iota_{P_1}}_*(K^{R_1}_{P_1})$ and ${\iota_{P_2}}_*(K^{R_2}_{P_2})$ are torsion-free and that each ${j_i}_*:H_2(X_i)\to H_2(X''_i)$ is injective. 
If the  genus functions of  ${X'_1}$ and $X'_2$ are algebraically equivalent, then those of ${X_1}$ and $X_2$ are algebraically equivalent. If furthermore $H_2(X_1)$ and $H_2(X_2)$  are torsion-free, then the restrictions of the genus functions of ${X'_1}$ and $X'_2$ respectively to ${\iota_{P_1}}_*(K^{R_1}_{P_1})$  and ${\iota_{P_2}}_*(K^{R_2}_{P_2})$ are also algebraically equivalent. 
\end{theorem}
\begin{proof}Due to Propositions~\ref{sec:bilinear:prop:genus} and \ref{sec:bilinear:prop:genus:zero-part}, the claims follow from Lemma~\ref{sec:quasi-invertible:subsec:algebraic:lem:X'}, Proposition~\ref{sec:quasi-invertible:subsec:algebraic:prop:decomposition:injective} and  Theorems~\ref{sec:genus:quasi-invertible:thm:genus preserving} and \ref{sec:genus:quasi-invertible:thm:estimate}.  
\end{proof}

For torsion-free genus functions, we give a similar stability under a much simpler assumption.  

\begin{theorem}\label{sec:genus:quasi-invertible:thm:stablity:non-degenerate:torsion-free}
	Suppose the intersection form of each $X_i$ is non-degenerate. If the torsion-free genus functions of ${X'_1}$ and $X'_2$ are algebraically equivalent, then those of ${X_1}$ and $X_2$ are algebraically equivalent, and those of ${X'_1}$ and ${X'_2}$  restricted respectively to ${\iota_{P_1}}_*(K^{R_1}_{P_1})/\mathrm{Tor}$  and ${\iota_{P_2}}_*(K^{R_2}_{P_2})/\mathrm{Tor}$ are algebraically equivalent. 
\end{theorem}
\begin{proof}
	Due to Proposition~\ref{sec:quasi-invertible:subsec:algebraic:prop:decomposition:injective:torsion},  Theorem~\ref{sec:genus:quasi-invertible:thm:estimate} and Corollary~\ref{sec:genus:quasi-invertible:cor:genus preserving:torsion-free}, we easily see that every $\alpha\in {\iota_i}_*(H_2(X_i))/\mathrm{Tor}$ and $\beta \in {\iota_{P_1}}_*(K^{R_1}_{P_1})/\mathrm{Tor}$ satisfy $g^*_{X'_i}(\alpha)\leq g^*_{X'_i}(\alpha+\beta)$ for each $i$ by an argument similar to the proof of  Corollary~\ref{sec:genus:quasi-invertible:cor:genus preserving:torsion-free}.  
	The claim thus follows from Propositions~\ref{sec:bilinear:prop:genus}, \ref{sec:bilinear:prop:genus:zero-part} and  \ref{sec:quasi-invertible:subsec:algebraic:prop:decomposition:injective:torsion} and  Corollary~\ref{sec:genus:quasi-invertible:cor:genus preserving:torsion-free}.  
\end{proof}

We remark that Theorems~\ref{sec:genus:quasi-invertible:thm:stablity:non-degenerate} and \ref{sec:genus:quasi-invertible:thm:stablity:non-degenerate:torsion-free} do not hold if the non-degeneracy assumption is removed as seen from Remark~\ref{sec:genus:remark:sum}.

\section{Genus functions, connected sums and boundary sums}\label{sec:genus:sum}
In this section, we study behavior of genus functions under connected sums and boundary sums with a certain type of 4-manifolds, by applying results in the last section. Throughout this section, \textit{we allow manifolds, submanifolds, and their boundaries to be disconnected unless otherwise stated.} 
Our results for boundary sums hold for any choice of a boundary component unless otherwise stated, and hence we do not specify the boundary component. 

\subsection{Submanifolds of $S^4(L)$}\label{sec:genus:sum:subsec:example}
Before discussing genus functions, we show that various 4-manifolds are realized as codimension-$0$ submanifolds of the surgered 4-manifold $S^4(L)$ (see Definition~\ref{sec:quasi-invertible:def:S^4(L)}) for some $S^2$-link $L$ in $S^4$. Note that $S^4(L)=S^4$ in the case where $L$ is an empty set.

\begin{lemma}\label{sec:genus:sum:lem:submanifold:generate} Let $L, L'$  be $S^2$-links in $S^4$, and let $V, V'$ be connected codimension-$0$ submanifolds of $S^4(L), S^4(L')$ with non-empty boundaries, respectively. Then, the following hold. 
	
	$(1)$ Any boundary sum $V\natural V'$ admits an embedding into $S^4(L\sqcup L')$, where $L\sqcup L'$ denotes the split union of $L$ and $L'$ in $S^4\#S^4\cong S^4$. 

$(2)$ If a $4$-manifold is obtained from $V$ by attaching $1$-handles along a connected component of $\partial V$, then the $4$-manifold admits an embedding into $S^4(L)$. 

$(3)$ If a $4$-manifold is obtained from $V$ by attaching homotopically canceling pairs of $1$- and $2$-handles, then the $4$-manifold admits an embedding into $S^4(L)$. 

$(4)$ If  a $4$-manifold is obtained from $V$ by attaching $2$-handles along a strongly slice link in $\partial V$ with $0$-framings, then the $4$-manifold admits an embedding into $S^4(L)$. 

$(5)$ Let $\ell$ be a link in $\partial(\natural_nS^2\times D^2)$ that bounds a disjoint union of smoothly and properly embedded disks in $\natural_nS^2\times D^2$.  Suppose a $4$-manifold is obtained from $\natural_n S^1\times D^3$ by attaching $2$-handles  along a link $\ell$ with the framings induced from the disks, where $\ell$ is regarded as a link in $\partial(\natural_nS^1\times D^3)$ using an arbitrary diffeomorphism $\partial(\natural_nS^2\times D^2)\to \partial(\natural_nS^1\times D^3)$. Then the  $4$-manifold admits an embedding into $S^4$. 
\end{lemma}
\begin{proof}The claim (1) is straightforward, and  the claims (2)--(4) are proved similarly to the proof of  Lemma~\ref{sec:quasi-invertible:lem:quasi-invertible:generate}, so  we prove  (5). The tubular neighborhood of the disks bounding $\ell$ form 2-handles attached  $\natural_n S^1\times D^3$. Since $\natural_nS^2\times D^2$ becomes diffeomorphic to $S^4$ by attaching $\natural_n S^1\times D^3$ for any gluing map (\cite{LP72}), the claim (5) follows. 
\end{proof}
Using the above lemma, we can generate numerous examples. Here we provide simple examples of codimension-$0$ submanifolds of $S^4$. 
\begin{example}\label{sec:genus:sum:example:simple}
	Let $\Sigma_g$ denote a closed surface of genus $g$. Then $\Sigma_g\times D^2$ $(g\geq 0)$,  $\natural_n S^1\times D^3$ $(n\geq 0)$, Mazur-type contractible 4-manifolds, their boundary sums, and their connected sums admit embeddings into $S^4$ and hence are realized as codimension-$0$ submanifolds of $S^4$. Furthermore, these 4-manifolds with 2-handles attached along strongly slice links with 0-framings are realized as codimension-$0$ submanifolds of $S^4$. 
\end{example}

We here give sufficient conditions that the second homology group of a submanifold of $S^4(L)$ is torsion-free. 
\begin{lemma}\label{sec:genus:sum:lem:torsion-free}
	$(1)$ $H_2(S^4(L))=0$. 
	
	$(2)$ For any codimension-$0$ submanifold $Z$ of $S^4$, $H_2(Z)$ is torsion-free. 
	
	$(3)$  If the boundary of a codimension-$0$ submanifold $Z$ of $S^4(L)$ is a disjoint union of rational homology $3$-spheres, then $H_2(Z)$ is torsion-free. 
\end{lemma}
\begin{proof}
 We first prove the claim (1). Let $E(L)$ be the exterior of $L$ in $S^4$. By the homology exact sequence for the pair $(S^4, E(L))$ and the excision theorem, we see that $H_2(E(L))=0$. We can thus easily check the claim using the homology exact sequence for the pair $(S^4(L), E(L))$ and the excision theorem. 

 We next prove the claims (2) and (3). In the case where $Z=S^4$, the claim (2) is trivial. So, we assume $\partial Z$ is non-empty. Then, by Lemma~\ref{sec:quasi-invertible:lem:quasi-invertible:generate}, $Z-\mathrm{int}\, D^4$ is a quasi-invertible cobordism from $S^3$ to $\partial Z$, and this cobordism is invertible in the case where $Z$ is a submanifold of $S^4$.  The claims (2) and (3) thus follow from  Lemma~\ref{sec:quasi-invertible:subsec:algebraic:lem:torsion-free}, since $H_2(Z-\mathrm{int}\, D^4)$ is isomorphic to $H_2(Z)$. 
\end{proof}

\subsection{Behavior of genus functions}
Now we discuss genus functions. Throughout this subsection, let $X$ be a $($possibly closed$)$ connected 4-manifold, let $L$ be an $S^2$-link in $S^4$, and let $Z$ be a connected summand or  a codimension-$0$ connected submanifold of the surgered 4-manifold $S^4(L)$. We show that genus functions are preserved under connected sums and boundary sums with $Z$.

\begin{theorem}\label{sec:genus:sum:thm:genus preserving}
The inclusion $H_2(X)\hookrightarrow  H_2(X\#Z)=H_2(X)\oplus H_2(Z)$ preserves the genus functions of $X$ and $X\#Z$. 
Furthermore, in the case where $\partial X$ and $\partial Z$ are non-empty, the inclusion $X\hookrightarrow X\natural Z$ also preserves the genus functions.  
\end{theorem}
Note that, due to the definition $X\#Z=(X-\mathrm{int}\, D^4)\cup_{\partial D^4}(Z-\mathrm{int}\, D^4)$, we identify  $H_2(X\#Z)$ with $H_2(X)\oplus H_2(Z)$ by using the isomorphisms $H_2(X-\mathrm{int}\, D^4)\to H_2(X)$ and $H_2(Z-\mathrm{int}\, D^4)\to H_2(Z)$ induced by the inclusions. 

\begin{proof}[Proof of Theorem~$\ref{sec:genus:sum:thm:genus preserving}$]
	
Since $Z-\mathrm{int}\, (D^4\sqcup D^4)$ is diffeomorphic to $Z\#(I\times S^3)$,  Lemma~\ref{sec:quasi-invertible:lem:quasi-invertible:generate} shows that $Z-\mathrm{int}\, (D^4\sqcup D^4)$ is a strongly quasi-invertible cobordism from $S^3$ to the disjoint union $S^3\sqcup \partial Z$.  
The 4-manifold  $X\#(Z-\mathrm{int}\, D^4)$ is constructed from $X-\mathrm{int}\, D^4$ by attaching this cobordism, 
so the inclusion $X-\mathrm{int}\, D^4 \hookrightarrow X\#(Z-\mathrm{int}\, D^4)$ preserves the genus functions due to Theorem~\ref{sec:genus:quasi-invertible:thm:genus preserving}.  We note that the inclusions $X-\mathrm{int}\, D^4\hookrightarrow X$ and $X\#(Z-\mathrm{int}\, D^4)\hookrightarrow X\#Z$ preserve the genus functions by Lemma~\ref{sec:genus:quasi-invertible:lem:genus:3-handle}. Therefore, the  inclusion $H_2(X)\hookrightarrow  H_2(X\#Z)$ preserves the genus functions. 

Now we assume $\partial X$ and $\partial Z$ are non-empty.  By  Lemma~\ref{sec:quasi-invertible:lem:quasi-invertible:generate},  any boundary sum $(I\times \partial X)\natural Z$ is a strongly quasi-invertible cobordism, where the sum is taken along $1\times \partial X$. Since  $X\natural Z$ is constructed from $X$ by attaching this cobordism, the inclusion $X\hookrightarrow X\natural Z$ preserves the genus functions due to Theorem~\ref{sec:genus:quasi-invertible:thm:genus preserving}. 
\end{proof}
In the case where $H_2(Z)=0$,  the above inclusions in Theorem~\ref{sec:genus:sum:thm:genus preserving} are isomorphisms, and hence $g_{X\#Z}$ and $g_X$ are completely determined by $g_X$. In the case where $H_2(Z)\neq 0$, we give the following estimates. 

\begin{theorem}\label{sec:genus:sum:thm:estimate}
	Every $\alpha\in H_2(X)$ and $\beta\in H_2(Z)$ satisfy  
	\begin{equation*}
		g_X(\alpha)\leq g_{X\#Z}(\alpha+\beta)\leq g_X(\alpha)+g_Z(\beta).
	\end{equation*}
	Furthermore, in the case where  $\partial X$ and $\partial Z$ are non-empty, they also satisfy  
	\begin{equation*}
		g_X(\alpha)\leq g_{X\natural Z}(\alpha+\beta)\leq g_X(\alpha)+g_Z(\beta).
	\end{equation*}
\end{theorem}
\begin{proof}Due to the constructions of $X\#Z$ and $X\natural Z$ in the proof of Theorem~\ref{sec:genus:sum:thm:genus preserving}, the claim follows from Theorem~\ref{sec:genus:quasi-invertible:thm:estimate}. 
\end{proof}

Due to the constructions of $X\#Z$ and $X\natural Z$ in the proof of Theorem~\ref{sec:genus:sum:thm:genus preserving}, Corollary~\ref{sec:genus:quasi-invertible:cor:genus preserving:torsion-free} shows that the torsion-free genus functions are preserved under connected sums and boundary sums with $Z$. 
\begin{corollary}\label{sec:genus:sum:cor:genus preserving:torsion-free}
The inclusion $H_2(X)\hookrightarrow  H_2(X\#Z)=H_2(X)\oplus H_2(Z)$ preserves the torsion-free genus functions of $X$ and $X\#Z$. 
Furthermore, in the case where both $\partial X$ and $\partial Z$ are non-empty, the inclusion $X\hookrightarrow X\natural Z$ also preserves  the torsion-free genus functions.  
\end{corollary}

We remark that the inclusions $H_2(Z)\hookrightarrow  H_2(X\#Z)=H_2(X)\oplus H_2(Z)$ and $Z\hookrightarrow X\natural Z$ do not necessarily preserve the genus functions. For example, in the case where $X$ is $S^2\times S^2$ or $S^2\times S^2-\mathrm{int}\, D^4$, we have many counterexamples of $Z$ as easily seen from results of \cite{AY12}, \cite{AY13}.  See also Remark~\ref{sec:genus:remark:sum}. 
\subsection{Stabilities for algebraic inequivalences of genus functions}\label{sec:genus:sum:subsec:stability}
Throughout this subsection, for each $i=1,2$, let $X_i$ be a $($possibly closed$)$ connected 4-manifold, let $L_i$ be an $S^2$-link in $S^4$, and 
let $Z_i$ be a connected summand or a codimension-$0$ connected submanifold of $S^4(L_i)$. 

We will show that algebraic inequivalences of genus functions are stable under connected sums and boundary sums. Our results thus show that exotic 4-manifolds with algebraically inequivalent genus functions always yield many more (orientedly) exotic 4-manifolds by taking connected sums and/or boundary sums.

We first show stabilities of algebraic inequivalences in the case where $H_2(Z_i)=0$. Specifically, we prove the following theorem, which is a restatement of Theorem~\ref{sec:intro:thm:sum:iff} and shows the stabilities by taking the contrapositions. 

\begin{theorem}\label{sec:genus:sum:thm:equivalent:sum}
Suppose each $H_2(Z_i)=0$. Then  the following hold. 

$(1)$ The genus functions of $X_1\#Z_1$ and $X_2\#Z_2$ are algebraically equivalent if and only if those of $X_1$ and $X_2$ are algebraically equivalent. 

$(2)$ Assume each $\partial X_i$ and $\partial Z_i$ are non-empty. Then, the genus functions of $X_1\natural Z_1$ and $X_2\natural Z_2$ are algebraically equivalent if and only if those of $X_1$ and $X_2$ are algebraically equivalent. 
\end{theorem}
\begin{proof}Under the assumptions, Theorem~\ref{sec:genus:sum:thm:genus preserving} shows that the genus functions  of each $X_i\#Z_i$ and $X_i\natural Z_i$ are algebraically equivalent to that of  $X_i$, showing the claims. 
\end{proof}

We also prove stabilities for torsion-free genus functions.  
\begin{theorem}\label{sec:genus:sum:thm:equivalent:sum:torsion-free}
	Suppose each $b_2(Z_i)=0$. Then  the following hold. 
	
	$(1)$ The torsion-free genus functions of $X_1\#Z_1$ and $X_2\#Z_2$ are algebraically equivalent if and only if those of $X_1$ and $X_2$ are algebraically equivalent. 
	
	$(2)$ Assume each $\partial X_i$ and $\partial Z_i$ are non-empty.  Then, the torsion-free genus functions of $X_1\natural Z_1$ and $X_2\natural Z_2$ are algebraically equivalent if and only if those of $X_1$ and $X_2$ are algebraically equivalent. 
\end{theorem}
\begin{proof}The claims follow from Corollary~\ref{sec:genus:sum:cor:genus preserving:torsion-free} by the same argument as in the proof of Theorem~\ref{sec:genus:sum:thm:equivalent:sum}. 
\end{proof}

In the case where $b_2(Z_i)\neq 0$, we prove similar stabilities of algebraic inequivalences under additional assumptions. 
Specifically, we prove the following theorem which is a refinement of Theorem~\ref{sec:intro:thm:equivalent:non-degenerate:sum}. See Lemma~\ref{sec:genus:sum:lem:torsion-free} for sufficient conditions that each $H_2(Z_i)$ is torsion-free.  

\begin{theorem}\label{sec:genus:sum:thm:equivalent:non-degenerate}
Assume that the intersection form of each $X_i$ is non-degenerate and that the groups $H_2(X_1)$ and $H_2(X_2)$  or the groups  $H_2(Z_1)$ and  $H_2(Z_2)$ are torsion-free. Then the following hold. 

$(1)$  Suppose that the genus functions of $X_1\#Z_1$ and $X_2\#Z_2$ are algebraically equivalent. Then, those of $X_1$ and $X_2$ are algebraically equivalent. If furthermore $H_2(X_1)$ and $H_2(X_2)$ are torsion-free, then the genus functions of $X_1\#Z_1$ and $X_2\#Z_2$ restricted respectively to $H_2(Z_1)$ and $H_2(Z_2)$ are algebraically equivalent. 

$(2)$ Suppose that each $\partial X_i$ and $\partial Z_i$ are non-empty and that the genus functions of $X_1\natural Z_1$ and $X_2\natural Z_2$ are algebraically equivalent. Then, those of $X_1$ and $X_2$ are algebraically equivalent. If furthermore $H_2(X_1)$ and $H_2(X_2)$ are torsion-free,
then the genus functions  of $X_1\natural Z_1$ and $X_2\natural Z_2$ restricted respectively to $H_2(Z_1)$ and $H_2(Z_2)$  are algebraically equivalent.   
\end{theorem}
\begin{proof}Since $H_2(S^4(L_i))=0$ for each $i$, the intersection form of each $Z_i$ is the zero-map. 
	Due to Propositions~\ref{sec:bilinear:prop:genus} and \ref{sec:bilinear:prop:genus:zero-part}, the claims thus follow from Theorems~\ref{sec:genus:sum:thm:genus preserving} and \ref{sec:genus:sum:thm:estimate}.  
\end{proof}

For torsion-free genus functions, we can prove stabilities without the torsion-free assumption.  
\begin{theorem}\label{sec:genus:sum:thm:equivalent:non-degenerate:torsion-free}
	Assume the intersection form of each $X_i$ is non-degenerate. Then the following hold. 
	
	$(1)$ If the torsion-free genus functions of  $X_1\#Z_1$ and $X_2\#Z_2$ are algebraically equivalent, then those of $X_1$ and $X_2$ are algebraically equivalent, and those of $X_1\#Z_1$ and $X_2\#Z_2$ restricted respectively to $H_2(Z_1)/\mathrm{Tor}$ and $H_2(Z_2)/\mathrm{Tor}$   are algebraically equivalent. 
	
	$(2)$ If each $\partial X_i$ and $\partial Z_i$ are non-empty, and the torsion-free genus functions of $X_1\natural Z_1$ and $X_2\natural Z_2$ are algebraically equivalent, then those of $X_1$ and $X_2$ are algebraically equivalent, and those of $X_1\natural Z_1$ and $X_2\natural Z_2$ restricted respectively to $H_2(Z_1)/\mathrm{Tor}$ and $H_2(Z_2)/\mathrm{Tor}$ are algebraically equivalent. 
\end{theorem}
\begin{proof}
	Since the intersection form of each $Z_i$ is the zero-map, the claims easily follow from Propositions~\ref{sec:bilinear:prop:genus} and \ref{sec:bilinear:prop:genus:zero-part}, Theorem~\ref{sec:genus:sum:thm:estimate} and Corollary~\ref{sec:genus:sum:cor:genus preserving:torsion-free}.  
\end{proof}

\begin{remark}\label{sec:genus:remark:sum}
	Theorems~\ref{sec:genus:sum:thm:equivalent:non-degenerate} and \ref{sec:genus:sum:thm:equivalent:non-degenerate:torsion-free} do not hold in the case where the intersection forms of $X_1,X_2$ are degenerate by the following simple reason. Let $X_1,X_2$ be an exotic pair of 4-manifolds with algebraically inequivalent genus functions, and assume they admit embeddings into $S^4$. Note that there are various examples of such pairs  (\cite{AY12}, \cite{AY13}), see also Theorem~\ref{sec:cork:exotic Stein}. Set $Z_1=X_2$ and $Z_2=X_1$. In this case, $X_1,X_2,Z_1,Z_2$ satisfy the assumptions of Theorems~\ref{sec:genus:sum:thm:equivalent:non-degenerate} and \ref{sec:genus:sum:thm:equivalent:non-degenerate:torsion-free} except the non-degeneracy assumption, but clearly the conclusions do not hold.  
\end{remark}

As a simple application, we show that algebraic equivalent classes of genus functions do not determine homeomorphism types even within closed 4-manifolds having isomorphic cohomology rings and arbitrarily large $b_2$.  

\begin{theorem}\label{sec:genus:sum:thm:closed:algebraically equivalent}For any closed connected $4$-manifold $X$, there exist infinitely many pairwise non-homeomorphic closed $4$-manifolds with pairwise algebraically equivalent genus functions such that their cohomology rings are isomorphic to that of $X\# S^1\times S^3$. 
\end{theorem}
\begin{proof}Let $K_1,K_2,\dots$ be infinitely many $S^2$-knots in $S^4$ such that their complements $S^4-K_1, S^4-K_2, \dots$ have pairwise non-isomorphic fundamental groups. As is well-known, such families of $S^2$-knots exist in abundance and easily produced, for example, as follows. Let $k$ be a knot in $S^3$, and let $K$ be the $S^2$-knot in $S^4$ derived from $k$ by Artin's spinning construction (\cite{Ar25}). Then $\pi_1(S^4-K)$ is isomorphic to $\pi_1(S^3-k)$ (see e.g.\ \cite[Section 2.1]{CKS}), and thus we obtain a desired infinite family of $S^2$-knots by varying a knot $k$. By the Seifert van-Kampen theorem, we see that each $\pi_1(S^4(K_i))$ is isomorphic to $\pi_1(S^4-K_i)$. It thus easily follows from Lemma~\ref{sec:genus:sum:lem:torsion-free} that the homology groups of $S^4(K_i)$ are isomorphic to those of $S^1\times S^3$. Furthermore, we can easily show that the cohomology ring of $S^4(K_i)$ is isomorphic to that of $S^1\times S^3$ (cf.\ \cite[Section~11]{Br97}). Note that every $S^2$-knot in $S^4$ admits a Seifert hypersurface (see e.g.\ \cite[Section 2.5]{CKS}). 
	
	Now let $X$ be a closed connected 4-manifold, and let $Y_i$ be the connected sum $X\#S^4(K_i)$. Since each $\pi_1(Y_i)$ is isomorphic to the free product $\pi_1(X)*\pi_1(S^4(K_i))$, it follows from the Grushko decomposition theorem (see e.g.\ \cite{Sta}) that the fundamental groups $\pi_1(Y_i)$ are pairwise non-isomorphic, showing that the 4-manifolds $Y_i$ are pairwise  non-homeomorphic. On the other hand, it is clear that their  cohomology rings  are isomorphic to that of $X\#S^1\times S^3$. Furthermore, by Theorem~\ref{sec:genus:sum:thm:genus preserving}, their genus functions are algebraically equivalent to that of $X$ and hence to that of $X\#S^1\times S^3$. The claim thus follows.  
\end{proof}
\section{Exotic 4-manifolds and genus functions}\label{sec:exotic}
In this section, we produce various families of exotic 4-manifolds with respect to genus functions. Note that we follow  Convention~\ref{sec:intro:convention}. 

We first show that any $4$-dimensional $2$-handlebody can be modified into infinite order generalized corks. Specifically, we prove the following theorem which is a refinement of Theorem~\ref{sec:cork:thm:2-handlebody:cork}.  

\begin{theorem}\label{sec:exotic:thm:2-handlebody:cork}
	For any $4$-dimensional $2$-handlebody $X$, there exists a generalized cork $(X',f)$ of infinite order such that $X'$ is HIHC-equivalent to  $X$, admits an embedding into $X$, and is homeomorphic to a $4$-manifold that deformation retracts onto an embedded copy of $X$.  Moreover, in the case where $H_2(X)\neq 0$,  there exist arbitrarily many  generalized corks $(X_i,f_i)$ of infinite order with the properties of  $(X',f)$ such that the $4$-manifolds $X_i$ are pairwise exotic. Furthermore, their genus functions are pairwise algebraically inequivalent for any orientations.  
\end{theorem}

\begin{proof}Let $X$ be a 4-dimensional 2-handlebody. By Proposition~\ref{sec:cork:prop:2-handlebody:Stein}, we have a 4-manifold $Y$ admitting a Stein structure such that $Y$ is HIHC-equivalent to  $X$, admits an embedding into $X$, and is homeomorphic to a 4-manifold that deformation retracts onto an embedded copy of $X$. Therefore, by Proposition~\ref{sec:cork:prop:construction}, we obtain an infinite order generalized cork $(X',f)$ having the required properties. 
	
	Now we assume $H_2(X)\neq 0$. Then, by Theorem~\ref{sec:cork:exotic Stein},  we obtain arbitrarily many pairwise exotic 4-manifolds $Y_i$ such that each $Y_i$ have the above properties of $Y$ and that the genus functions of $Y_i$ are pairwise algebraically inequivalent for any orientations. Therefore, by the construction in the proof of Proposition~\ref{sec:cork:prop:construction}, we obtain arbitrarily many generalized corks $(X_i,f_i)$ of infinite order such that the 4-manifolds $X_i$ have the above properties of $X'$ and are pairwise homeomorphic. Since each $X_i$ is obtained by attaching an invertible homology cobordism to $Y_i$, the genus function of each $X_i$ is algebraically equivalent to that of $Y_i$ due to Corollary~\ref{sec:genus:quasi-invertible:cor:inclusion:invertible homology cobordism}, and thus the 4-manifolds  $X_i$ are pairwise exotic. This completes the proof. 
\end{proof}
\begin{remark}\label{sec:exotic:thm:2-handlebody:cork:remark}
	As seen from the construction, the self-homeomorphism $f$ of $\partial X'$ extends to a self-homeomorphism of $X'$ that induces the identity on $H_2(X')$, and the same claim holds for the self-homeomorphism $f_i$ of each $\partial X_i$. 
\end{remark}

We here give a method for producing exotic 4-manifolds with equivalent genus functions from generalized corks. 

\begin{proposition}\label{sec:exotic:prop:cork:genus:method}
	Let $\{(C,f)\}_{f\in A}$ be a family of pairwise inequivalent generalized corks, and assume each $f\in A$ extends to a self-homeomorphism of $C$ that preserves the genus function. Then, there exists a family $\{C_f\}_{f\in A\sqcup \{\mathrm{id}_{\partial C}\}}$ of pairwise exotic $4$-manifolds HIHC-equivalent to  $C$ such that their genus functions are pairwise orientedly equivalent and are algebraically equivalent to that of $C$.  Furthermore, they admit deformation retractions onto an embedded copy of $C$ and embeddings into $C$.  
\end{proposition}
\begin{proof}For each $f\in A\sqcup \{\mathrm{id}_{\partial C}\}$, let $C_f$ denote the 4-manifold $C\cup_f \mathbb{P}$ given in the proof of Corollary~\ref{cor:AR:arbitrary}. By Corollary~\ref{cor:AR:arbitrary} and Lemma~\ref{lem:AR:homotopy}, the family $\{C_f\}_{f\in A\sqcup \{\mathrm{id}_{\partial C}\}}$ has all the desired properties possibly except those on genus functions.  By Corollary~\ref{sec:genus:quasi-invertible:cor:inclusion:invertible homology cobordism}, the genus function of each $C_f$ is algebraically equivalent to that of $C$, and this equivalence is induced by the inclusion $C\hookrightarrow C_f$. Since each $f\in A$ preserves the orientation and extends to a self-homeomorphism of $C$ that preserves the genus function, this self-homeomorphism of $C$ extends to an orientation-preserving homeomorphism $C_{\mathrm{id}_{\partial C}}\to C_f$  that preserves the genus functions. Therefore, the family has the desired properties on genus functions as well. 
\end{proof}

Now we are ready to produce various infinite families of  exotic 4-manifolds with pairwise equivalent genus functions. We here prove the following theorem which is a refinement of Theorem~\ref{sec:intro:thm:2-handlebody:infinite}. 
\begin{theorem}\label{sec:exotic:thm:2-handlebody:infinite}
	For any $4$-dimensional $2$-handlebody $X$, there exist infinitely many pairwise exotic $4$-manifolds HIHC-equivalent to  $X$ such that their genus functions are pairwise orientedly equivalent. Furthermore, these $4$-manifolds admit embeddings into $X$ and are homeomorphic to a $4$-manifold that deformation retracts onto an embedded copy of $X$.
\end{theorem}
\begin{proof}This is straightforward from Theorem~\ref{sec:exotic:thm:2-handlebody:cork}, Remark~\ref{sec:exotic:thm:2-handlebody:cork:remark}, and Proposition~\ref{sec:exotic:prop:cork:genus:method}. 
\end{proof}

The theorem below is a refinement of Theorem~\ref{sec:intro:thm:symplectic:infinitely}. 
\begin{theorem}\label{sec:exotic:thm:symplectic:infinitely}
For any $4$-manifold $X$ admitting an embedding into a symplectic $4$-manifold with weakly convex boundary, there exist infinitely many pairwise exotic $4$-manifolds HIHC-equivalent to  $X$ such that their genus functions are pairwise orientedly equivalent and are algebraically equivalent to that of $X$. Furthermore, these $4$-manifolds admit deformation retractions onto embedded copies of $X$ and embeddings into $X$.
\end{theorem}
\begin{proof}This is straightforward from Propositions~\ref{sec:cork:prop:construction} and \ref{sec:exotic:prop:cork:genus:method}. 
\end{proof}

The above theorem implies the following corollary, which is a slight refinement of Corollary~\ref{sec:intro:cor:nCP2}. 
\begin{corollary}\label{sec:exotic:cor:nCP2} $(1)$ For each positive integer $n$, there exist infinitely many pairwise exotic simply connected $4$-manifolds  $($with boundary$)$ such that their genus functions are pairwise orientedly equivalent and are algebraically equivalent to that of $n\mathbb{C}P^2$.  
	
	$(2)$ For any $D^2$-bundle over a closed surface, there exist infinitely many pairwise exotic $4$-manifolds HIHC-equivalent to  the total space $X$ of the bundle such that their genus functions are pairwise orientedly equivalent and are algebraically equivalent to that of $X$. 
\end{corollary}
\begin{proof}(1) Since $D^4$ admits a Stein structure, its blow-ups $n\overline{\mathbb{C}P^2}-\textnormal{int} D^4$ is a symplectic $4$-manifold with weakly convex boundary. Also, the genus function of $n\overline{\mathbb{C}P^2}-\textnormal{int} D^4$ is equivalent to that of $n\overline{\mathbb{C}P^2}$ due to Lemma~\ref{sec:genus:quasi-invertible:lem:genus:3-handle}. The claim (1) thus follows from Theorem~\ref{sec:exotic:thm:symplectic:infinitely}. 
	
	(2) Let $X$ denote the total space of the $D^2$-bundle over a closed surface of genus $g$ with Euler number $n$. Then $X$ admits a Stein structure for $n\leq 2g-2$ (see \cite[Exercises~11.2.5.(a)]{GS}). Hence the claim (2) follows from Theorem~\ref{sec:exotic:thm:symplectic:infinitely} in the case where $n\leq 2g-2$. In the case where $n\geq 2-2g$, the claim (2) follows from the case where $n\leq 2g-2$ by reversing the orientation. The remaining cases are $(g,n)=(0,0), (0,\pm1)$. In the case where $(g,n)=(0,\pm 1)$, the claim (2) follows from (1). In the case where $(g,n)=(0,0)$, the 4-manifold $X$ is diffeomorphic to $S^2\times D^2$ and hence admits an embedding into $D^4$. The claim (2) thus follows from Theorem~\ref{sec:exotic:thm:symplectic:infinitely} in this case. This completes the proof.  
\end{proof}
	
The theorem below is a slight refinement of Theorem~\ref{sec:intro:thm:Stein nuclei:infinitely:infinitely}. 

\begin{theorem}\label{sec:exotic:thm:Stein nuclei:infinitely:infinitely}
	There exist infinitely many pairwise exotic simply connected $4$-manifolds whose  genus functions are pairwise algebraically inequivalent for any orientations such that, for each $4$-manifold $X$ among them, there exist infinitely many pairwise exotic $4$-manifolds homeomorphic to $X$ whose genus functions are all orientedly equivalent to that of $X$. 
\end{theorem}
\begin{proof}Let $\{X_i\}_{i\in \mathbb{Z}}$ be an infinite family of pairwise exotic 4-manifolds admitting Stein structures such that the genus functions of the members are algebraically inequivalent for any orientations. Such examples were obtained by Akbulut and the author in \cite{AY_JSG14} and the author in \cite{Y14}. By the construction in the proof of  Proposition~\ref{sec:cork:prop:construction}, we obtain an infinite family $\{(Y_i,f_i)\}_{i\in \mathbb{Z}}$ of infinite order generalized corks such that each $Y_i$ is obtained from $X_i$ by attaching the same invertible homology cobordism, and that the 4-manifolds $Y_i$ are pairwise homeomorphic. By Corollary~\ref{sec:genus:quasi-invertible:cor:inclusion:invertible homology cobordism}, the genus function of each $Y_i$ is algebraically equivalent to that of $X_i$. Due to the proof of Proposition~\ref{sec:exotic:prop:cork:genus:method}, by attaching another fixed invertible homology cobordism to each $Y_i$ and twisting $Y_i$ in the resulting 4-manifold $X'_{i,0}$ using each $f_i^j$, we obtain an infinite family $\{X'_{i,j}\}_{(i,j)\in \mathbb{Z}\times \mathbb{Z}}$ of pairwise exotic simply connected 4-manifolds. By Proposition~\ref{sec:exotic:prop:cork:genus:method}, for each fixed $i$, the genus functions of the 4-manifolds $X'_{i,j}$ are pairwise equivalent, and furthermore for each fixed $j$, the genus functions of $X'_{i,j}$ are pairwise algebraically inequivalent for any orientations.  	 
\end{proof}	

By the proof of the above theorem, we can construct similar infinite families for varieties of homeomorphism types by using exotic Stein manifolds of the author~\cite{Y14} as a building block. For example, any finitely presented group is realized as the fundamental group of infinitely many exotic 4-manifolds satisfying the properties in the above theorem. 

Here we construct similar examples for much more homeomorphism types, by relaxing the condition ``infinitely many'' to ``arbitrarily many''.  Specifically, we prove the following theorem which is a refinement of Theorem~\ref{sec:intro:thm:2-handlebody:infinitely:arbitrarily}.  

\begin{theorem}\label{sec:exotic genus:thm:2-handlebody:infinitely:arbitrarily}
For any $4$-dimensional $2$-handlebody $X$ with $H_2(X)\neq 0$, there exist arbitrarily many pairwise exotic $4$-manifolds HIHC-equivalent to  $X$ whose genus functions are pairwise algebraically inequivalent for any orientations, such that each $4$-manifold $X'$ among them satisfies the following.
\begin{itemize}
	\item There exist infinitely many pairwise exotic $4$-manifolds homeomorphic to $X'$ whose genus functions are all orientedly equivalent to that of $X'$. 
	\item $X'$ and the above infinitely many $4$-manifolds admit embeddings into $X$ and are homeomorphic to a  $4$-manifold that deformation retracts onto an embedded copy of $X$. 
\end{itemize}
\end{theorem}
\begin{proof} Using infinite order generalized corks of Theorem~\ref{sec:exotic:thm:2-handlebody:cork} as a building block, we can prove the claim similarly to the proof of Theorem~\ref{sec:exotic:thm:Stein nuclei:infinitely:infinitely}. 
\end{proof}

As an application of our results and Theorem~\ref{sec:cork:exotic Stein}, we give exotic 4-manifolds which are stably exotic under boundary sums with arbitrary Stein manifolds having non-degenerate intersection forms.  Specifically, we prove the following theorem which is a refinement of Theorem~\ref{sec:intro:thm:sum:Stein}. 
\begin{theorem}\label{sec:exotic:thm:sum:Stein}
	Let $X$ be a $4$-dimensional $2$-handlebody with $H_2(X)\neq 0$ whose intersection form is represented by a zero-matrix. Then	for any positive integer $n$, there exist pairwise exotic $4$-manifolds $X_1,X_2,\dots, X_n$ HIHC-equivalent to  $X$ and admitting Stein structures such that, for any pairwise orientedly homeomorphic $4$-manifolds $Y_1,Y_2,\dots, Y_n$ having non-degenerate intersection forms and admitting Stein structures, the boundary sums $X_1\natural Y_1$, $X_2\natural Y_2$, $\dots$, $X_n\natural Y_n$ are  pairwise exotic. 
\end{theorem}
\begin{proof}By Theorem~\ref{sec:cork:exotic Stein}, we have pairwise exotic $4$-manifolds $X_1,\dots, X_n$ HIHC-equivalent to  $X$ and admitting Stein structures such that their genus functions are pairwise algebraically inequivalent for any orientations. Let $Y_1,\dots, Y_n$ be pairwise orientedly homeomorphic 4-manifolds having non-degenerate intersection forms and admitting Stein structures. Since each $X_i$ and $Y_i$ admit Stein structures, we can apply the argument on the genus functions of the 4-manifolds $X_i$ in the proof of Theorem~\ref{sec:cork:exotic Stein} to the restrictions of those of  the 4-manifolds  $X_i\natural Y_i$ to $H_2(X_i)$ in the exactly same way. We thus see that the restrictions of the genus functions of  $X_i\natural Y_i$ to $H_2(X_i)$ remain pairwise algebraically inequivalent for any orientations. Since  the second homology group of any 2-handlebody is torsion-free, it follows from Proposition~\ref{sec:bilinear:prop:genus:zero-part} that the genus functions of $X_i\natural Y_i$ are pairwise algebraically inequivalent for any orientations. Therefore,  the 4-manifolds $X_i\natural Y_i$ are pairwise exotic. 
\end{proof}
There are rich examples of 4-manifolds satisfying the above assumption of $X$, hence we can realize various homeomorphism types including non-simply connected ones as those of the above $X_1,\dots, X_n$.

\section{Exotically knotted submanifolds}\label{sec:knotted}
As an application of our exotic 4-manifolds in the last section, we modify a submanifold of a 4-manifold with codimension at most one into infinitely many pairwise exotically knotted submanifolds. Note that we follow Convention~\ref{sec:intro:convention}. 

For a codimension-$0$ submanifold $X$ of a 4-manifold $Z$,  the closure of $Z-X$ is called the \textit{exterior} of $X$ and denoted by $E(X)$. 
To state our results shortly, we use the following terminologies. 

\begin{definition}\label{sec:knotted:def:exotically embedded}Let $X_1,X_2$ be (possibly non-orientable) submanifolds of a (possibly non-orientable) 4-manifold $Z$. 
	
			$(1)$ The pair $X_1, X_2$ is called \textit{exotically knotted} if they have the following properties. 
	\begin{itemize}
				\item [\textnormal{(i)}] $X_1$ is diffeomorphic to $X_2$. In the case where $Z$ is orientable, and their codimensions are zero, the diffeomorphism is assumed to preserve the orientations induced from $Z$. 
		\item [\textnormal{(ii)}] $X_1$ is topologically ambient isotopic to $X_2$. In the case where $\partial Z$ is non-empty, the topological isotopy is assumed to fix $\partial Z$ pointwise. 
		\item [\textnormal{(iii)}] $X_1$ is not smoothly isotopic to $X_2$. (We do not assume a smooth isotopy fixes $\partial Z$ pointwise.)
	\end{itemize}
	
	$(2)$ The pair $X_1, X_2$ is called \textit{exotically embedded} if they have the following properties. 
	\begin{itemize}
		\item [\textnormal{(i)}] The condition (i) of (1) holds.  
	\item [\textnormal{(ii)}] There exists a self-homeomorphism of $Z$ that maps $X_1$ to $X_2$. In the case where $Z$ is orientable, the self-homeomorphism is assumed to be orientation-preserving. 
	\item [\textnormal{(iii)}] There exists no self-diffeomorphism of $Z$ that maps $X_1$ to $X_2$. 
	\end{itemize}
\end{definition}
In the above definition, we require the conditions on orientations and the boundary so that our examples have stronger properties, but our results are still new without these extra conditions. Also, our results are new as exotically knotted submanifolds and also as exotically embedded submanifolds. 
 
We first produce exotically knotted submanifolds of codimension-0. Specifically, we prove the following theorem which is a refinement of Theorem~\ref{sec:intro:thm:exotic:embedded}. 

\begin{theorem}\label{sec:knotted:thm:exotic:embedded}
Let $X$ be a codimension-$0$ submanifold of a  $4$-manifold $Z$.  These manifolds may be non-orientable and may have $($possibly disconnected$)$ boundaries. Suppose the exterior of $X$  either is a  $($orientable$)$ $2$-handlebody or admits an embedding into a symplectic $4$-manifold with weakly convex boundary, and assume the exterior has connected boundary. Then,  $Z$ admits infinitely many pairwise exotically knotted and exotically embedded codimension-$0$ submanifolds HIHC-equivalent to  $X$. Furthermore, each of these submanifolds contains $X$ as a submanifold and admits a deformation retraction onto $X$.
\end{theorem}

To prove this theorem, we give a method for producing exotically knotted submanifolds by using generalized corks obtained from contractible corks.  We use the following setting only for the forthcoming Lemma~\ref{sec:knotted:lem:twisted_double} and Proposition~\ref{sec:knotted:prop:exotic:embedded:method}. Let $\{(C_0,f)\}_{f\in A}$ be a family of pairwise inequivalent contractible corks having a fixed embedded 3-ball $B_0$ in $\partial C_0$ such that every $f\in A$ is identity on a neighborhood of $B_0$. Let $Y$ be a 4-manifold with non-empty boundary $M$. Take the boundary sum $(I\times M)\natural C_0$ for the boundary component $1\times M$ using the 3-ball $B_0$ of $\partial C_0$. This manifold $P$ is a homology cobordism from $M$ to $M\#\partial C_0$. We set $M'=M\#\partial C_0$. 
For each $f\in A\sqcup\{\mathrm{id}_{\partial C_0}\}$, let $\widetilde{f}$ denote the self-diffeomorphism of $M'$ that extends the self-diffeomorphism $f$ of $\partial C_0$ as the identity on the rest of $M'$. We denote the glued 4-manifold $Y\cup_{\mathrm{id}_{M}} P$ by $C$, and assume that $\{(C,\widetilde{f})\}_{f\in {A}}$ is a family of pairwise inequivalent generalized corks. We furthermore assume that, for each $f\in A\sqcup\{\mathrm{id}_{\partial C_0}\}$, the twisted double $C_0\cup_{{f}} \overline{C_0}$ of $C_0$ is diffeomorphic to $S^4$. 

Let  $D_{C_0}$ be the ``partial double'' $C_0\cup_{\mathrm{id}_{\partial C_0-\mathrm{int}\, B_0}} \overline{C_0}$ of $C_0$,  which is obtained from $C_0$ by attaching $\overline{C_0}$ using the identity map on the exterior of $B_0$. Clearly,  $D_{C_0}$ is a submanifold of  the double $P\cup_{\mathrm{id}_{M'}} \overline{P}$. Due to the assumption on the double of $C_0$, the disk theorem implies that $D_{C_0}$ is diffeomorphic to the 4-ball.  
Since each $f\in A\sqcup\{\mathrm{id}_{\partial C_0}\}$ is identity on a neighborhood of the 3-ball $B_0$ in $\partial C_0$, it follows that  the twisted double $P\cup_{\widetilde{f}} \overline{P}$ is obtained from $(P\cup_{\mathrm{id}_{M'}} \overline{P})-\textnormal{int}\, D_{C_0}$ by attaching the ``twisted partial double'' $C_0\cup_{f|_{\partial C_0-\mathrm{int}\, B_0}} \overline{C_0}$ of $C_0$, which is obtained from $C_0$ by attaching $\overline{C_0}$ using the restriction of $f$ to the exterior of $B_0$. We note that each twisted double $P\cup_{\widetilde{f}} \overline{P}$ is diffeomorphic to the connected sum of $I\times M$ and the twisted double $C_0\cup_f \overline{C_0}$. 

\begin{lemma}\label{sec:knotted:lem:twisted_double}
	$(1)$ $P$ admits a deformation retraction onto the boundary component $M$. 
	
	$(2)$ For each $f\in A\sqcup\{\mathrm{id}_{\partial C_0}\}$, there exists a diffeomorphism $P\cup_{\widetilde{f}} \overline{P} \to I\times M$ such that it fixes the boundary $M\sqcup \overline{M}$ pointwise and that its restriction to the submanifold $(P\cup_{\mathrm{id}_{M'}} \overline{P})-\textnormal{int}\, D_{C_0}$ is independent of $f$. 
\end{lemma}
\begin{proof}The claim $(1)$ is straightforward from the construction of $P$. Since each twisted double $C_0\cup_f \overline{C_0}$ is diffeomorphic to $S^4$, the disk theorem shows that the twisted partial double of $C_0$ is diffeomorphic to $D^4$. Thus, every diffeomorphism between the boundary of the twisted partial double and $\partial D^4$ extends to a diffeomorphism between their interiors. Therefore, the claim (2) follows from the aforementioned decomposition of the twisted double of $P$.   
\end{proof}
The following proposition gives exotically knotted codimension-$0$ submanifolds. 
\begin{proposition}\label{sec:knotted:prop:exotic:embedded:method}
 Let $X$ be a codimension-$0$ submanifold of a  $4$-manifold $Z$. These manifolds may be non-orientable and may have $($possibly disconnected$)$ boundaries. If the exterior  $E(X)$ is  diffeomorphic to the $4$-manifold $Y$, then $Z$ admits a family $\{X_{
 \widetilde{f}}\}_{f\in A\sqcup\{\mathrm{id}_{\partial C_0}\}}$ of pairwise exotically knotted and exotically embedded codimension-$0$ submanifolds HIHC-equivalent to  $X$. Furthermore, each $X_{\widetilde{f}}$ contains $X$ as a submanifold and admits a deformation retraction  onto $X$. 
 \end{proposition}
 \begin{proof}
 	We identify the boundary $\partial E(X)$ with $M$ using a fixed diffeomorphism $Y\to E(X)$. Note that $Y$ is oriented due to our convention. Let $\mathbb{P}$ be the invertible homology cobordism from $M'$ to some closed (connected) 3-manifold $N$ used in Theorem~\ref{thm:AR:arbitrary}. We denote the 4-manifold $C\cup_{\mathrm{id}_{M'}} \mathbb{P}$ by $S_{\mathrm{id}_{M'}}$. We attach the upside down cobordism $\overline{\mathbb{P}}$ to $S_{\mathrm{id}_{M'}}$ along $N$ using $\mathrm{id}_{N}$, and then attach the upside down cobordism $\overline{P}$ of $P$ along $M'$ using $\mathrm{id}_{M'}$. We denote the resulting 4-manifold  $Y\cup_{\mathrm{id}_{M}} P\cup_{\mathrm{id}_{M'}} \mathbb{P}\cup_{\mathrm{id}_N} \overline{\mathbb{P}}\cup_{\mathrm{id}_{M'}} \overline{P}$ by $\widehat{Y}_{\mathrm{id}_{M'}}$. Since the double $\mathbb{P}\cup_{\mathrm{id}_N} \overline{\mathbb{P}}$  admits a diffeomorphism to $I\times M'$ that fixes the boundary $\overline{M'}\sqcup M'$ pointwise (see Remark~\ref{rem:AR:cobordism}), the submanifold  $P\cup_{\mathrm{id}_{M'}} \mathbb{P}\cup_{\mathrm{id}_N} \overline{\mathbb{P}}\cup_{\mathrm{id}_{M'}} \overline{P}$ of $\widehat{Y}_{\mathrm{id}_{M'}}$ admits a diffeomorphism to the double $P\cup _{\mathrm{id}_{M'}}\overline{P}$ that fixes the boundary $\overline{M}\sqcup M$ pointwise. Let $\widehat{D}_{C_0}$ denote the preimage of $D_{C_0}$ under this diffeomorphism.  
By Lemma~\ref{sec:knotted:lem:twisted_double}, the 4-manifold $\widehat{Y}_{\mathrm{id}_{M'}}$ is diffeomorphic to $Y$ and hence to the exterior $E(X)$, and this diffeomorphism $\varphi_{\mathrm{id}_{M'}}: \widehat{Y}_{\mathrm{id}_{M'}}\to E(X)$ fixes the boundary $M$ pointwise. We set $D_{C_0}^0=\varphi_{\mathrm{id}_{M'}}(\widehat{D}_{C_0})$. 
 	
 	For each $f\in A$,  we apply generalized cork twists to $S_{\mathrm{id}_{M'}}$ and $\widehat{Y}_{\mathrm{id}_{M'}}$ along $(C,	\widetilde{f})$, and denote the resulting 4-manifolds $C\cup_{\widetilde{f}} \mathbb{P}$ and  $C\cup_{\widetilde{f}} \mathbb{P}\cup_{\mathrm{id}_N} \overline{\mathbb{P}}\cup_{\mathrm{id}_{M'}} \overline{P}$ by $S_{\widetilde{f}}$ and $\widehat{Y}_{\widetilde{f}}$, respectively. 
 	Due to the diffeomorphism $P\cup_{\mathrm{id}_{M'}} \mathbb{P}\cup_{\mathrm{id}_N} \overline{\mathbb{P}}\cup_{\mathrm{id}_{M'}} \overline{P}\to P\cup _{\mathrm{id}_{M'}}\overline{P}$ in the last paragraph, each cobordism $P\cup_{\widetilde{f}} \mathbb{P}\cup_{\mathrm{id}_N} \overline{\mathbb{P}}\cup_{\mathrm{id}_{M'}} \overline{P}$ is obtained from the exterior of $\widehat{D}_{C_0}$ in $P\cup_{\mathrm{id}_{M'}} \mathbb{P}\cup_{\mathrm{id}_N} \overline{\mathbb{P}}\cup_{\mathrm{id}_{M'}} \overline{P}$ by attaching a 4-manifold diffeomorphic to the twisted partial double of $C_0$. 
 	Thus, by  Lemma~\ref{sec:knotted:lem:twisted_double}, 
 	we have a diffeomorphism $P\cup_{\widetilde{f}} \mathbb{P}\cup_{\mathrm{id}_N} \overline{\mathbb{P}}\cup_{\mathrm{id}_{M'}} \overline{P}\to P\cup_{{\mathrm{id}_{M'}}} \mathbb{P}\cup_{\mathrm{id}_N} \overline{\mathbb{P}}\cup_{\mathrm{id}_{M'}} \overline{P}$ that fixes the exterior of $\widehat{D}_{C_0}$ and hence the boundary $\overline{M}\sqcup M$ pointwise. 
 	Extending this diffeomorphism as the identity on the rest of $\widehat{Y}_{\widetilde{f}}$, we obtain a diffeomorphism $\psi_{\widetilde{f}}:\widehat{Y}_{\widetilde{f}}	\to \widehat{Y}_{\mathrm{id}_{M'}}$ that is identity on $\widehat{Y}_{\mathrm{id}_{M'}}-\textnormal{int}\, \widehat{D}_{C_0}$. We define a diffeomorphism   $\varphi_{{\widetilde{f}}}:\widehat{Y}_{\widetilde{f}}\to E(X)$ by $\varphi_{\widetilde{f}}=\varphi_{\mathrm{id}_{M'}}\circ \psi_{\widetilde{f}}$.

 	For each  $f\in A\sqcup\{{\mathrm{id}_{\partial C_0}}\}$, we denote the submanifolds $\varphi_{\widetilde{f}}(S_{\widetilde{f}})$,  $\varphi_{\widetilde{f}}(\overline{\mathbb{P}})$ and $\varphi_{\widetilde{f}}(\overline{P})$  of $E(X)$ respectively by $S^0_{\widetilde{f}}$,  $\overline{\mathbb{P}}^0_{\widetilde{f}}$ and $\overline{P}^0_{\widetilde{f}}$. 
 	Let $X_{\widetilde{f}}$ be the (possibly non-orientable) submanifold $\overline{\mathbb{P}}^0_{\widetilde{f}}	\cup \overline{P}^0_{\widetilde{f}} \cup X$ of $Z$. 
 	Then,  the 4-manifolds $X_{\widetilde{f}}$ are pairwise diffeomorphic, since each $\varphi_{\widetilde{f}}$ fixes the boundary $M$ of $\widehat{Y}_{\widetilde{f}}$ pointwise. Note that $X_{\widetilde{f}}$ are pairwise orientedly diffeomorphic in the case where $Z$ is orientable. Clearly, the submanifold $\overline{\mathbb{P}}^0_{\widetilde{f}}\cup \overline{P}^0_{\widetilde{f}}$ of $Z$ is an invertible homology cobordism  from $\overline{M}$ to $\overline{N}$, and each $X_{\widetilde{f}}$ is obtained from $X$ by attaching this cobordism  along $\overline{M}$. We thus see that each $X_{\widetilde{f}}$ deformation retracts onto $X$ and is HIHC-equivalent to  $X$ by Lemmas~\ref{lem:AR:homotopy} and \ref{sec:knotted:lem:twisted_double}.  
 	
 	We will show that the submanifolds $X_{\widetilde{f}}$ of $Z$ are pairwise exotically knotted and exotically embedded. We note that $Z=S^0_{{\widetilde{f}}}\cup X_{{\widetilde{f}}}$ for each $f\in A\sqcup\{{\mathrm{id}_{\partial C_0}}\}$. By Theorem~\ref{thm:AR:arbitrary}, the exteriors of  $X_{\widetilde{f}}$ are pairwise exotic, and thus the pairs $(Z, X_{\widetilde{f}})$ are pairwise non-diffeomorphic. Consequently, the submanifolds $X_{\widetilde{f}}$ are pairwise smoothly non-isotopic. 
 	
 	In the rest, we will prove that the submanifolds $X_{\widetilde{f}}$ are pairwise topologically ambient isotopic in $Z$ by showing that each $X_{\widetilde{f}}$ is topologically ambient isotopic to $X_{{\mathrm{id}_{M'}}}$.  
Due to the assumption, each ${\widetilde{f}}$ extends to a self-homeomorphism of $C$ that is identity on the exterior of the contractible 4-manifold $C_0$. Extending this homeomorphism as the identity on the rest of $\widehat{Y}_{\widetilde{f}}$, we get a homeomorphism $\widehat{Y}_{\widetilde{f}}\to \widehat{Y}_{\mathrm{id}_{M'}}$ that maps $S_{\widetilde{f}}$ to $S_{\mathrm{id}_{M'}}$. Note that this homeomorphism is identity on the exterior $\widehat{Y}_{\mathrm{id}_{M'}}-\textnormal{int}\, \widehat{D}_{C_0}$ of $\widehat{D}_{C_0}$. Composing this homeomorphism with the diffeomorphisms $\varphi_{\widetilde{f}}^{-1}$ and $\varphi_{\mathrm{id}_{M'}}$, we obtain a self-homeomorphism of $E(X)$ that maps $S^0_{\widetilde{f}}$ to $S^0_{{\mathrm{id}_{M'}}}$ and fixes $E(X)-\textnormal{int}\, D_{C_0}^0$ and hence the boundary $M$ of $E(X)$ pointwise. Extending this self-homeomorphism as the identity on the rest of $Z$, we obtain a self-homeomorphism ${\Phi}_{\widetilde{f}}$ of $Z$ that maps $X_{\widetilde{f}}$ to $X_{{\mathrm{id}_{M'}}}$ and fixes the exterior of $D_{C_0}^0$ pointwise. 
Note that ${\Phi}_{\widetilde{f}}$ fixes $\partial Z$ pointwise in the case where $\partial Z$ is non-empty. 

Clearly, the restriction of ${\Phi}_{\widetilde{f}}$ to $D_{C_0}^0$ fixes $\partial D_{C_0}^0$ pointwise. Since $D_{C_0}^0$ is diffeomorphic to $D^4$, the Alexander trick~\cite{Alex} shows that the restriction to $D_{C_0}^0$ is topologically isotopic to the identity map on $D_{C_0}^0$ fixing $\partial D_{C_0}^0$ pointwise. Thus, ${\Phi}_{\widetilde{f}}$ is  topologically isotopic to the identity map on $Z$. Note that this isotopy fixes $\partial Z$ pointwise in the case where $\partial Z$ is non-empty. 
Therefore, each submanifold $X_{\widetilde{f}}$ of $Z$ is topologically ambient isotopic to $X_{{\mathrm{id}_{M'}}}$. This completes the proof. 
\end{proof}

We can now easily prove Theorem~$\ref{sec:knotted:thm:exotic:embedded}$ using our generalized corks. 

\begin{proof}[Proof of Theorem~$\ref{sec:knotted:thm:exotic:embedded}$]Let  $Y$ be a 4-manifold diffeomorphic to the exterior $E(X)$. We first consider the case where $Y$ admits an embedding into a  symplectic $4$-manifold with weakly convex boundary.  Let $(C_0,\tau)$ be an infinite order contractible cork of Gompf (see Subsection~\ref{subsec:Gompf's corks}). By results of Akbulut~\cite{A20}, Tange~\cite{T21} and Naylor and Schwartz \cite{NS22},  the twisted double $C_0\cup_{\tau^i} \overline{C_0}$ is diffeomorphic to $S^4$ for each integer $i$. 
	Therefore, we can apply Proposition~\ref{sec:knotted:prop:exotic:embedded:method} using a generalized cork obtained from $Y$ in the proof of Proposition~\ref{sec:cork:prop:construction}. The claim thus follows in this case. 
	
	We next consider the case where $Y$ is a 2-handlebody. In this case, we modify $Y$ into a Stein handlebody $Y'$ by using Proposition~\ref{sec:cork:prop:2-handlebody:Stein}. Note that $Y$ is obtained from $Y'$ by attaching the upside down cobordism of an invertible cobordism $P$. Furthermore, $P$ admits a deformation retraction onto the negative boundary $\partial Y$. Since $Y$ is diffeomorphic to $E(X)$, we may regard $P$ as a submanifold of $E(X)$ such that the negative boundary of $P$ is $\partial E(X)$. Then, the union $X'=X\cup P$ is a  submanifold of $Z$. Clearly, $X'$ contains $X$ as a submanifold and admits a deformation retraction onto $X$. Since the exterior of $X'$ in $Z$ is diffeomorphic to the Stein manifold $Y'$, we can prove the desired claim by applying that of the first case to $X'$. 
\end{proof}
\begin{remark}\label{sec:knotted:rem:exotic:embedded}
	For Theorem~\ref{sec:knotted:thm:exotic:embedded} and Proposition~\ref{sec:knotted:prop:exotic:embedded:method}, the exotically knotted submanifolds of $Z$   are obtained from $X$ by attaching invertible homology cobordisms embedded in $Z$, and hence they admit embeddings into $X$. It thus follows from Corollary~\ref{sec:genus:quasi-invertible:cor:inclusion:invertible homology cobordism} that the genus functions of these exotically knotted submanifolds are algebraically equivalent to that of $X$. 
\end{remark}
We next discuss codimension-1 submanifolds. We use the following terminology. 
\begin{definition}\label{sec:knotted:def:cobordant_through}Let $Z$ be a (possibly non-orientable) 4-manifold which may have  (possibly disconnected) boundary, and  let $M_1,M_2$ be codimension-1 submanifolds of $Z$.  The submanifold $M_1$ is called \textit{homology cobordant to  $M_2$ through $Z$} if $Z$ admits a codimension-$0$ (oriented) submanifold $P$ with $\partial P=\overline{M_1}\sqcup M_2$ such that $P$ is a homology cobordism from $M_1$ to $M_2$. 
\end{definition}

As an application of Theorem~$\ref{sec:knotted:thm:exotic:embedded}$, we modify separating codimension-$1$ submanifolds into infinitely many pairwise exotically knotted and exotically embedded codimension-$1$ submanifolds. Specifically, we prove the following theorem which is a refinement of Theorem~\ref{sec:intro:thm:exotic:embedded:3-manifold}. 
\begin{theorem}\label{sec:knotted:theorem:exotic:embedded:3-manifold}Let $M$ be a closed $3$-manifold embedded in the interior of   a $($possibly non-orientable$)$ $4$-manifold $Z$ which may have $($possibly disconnected$)$ boundary. Suppose the exterior of the tubular neighborhood of $M$ has a connected component with connected boundary which either is a $4$-dimensional $($orientable$)$ $2$-handlebody or admits an embedding into a symplectic $4$-manifold with weakly convex boundary.   Then, $Z$ admits infinitely many pairwise exotically knotted and exotically embedded codimension-$1$ submanifolds which are homology cobordant to $M$ through $Z$.  
\end{theorem}
\begin{proof}Due to the assumption, we can apply Theorem~$\ref{sec:knotted:thm:exotic:embedded}$ to a component $X$ of the exterior of the tubular neighborhood for $M$. Note that $\partial E(X)$ is diffeomorphic to $M$. Let $X_i$ $(i\in \mathbb{Z})$ be the resulting infinitely many pairwise exotically knotted and exotically embedded codimension-0 (possibly non-orientable) submanifolds of $Z$, and let each $M_i$ be the boundary $\partial E(X_i)$ of the exterior $E(X_i)$. 
	
	We fix an orientation of the exterior $E(X)$ and equip each $E(X_i)\subset E(X)$ with the orientation induced from $E(X)$. In the case where the orientaion of $M$ is that of $\partial E(X)$ (resp.\ $\partial\overline{E(X)}$), we equip each $M_i$ with the orientaion of $\partial E(X_i)$ (resp.\ $\partial\overline{E(X_i)}$). Clearly, the 3-manifolds $M_i$ are pairwise orientedly diffeomorphic. 
	By the proof of  Theorem~$\ref{sec:knotted:thm:exotic:embedded}$, each $M_i$ is homology cobordant to $M$ through $Z$. Since the submanifolds $X_i$ are pairwise topologically ambient isotopic, their boundary components $M_i$ are also pairwise topologically ambient isotopic. 
	
	We will show that the pairs $(Z,M_i)$ are pairwise non-diffeomorphic. 
	By the proof of  Theorem~$\ref{sec:knotted:thm:exotic:embedded}$, the exteriors $E(X_i)$ are pairwise exotic, and the 4-manifolds $X_i$ are pairwise diffeomorphic. Note that the exterior of the tubular neighborhood of each $M_i$ is the disconnected 4-manifold $E(X_i)\sqcup X_i$. 
	Now suppose to the contrary that a pair $(Z,M_i)$ is diffeomorphic to some $(Z,M_j)$ with $i\neq j$. Then, $E(X_i)\sqcup X_i$ is diffeomorphic to $E(X_j)\sqcup X_j$ due to the uniqueness of the tubular neighborhood. Since $E(X_i)$ is not diffeomorphic to $E(X_j)$, we see that $E(X_i)$ is diffeomorphic to $X_j$ and hence to $X_i$ and that $X_i$ is diffeomorphic to $E(X_j)$. Consequently, $E(X_i)$ is diffeomorphic to $E(X_j)$, giving a contradiction. Thus, the pairs $(Z,M_i)$ are pairwise non-diffeomorphic. Therefore, the 3-manifolds $M_i$ are pairwise exotically embedded and exotically knotted in $Z$. 
\end{proof}

We obtain the following corollary which is a restatement of Corollary~\ref{sec:intro:cor:exotic:embedded:3-manifold:S4}. 
\begin{corollary}\label{sec:knotted:cor:exotic:embedded:3-manifold:S4}
For any  closed $3$-manifold $M$ embedded in $S^4$,  every $($possibly non-orientable$)$ $4$-manifold admits infinitely many pairwise exotically knotted and exotically embedded codimension-$1$ submanifolds which are homology cobordant to $M$. 
\end{corollary}
\begin{proof}Let $M$ be a closed 3-manifold embedded in $S^4$. It follows from the Mayer--Vietoris exact sequence that $M$ is separating. Let $Y$ be a connected component of the exterior of the tubular neighborhood for $M$. Since $M$ is separating, the boundary $\partial Y$ is diffeomorphic to $M$. Also, $Y$ admits an embedding into $D^4$ and hence into a symplectic 4-manifold with weakly convex boundary. Now let $Z$ be an arbitrary (possibly non-orientable) 4-manifold. Then $Y$ admits an embedding into the interior of $Z$, and thus we can apply Theorem~\ref{sec:knotted:theorem:exotic:embedded:3-manifold} to the submanifold $\partial Y$ of $Z$. Therefore, the claim follows.  
\end{proof}

The above corollary implies the following corollary which is a refinement of Corollary~\ref{sec:intro:cor:exotic:embedded:3-manifold:any}. 
\begin{corollary}\label{sec:knotted::cor:exotic:embedded:3-manifold:any}
			For each integer $n\geq 0$, there exists a closed $3$-manifold $M$ with $b_1=n$ such that every $($possibly non-orientable$)$ $4$-manifold admits infinitely many pairwise exotically knotted and exotically embedded submanifolds diffeomorphic to $M$. 
\end{corollary}
\begin{proof}As is well-known, for each integer $n\geq 0$, there are various closed 3-manifolds with $b_1=n$ admitting embeddings into $S^4$. For example, each $n(S^1\times S^2)$ is such a 3-manifold. The claim thus follows from Corollary~\ref{sec:knotted:cor:exotic:embedded:3-manifold:S4}. 
\end{proof}

For a given closed 3-manifold, here we consider the question of what is the minimum second Betti number of a closed 4-manifold which contains a pair of exotically knotted (or exotically embedded) separating submanifolds diffeomorphic to the 3-manifold. Since every closed 3-manifold embedded in a simply connected closed 4-manifold is separating as follows from the Mayer--Vietoris exact sequence, this question is an extension of the question mentioned in the paragraph preceding  Corollary~\ref{sec:intro:cor:exotic:embedded:3-manifold:some_4-manifold}, where the ambient 4-manifold is assumed to be simply connected. 
As an application of Theorem~\ref{sec:knotted:theorem:exotic:embedded:3-manifold}, we show that every non-negative integer is realized as such a minimum number, giving  a refinement of Corollary~\ref{sec:intro:cor:exotic:embedded:3-manifold:some_4-manifold}. 
\begin{corollary}\label{sec:knotted:cor:exotic:embedded:3-manifold:some_4-manifold}
	For each even integer $n\geq 0$ and each integer $m$ with $m\geq n/2$, there exists a closed $3$-manifold $M$ with $b_1=m$ such that a $($simply connected$)$ closed $4$-manifold with $b_2=n$ contains infinitely many pairwise exotically knotted and exotically embedded separating submanifolds diffeomorphic to $M$, but no closed $4$-manifold with $b_2 <n$ contains a separating submanifold diffeomorphic to $M$. 
\end{corollary}
\begin{proof}
	We set $r=n/2\geq 0$ and $k=m-r\geq 0$.  
	In the $n=0$ case, the claim follows from Corollary~\ref{sec:knotted::cor:exotic:embedded:3-manifold:any}, so we prove the $n>0$ case. Let $K$ be an invertible knot in $S^3$ with $|\sigma(K)|\geq 2n-1$, where $\sigma(K)$ denotes the signature of $K$. For example, many torus knots satisfy this condition. Let $Y$ be the 4-manifold obtained from $D^4$ by attaching a 2-handle along $K$ with $0$-framing, and  let $Y_{r,k}$ be the boundary sum of $r$ copies of  $Y$ and $k$ copies of $S^2\times D^2$. We denote by $M_{r,k}$ the boundary $\partial Y_{r,k}$.  
We can embed each 4-manifold  $Y_{r,k}$ and hence $M_{r,k}$ into the closed 4-manifold $r(S^2\times S^2)$ by taking the double of $Y_{r,0}$, which contains $Y_{r,k}$ as a submanifold. 

By Theorem~\ref{sec:knotted:theorem:exotic:embedded:3-manifold}, we obtain a closed 3-manifold $M'_{r,k}$ homology cobordant to $M_{r,k}$ such that $r(S^2\times S^2)$ admits infinitely many pairwise exotically knotted and exotically embedded submanifolds diffeomorphic to $M'_{r,k}$. As easily seen, we have $b_1(M'_{r,k})=b_1(M_{r,k})=r+k=m$. Note that $r(S^2\times S^2)$ is simply connected and satisfies $b_2=2r=n$.  
	
	Since $M'_{r,k}$ is homology cobordant to $M_{r,k}$, we have a homology cobordism $P_{r,k}$ from $M'_{r,k}$ to $M_{r,k}=M_{r,0}\#k(S^2\times S^1)$. By attaching 3-handles to $P_{r,k}$, we obtain a cobordism from $M'_{r,k}$ to $M_{r,0}$. By further attaching 1-handles to this cobordism and then attaching the upside down cobordism of $P_{r,k}$, we obtain a cobordism $P'_{r,k}$ from $M'_{r,k}$ to itself whose interior contains $M_{r,0}$ as a separating submanifold. 
	
	Now we show that each $M'_{r,k}$ admits no embedding into a closed 4-manifold with $b_2<n$ as a separating submanifold. 
	Suppose to the contrary that $M'_{r,k}$ is diffeomorphic to a separating submanifold of a closed 4-manifold with $b_2<n$. We remove the tubular neighborhood of the submanifold from the closed 4-manifold and glue the cobordism $P'_{r,k}$ along the resulting boundary. The resulting closed 4-manifold still satisfies $b_2<n$ as follows from the Mayer--Vietoris exact sequence,  and furthermore it contains $M_{r,0}$ as a separating submanifold. This contradicts a result of Kawauchi that  $M_{r,0}$  admits no such embedding (\cite[the proof of Theorem~2.5]{Kaw88}). Therefore, the desired claim follows.  
\end{proof}

As shown in the above proof, the simply connected closed $4$-manifold with $b_2=n$ in Corollary~\ref{sec:knotted:cor:exotic:embedded:3-manifold:some_4-manifold} can be taken as $\frac{n}{2}(S^2\times S^2)$ for each $n$.

\section{Diffeomorphism invariants of genus function type}\label{sec:genus type}
In this section, we introduce a notion of genus function type for diffeomorphism invariants of 4-manifolds and  give examples of such invariants which can detect exotic 4-manifolds. Furthermore, we show that any invariant of genus function type yields lower bounds for values of genus functions and shares many properties of genus functions. Throughout this section, \textit{we allow boundaries of $4$-manifolds to be disconnected.} 

\subsection{Definitions and examples}
To state the notion, we consider the category of 4-manifolds whose morphisms are orientation-preserving diffeomorphisms. Note that manifolds are oriented due to our convention. For a fixed ordered set $R$, we also consider the category of $R$-valued maps whose morphisms are maps between the source sets of $R$-valued maps.

 \begin{definition}\label{sec:genus type:def:invariant}
 	Let $R$ be an ordered set. A functor $G$ from the category of 4-manifolds to the category of $R$-valued maps is called an oriented diffeomorphism invariant of \textit{genus function type} if the following conditions hold.  
 	\begin{itemize}
 		\item [(i)] $G$ assigns to each $4$-manifold $X$ a map $G_X:H_2(X)\to R$, which we call the \textit{$G$ invariant} of $X$.  
 		\item [(ii)] Let $X,Y$ be arbitrary 4-manifolds. For any orientation-preserving embedding $\varphi:X\to Y$ and any class $\alpha\in H_2(X)$,  the induced homomorphism $\varphi_*:H_2(X)\to H_2(Y)$ satisfies $G_X(\alpha)\geq G_Y(\varphi_*(\alpha))$. 
 	\end{itemize}
 If furthermore $G$ satisfies the condition (ii) for any embedding $\varphi$ that does not necessarily preserves the orientations, then $G$ is called a $($unoriented$)$ diffeomorphism invariant of genus function type.  
 \end{definition}
 \begin{remark}\label{sec:genus type:remark:invariant}
 	It is straightforward from the condition (ii) that any oriented diffeomorphism $\varphi: X\to Y$ preserves the values of $G$ invariants, that is, the equality $G_X=G_Y\circ \varphi_*$ holds. The above $G$ is an oriented diffeomorphism invariant in this  sense. Similarly, if $G$ is an unoriented diffeomorphism invariant, then  any unoriented diffeomorphism preserves the values of $G$ invariants. 
 \end{remark}
  
 Clearly, genus functions of 4-manifolds form an unoriented diffeomorphism invariant of genus function type. 
 As in the case of genus functions, we use the following terminologies. 
 \begin{definition}Let $G$ be an oriented diffeomorphism invariant of genus function type, and let $X, Y$ be $4$-manifolds. 
 
  	(1) We will say that a continuous map $\varphi:X\to Y$ \textit{preserves the $G$ invariants} if the induced homomorphism $\varphi_*:H_2(X)\to H_2(Y)$ preserves the values of the $G$ invariants. 
  	
 	(2) The $G$ invariants of $X$ and $Y$ are called \textit{equivalent} if there exists a homeomorphism $\varphi:X\to Y$ that  preserves the $G$ invariants. If $\varphi$ furthermore preserves the orientations, then the $G$ invariants of $X$ and $Y$  are called \textit{orientedly equivalent}. 
 	
 	(3) The $G$ invariants of $X$ and $Y$ are called \textit{algebraically equivalent} if there exists an isomorphism $H_2(X)\to H_2(Y)$ that  preserves the intersection forms and $G$ invariants. 
 \end{definition}
 
Due to Remark~\ref{sec:genus type:remark:invariant}, $G$ invariants of orientedly diffeomorphic 4-manifolds are orientedly equivalent. Also, if $G$ invariants of two 4-manifolds are orientedly equivalent, then they are algebraically equivalent. Furthermore, if $G$ is an unoriented diffeomorphism invariant of genus function type, then $G$ invariants of (unorientedly)
 diffeomorphic 4-manifolds are equivalent. We also define ``torsion-free $G$ invariants'' and related terminologies in the same way as for the torsion-free genus functions. 

We introduce the following terminologies extracting properties of genus functions. 
 
\begin{definition}\label{sec:genus type:def:subadditivity}
	Let $G$ be an oriented diffeomorphism invariant of genus function type.  
	
	(1) Let $R$ be an additive magma (i.e.,  a set equipped with an additve binary operation $+$)  and suppose that the $G$ invariant of each 4-manifold is an $R$-valued map. We will say that \textit{$G$ is subadditive with respect to gluing} if $G$ has the following property. Let $X'$ be a 4-manifold obtained from a 4-manifold $X$ by attaching a 4-dimensional cobordism $P$, and let $\iota:X\hookrightarrow X'$ and $\iota_P:P\hookrightarrow X'$ be the inclusions. Then every $\alpha\in H_2(X)$ and $\beta\in H_2(P)$ satisfy  $G_{X'}(\iota_*(\alpha)+{\iota_P}_*(\beta))\leq G_X(\alpha)+G_P(\beta)$. 
	
	(2) We will say that \textit{$G$ is preserved under attachment of a $3$-handle $($resp.\ a $4$-handle$)$} if   for any 4-manifold $X'$ obtained from a 4-manifold $X$ by attaching a 3-handle (resp.\ a 4-handle), the inclusion $X\hookrightarrow X'$ preserves the $G$ invariants.  
	
	(3)  Let $P$ be a 4-dimensional cobordism. We will say that \textit{$G$ is preserved under attachment of $P$} if  for any 4-manifold $X'$ obtained from a 4-manifold $X$ by attaching the cobordism $P$,  the inclusion $X\hookrightarrow X'$ preserves the $G$ invariants.  
	
	(4) Let $Z$ be a $4$-manifold. 	We will say that \textit{$G$ is preserved under connected sum with $Z$} if for any 4-manifold $X$,  the inclusion  $H_2(X)\hookrightarrow  H_2(X\#Z)=H_2(X)\oplus H_2(Z)$ preserves the $G$ invariants of $X$ and $X\#Z$. In the case where $\partial Z$ is non-empty,  we will say that \textit{$G$ is preserved under boundary sum with $Z$} if for any 4-manifold $X$ with non-empty boundary and for any choice of a boundary sum $X\natural Z$, the inclusion  $X\hookrightarrow X\natural Z$ preserves the $G$ invariants.  
\end{definition}
 The genus function clearly has the above subadditivity and is preserved under attachment of a 3- and a 4-handle due to   Lemma~\ref{sec:genus:quasi-invertible:lem:genus:3-handle}. All the properties of the genus function used in the proofs of the results in Sections~\ref{sec:Genus functions and quasi-invertible cobordisms} and \ref{sec:genus:sum} are only those in Definitions~\ref{sec:genus type:def:invariant} and \ref{sec:genus type:def:subadditivity}.  It might be worth to remark that the genus function furthermore has the property that every 4-manifold $X$ satisfies $g_X(0)=g_{D^4}(0)$. 
 
We give other examples of oriented diffeomorphism invariants of genus function type which can detect exotic 4-manifolds. 
 \begin{example}\label{sec:genus type:example:invariant}
 	(1) For a 4-manifold $X$ and a class $\alpha\in H_2(X)$, let $E^{SW}_X(\alpha)$ be the maximum value of an  integer $|\langle K, \varphi_*(\alpha)\rangle|$, where $\varphi: X\to Z$ is an orientation-preserving embedding of $X$ into a closed 4-manifold $Z$ with $b_2^+>1$ having non-vanishing Seiberg--Witten invariant, and  $K$ is a Seiberg--Witten basic class of $Z$. We furthermore assume $Z$ is of Seiberg--Witten simple type in the sense of \cite{OzSz_ad} so that the adjunction inequalities (\cite{KM94}, \cite{OzSz_ad}) guarantee $E^{SW}_X(\alpha)<\infty$. In the case where $X$ admits no embedding into such a closed 4-manifold $Z$, we set $E^{SW}_X(\alpha)=-\infty$. Then, the resulting map $E^{SW}_X:H_2(X)\to \mathbb{Z}\cup\{-\infty\}$ gives an oriented diffeomorphism invariant of genus function type by varying $X$. 	
 	
 	Algebraic equivalent classes of this invariant ${E}^{SW}$ can detect various families of orientedly exotic 4-manifolds, as implicitly shown in \cite{AY12}, \cite{AY13}, \cite{AY_JSG14}, \cite{Y11}, \cite{Y14} and \cite{Y19Tr}. This invariant ${E}^{SW}$ clearly has the subadditivity with respect to gluing, but it is not preserved under attachment of a 4-handle, since $E^{SW}_{D^4}(0)=0$ and $E^{SW}_{S^4}(0)=-\infty$. We can define several variations of this invariant by using the Seiberg--Witten invariant and the Bauer--Furuta invariant. These invariants will be further discussed in \cite{Y_sta}.

 (2)  Recently, Ren and Willis \cite{RW} introduced an interesting new invariant of 4-manifolds called the lasagna $s$-invariant by using skein lasagna modules. The lasagna $s$-invariant  (with empty boundary link) of a 4-manifold $X$ is a map $s_X:H_2(X)\to \mathbb{Z}\cup \{-\infty\}$.  According to their result \cite[Proposition~1.11]{RW}, this invariant is an oriented diffeomorphism invariant of genus function type. Furthermore, they showed that this invariant can detect many orientedly exotic pairs of 4-manifolds. 
 
 (3) In the rest of this section, an immersed 2-sphere means the one whose singular points are only finitely many transverse double points. For a 4-manifold $X$ and a class $\alpha\in H_2(X)$, let $g^{\mathrm{dp}}_X(\alpha)$ be the minimum number of double points of an immersed 2-sphere in $X$ representing $\alpha$. In the case where $\alpha$ is not represented by any immersed 2-sphere, we set $g^{\mathrm{dp}}_X(\alpha)=\infty$. Then, the resulting map $g^{\mathrm{dp}}_X: H_2(X)\to \mathbb{Z}\cup\{\infty\}$ gives an unoriented diffeomorphism invariant of genus function type by varying $X$. Similarly to the genus function, determining  the values of $g^{\mathrm{dp}}_X$ has been an interesting problem, and it is in general quite hard (see e.g\ \cite{FS95}, \cite{Sch19}). 
 
 As is well-known, in the case where $X$ is simply connected, every $\alpha\in H_2(X)$ satisfies $g^{\mathrm{dp}}_X(\alpha)<\infty$ (see \cite[Remark 1.2.4.]{GS}). By contrast, in the case where $\pi_2(X)=0$, every non-zero class $\alpha\in H_2(X)$ satisfies $g^{\mathrm{dp}}_X(\alpha)=\infty$. In either case, every $\alpha\in H_2(X)$ satisfies $g^{\mathrm{dp}}_X(\alpha)\geq g_X(\alpha)$, since any double point can be eliminated by increasing the genus (see e.g.\ \cite[p.38]{GS}). Using this fact, we easily see that algebraic equivalent classes of $g^{\mathrm{dp}}$ can detect many of infinite families of exotic 4-manifolds in \cite{AY_JSG14}, \cite{Y11}, \cite{Y19Tr}.  
 The invariant $g^{\mathrm{dp}}$ clearly has the subadditivity with respect to gluing, and we can show that this invariant is preserved under attachment of a 3- and a 4-handle, similarly to the proof of Lemma~\ref{sec:genus:quasi-invertible:lem:genus:3-handle}. 
 
 (4) Let $G$ be an oriented diffeomorphism invariant of genus function type. We can modify $G$ into an unoriented invariant, for example, as follows. For a 4-manifold $X$, we set $\overline{G}_X=G_{\overline{X}}$. Note that the resulting orientation reversal  invariant $\overline{G}$ is also of genus function type. For a class $\alpha\in H_2(X)$, we set $|G|_X(\alpha)=\max\{G_X(\alpha),  G_{\overline{X}}(\alpha)\}$. The resulting $|G|$ is an unoriented diffeomorphism invariant of genus function type. For example, in the case of the invariant $E^{SW}$ defined in (1), the unoriented version $|E^{SW}|$ can distinguish various families of unorientedly exotic 4-manifolds as implicitly shown in \cite{AY13} and  \cite{Y14}. 
 \end{example}

\subsection{Properties}
Now we will show that invariants of genus function type share many properties of genus functions. 
Throughout this subsection, let $G$ be an oriented diffeomorphism invariant of genus function type, and assume that the invariant $G_X$ of  each 4-manifold $X$ is an $R$-valued map.   

We first show that $G$ yields a function $A_G:R\times \mathbb{Z}\to \mathbb{Z}\cup\{\infty\}$ that gives lower bounds for values of genus functions. 
For integers $g, n$ with $g\geq 0$, let $S(g,n)$ be the total space of the $D^2$-bundle over the closed surface of genus $g$ with Euler number $n$, and let $\gamma(g,n)$ denote a generator of $H_2(S(g,n))(\cong \mathbb{Z})$. 
To define the function $A_G$, we observe the following lemma. 
\begin{lemma}\label{sec:genus type:lem:genus:D2-bundle}
	$(1)$ $G_{S(g,n)}(-\gamma(g,n))=G_{S(g,n)}(\gamma(g,n))$ for any integers $g, n$ with $g\geq 0$. 
	
	$(2)$ For any integers $g_1, g_2, n$ with $g_2\geq g_1\geq 0$, the inequality $G_{S(g_2,n)}(\gamma(g_2,n))\geq G_{S(g_1,n)}(\gamma(g_1,n))$ holds. 
\end{lemma}
\begin{proof}As easily seen from a handlebody diagram, each $S(g,n)$ admits an orientation-preserving self-diffeomorphism whose induced isomorphism maps $\gamma(g,n)$ to $-\gamma(g,n)$. The claim (1) thus follows.  Due to the assumption $g_2\geq g_1$, we easily see that $\gamma(g_1,n)$ is represented by a closed surface of genus $g_2$ with self-intersection number $n$, which is obtained from that of genus $g_1$ by taking a connected sum with a locally embedded surface. We thus have an embedding of $S(g_2,n)$ into $S(g_1,n)$ whose  induced homomorphism maps $\gamma(g_2,n)$ to $\gamma(g_1,n)$ or to $-\gamma(g_1,n)$.  The claim (2) thus follows from the definition of genus function type and the claim (1). 
\end{proof}

For $(r,n)\in R\times \mathbb{Z}$, we define an element $A_G(r,n)$ of $\mathbb{Z}\cup\{\infty\}$  as follows. In the case where $r=G_{S(g,n)}(\gamma(g,n))$ for some integer $g\geq 0$, we set 
\begin{equation*}
	A_G(r,n)=\min\{g\in \mathbb{Z}_{\geq 0}\mid r=G_{S(g,n)}(\gamma(g,n))\}. 
\end{equation*}
In the case where $G_{S(g,n)}(\gamma(g,n))< r < G_{S(g+1,n)}(\gamma(g+1,n))$ for some integer $g\geq 0$, we set $A_G(r,n)=g+1$. In the case where $r<G_{S(g,n)}(\gamma(g,n))$ (resp.\ $r>G_{S(g,n)}(\gamma(g,n))$) for any  integer $g\geq 0$, we set $A_G(r,n)=0$ (resp.\ $A_G(r,n)=\infty$).  The resulting function $A_G:R\times \mathbb{Z}\to \mathbb{Z}\cup\{\infty\}$ is well-defined due to Lemma~\ref{sec:genus type:lem:genus:D2-bundle}. 

We prove that the function $A_G$ gives lower bounds for values of genus functions.
  \begin{proposition}\label{sec:genus type:prop:lower_bound} 
 	Every second homology class $\alpha$ of a $4$-manifold $X$ satisfies $g_X(\alpha)\geq A_G(G_X(\alpha),\alpha\cdot \alpha)$. 
  \end{proposition}
  In the above claim, $\alpha\cdot \alpha$ denotes the self-intersection number of $\alpha$.  

\begin{proof}[Proof of Proposition~$\ref{sec:genus type:prop:lower_bound}$]
	We set $g_0=g_X(\alpha)$. Then $\alpha$ is represented by a closed surface of genus $g_0$. Let $n$ be its self-intersection number. We have an embedding of $S(g_0,n)$ into $X$ whose induced homomorphism maps $\gamma(g_0,n)$ to $\alpha$, showing that  $G_{S(g_0,n)}(\gamma(g_0,n))\geq G_{X}(\alpha)$.  This inequality and the definition of $A_G(G_X(\alpha),\alpha\cdot \alpha)$ shows $g_0\geq A_G(G_X(\alpha),\alpha\cdot \alpha)$. The claim thus follows. 
\end{proof}

The inequality in Proposition~$\ref{sec:genus type:prop:lower_bound}$ can be sharp, depending on the invariant $G$. For example, in the case of the genus function, the equality always holds. For the invariant $E^{SW}$ defined in Example~\ref{sec:genus type:example:invariant}.(1)  by using the Seiberg--Witten invariant, we can easily determine all the values of  the function $A_{E^{SW}}$ using the adjunction inequality. Furthermore, we easily see that there are many examples of 4-manifolds and their homology classes for which the above bound is sharp.

We next show that $G$ is preserved under attachment of a (quasi-)invertible cobordism. 
  \begin{theorem}\label{sec:genus type:thm:invertible} 	
	$(1)$ $G$ is preserved under attachment of any invertible cobordism and also under boundary sum with any codimension-$0$ submanifold of $S^4$ having non-empty boundary.  
	
	$(2)$ If $G$ is preserved under attachment of a $4$-handle, then $G$ is preserved under connected sum with any codimension-$0$ submanifold of $S^4$. 
	
	$(3)$ Suppose $G$ is preserved under attachment of a $3$- and a $4$-handle. Then, $G$ is preserved under attachment of any quasi-invertible cobordism. Furthermore, for any $S^2$-link $L$, the invariant $G$ is preserved under boundary sum with any codimension-$0$ submanifold of $S^4(L)$ having non-empty boundary and also under connected sum with any codimension-$0$ submanifold of $S^4(L)$.   
\end{theorem}
 \begin{proof}We can prove these claims similarly to the proofs of Theorems~\ref{sec:genus:quasi-invertible:thm:genus preserving} and \ref{sec:genus:sum:thm:genus preserving}. See also Remark~\ref{sec:genus:quasi-invertible:remark:genus property}. 
 \end{proof}
 \begin{remark}\label{sec:genus type:remark:1-handle}
 (1) By Theorem~\ref{sec:genus type:thm:invertible}.(1) and Lemma~\ref{sec:quasi-invertible:lem:quasi-invertible:generate}, we see that $G$ is preserved under attachment of a 1-handle whose attaching region belongs to a connected boundary component and also under attachment of 2-handles along a strongly slice link with 0-framings.  
 
 (2) The conclusion of Theorem~\ref{sec:genus type:thm:invertible}.(3) still holds even if we replace the assumption on 3- and 4-handles with the weaker assumption that $G$ is preserved under attachment of $\natural_n(S^1\times D^3)$ along $n(S^1\times S^2)$ for any $n\geq 0$. 
\end{remark} 	

Using oriented diffeomorphism invariants of genus function type, we can characterize $S^4$ and its submanifolds as follows. 
 \begin{proposition}\label{sec:genus type:prop:characterization} 
 	$(1)$ A $4$-manifold is diffeomorphic to $S^4$  if and only if any oriented diffeomorphism invariant of genus function type  is preserved under connected sum with the $4$-manifold.  
 	
	$(2)$ A $4$-manifold with non-empty boundary admits an embedding into $S^4$  if and only if any oriented diffeomorphism invariant of genus function type  is preserved under boundary sum with the $4$-manifold.  
\end{proposition}
\begin{proof}The ``only if'' part of the claim (1) is straightforward, and that of (2) is Theorem~\ref{sec:genus type:thm:invertible}. In the rest, we prove the ``if'' part of each claim. 
	
	(1) For a 4-manifold $X$ and a class $\alpha\in H_2(X)$, we set $D^{S^4}_X(\alpha)=0$ (resp.\ $D^{S^4}_X(\alpha)=\infty$) in the case where $X$  is diffeomorphic (resp.\ not diffeomorphic) to $S^4$. 
	Then, the resulting map $D^{S^4}_X:H_2(X)\to \{0, \infty\}$ gives an oriented diffeomorphism invariant of genus function type. 
   Now let $Z$ be a 4-manifold such that any oriented diffeomorphism invariant of genus function type  is preserved under connected sum with $Z$.  Then,  $Z$ satisfies $D^{S^4}_{S^4\#Z}(0)=D^{S^4}_{S^4}(0)=0$, and hence $Z$ is  diffeomorphic to $S^4$. Therefore, the claim (1) follows. 
   
    (2) For a 4-manifold $X$ and a class $\alpha\in H_2(X)$, we set $E^{S^4}_X(\alpha)=0$ (resp.\ $E^{S^4}_X(\alpha)=-\infty$) in the case where $X$  admits an embedding (resp.\ no embedding) into $S^4$ that preserves the orientations. 
 Then, the resulting map $E^{S^4}_X:H_2(X)\to \{-\infty, 0\}$ gives an oriented diffeomorphism invariant of genus function type. Now let $Z$ be a 4-manifold with non-empty boundary such that any oriented diffeomorphism invariant of genus function type  is preserved under boundary sum with $Z$.  Then, $Z$ satisfies $E^{S^4}_{D^4\natural Z}(0)=E^{S^4}_{D^4}(0)=0$, and hence $Z$  admits an embedding into $S^4$. Therefore, the claim (2) follows. 
 \end{proof}
 \begin{remark}\label{sec:genus type:remark:characterization} 
 (1) Similarly to the above proof, we see that for any closed 4-manifold, its diffeomorphism type and oriented diffeomorphism type are determined by oriented diffeomorphism invariants of genus function type similar to $D^{S^4}$. 
 
 (2) For a closed (possibly disconnected) 3-manifold $M$, we define an oriented diffeomorphism invariant $E^{I\times M}$ in the same way as the above $E^{S^4}$. Using this invariant similarly to the above proof, we can prove that a cobordism with negative boundary $M$ admits an embedding into $I\times M$ that preserves the orientations if any oriented diffeomorphism invariant of genus function type  is preserved under attachment of the cobordism. 
 \end{remark}

We show that the stabilities of algebraic inequivalences also hold for $G$ invariants and their torsion-free versions. More precisely, we prove the following theorem for attachments of (quasi-)invertible cobordisms.  

\begin{theorem}\label{sec:genus type:thm:quasi-invertible:many_results}
	$(1)$ Theorems~$\ref{sec:genus:quasi-invertible:thm:stablity:H2zero}$,  $\ref{sec:genus:quasi-invertible:thm:stablity:non-degenerate}$ and $\ref{sec:genus:quasi-invertible:thm:stablity:non-degenerate:torsion-free}$ and Corollaries~$\ref{sec:genus:quasi-invertible:cor:inclusion:isomorphism}$, $\ref{sec:genus:quasi-invertible:cor:inclusion:invertible homology cobordism}$ and $\ref{sec:genus:quasi-invertible:cor:algebraic equivalence:sufficient}$ hold even when ``genus functions'' in each statement is replaced with ``$G$ invariants'', in the case where the cobordism $P$ in Subsection~$\ref{subsec:Genus functions and quasi-invertible cobordisms:behavior}$ is invertible. 
	
	$(2)$ Suppose that $G$ is preserved under attachment of a $3$- and a $4$-handle. Then,  Theorems~$\ref{sec:genus:quasi-invertible:thm:stablity:H2zero}$,  $\ref{sec:genus:quasi-invertible:thm:stablity:non-degenerate}$ and $\ref{sec:genus:quasi-invertible:thm:stablity:non-degenerate:torsion-free}$ and Corollaries~$\ref{sec:genus:quasi-invertible:cor:inclusion:isomorphism}$ and $\ref{sec:genus:quasi-invertible:cor:algebraic equivalence:sufficient}$ hold even when ``genus functions'' in each statement is replaced with ``$G$ invariants''. 
\end{theorem}
\begin{proof}Due to Theorem~\ref{sec:genus type:thm:invertible}, we see that the left inequality in Theorem~\ref{sec:genus:quasi-invertible:thm:estimate} and Corollary~\ref{sec:genus:quasi-invertible:cor:genus preserving:torsion-free} holds for $G$ invariants under the above assumptions by the identical argument. Therefore, we can prove the claims by the identical arguments. 
\end{proof}
For boundary sums and connected sums, we prove the following theorem. 

\begin{theorem}\label{sec:genus type:thm:sum:many_results}
	$(1)$ For boundary sums, Theorems~$\ref{sec:genus:sum:thm:equivalent:sum}$, $\ref{sec:genus:sum:thm:equivalent:sum:torsion-free}$, $\ref{sec:genus:sum:thm:equivalent:non-degenerate}$ and $\ref{sec:genus:sum:thm:equivalent:non-degenerate:torsion-free}$ hold even when ``genus functions'' in each statement is replaced with ``$G$ invariants'', in the case where $L$, $L_1$ and $L_2$  are empty sets.  
	
	$(2)$ Suppose that $G$ is preserved under attachment of a $4$-handle. Then, Theorems~$\ref{sec:genus:sum:thm:equivalent:sum}$, $\ref{sec:genus:sum:thm:equivalent:sum:torsion-free}$, $\ref{sec:genus:sum:thm:equivalent:non-degenerate}$ and $\ref{sec:genus:sum:thm:equivalent:non-degenerate:torsion-free}$   hold even when ``genus functions'' in each statement is replaced with ``$G$ invariants'', in the case where $L$, $L_1$ and $L_2$ are empty sets.  If furthermore $G$ is preserved under  attachment of a $3$-handle, then these claims and Theorem~$\ref{sec:genus:sum:thm:closed:algebraically equivalent}$ for $G$ invariants hold without any condition on $L$, $L_1$ and $L_2$. 
\end{theorem}
\begin{proof}Due to Theorem~\ref{sec:genus type:thm:invertible}, we see that the left inequality in the second claim of Theorem~\ref{sec:genus:sum:thm:estimate} and the second claim of Corollary~\ref{sec:genus:sum:cor:genus preserving:torsion-free} hold for $G$ invariants by the identical arguments. Therefore, we can prove the claim (1) by the identical arguments. We can similarly prove the claim (2) under the above assumptions. 
\end{proof}

The above theorems imply that if algebraic equivalent classes of $G$ invariants can detect a family of orientedly exotic 4-manifolds, then they can detect infinitely many families of orientedly exotic 4-manifolds. 
\begin{corollary}\label{sec:genus type:cor:detect_infinite}
	Suppose $G$ admits a family $\{X_i\}_{i\in A}$ of pairwise orientedly homeomorphic $4$-manifolds whose $G$ invariants are pairwise algebraically inequivalent. Then the following hold. 
	
	$(1)$ If the members $X_i$ have non-empty boundary, then there exists a doubly indexed family $\{X_{i,j}\}_{(i,j)\in A\times \mathbb{Z}}$ of $4$-manifolds having the following properties. 
\begin{itemize}
	\item For each fixed $j\in \mathbb{Z}$, the members of $\{X_{i,j}\}_{i\in A}$ are pairwise orientedly exotic and have pairwise algebraically inequivalent $G$ invariants. 
	\item For each fixed $i\in A$, the members of  $\{X_{i,j}\}_{j\in \mathbb{Z}}$ are pairwise not homeomorphic. 
\end{itemize}	
	 
	$(2)$ If $G$ is preserved under attachment of a $4$-handle, then the conclusion of the above $(1)$ holds even in the case where each $X_i$ has no boundary.  
	
	$(3)$ If $G$ is preserved under attachment of a $3$- and a $4$-handle, and each $X_i$ has no boundary, then the conclusion of the above $(1)$ still holds, and furthermore the members $X_{i,j}$ can be chosen to be closed $4$-manifolds.   
\end{corollary}
\begin{proof}We first consider the case (1). By Theorems~\ref{sec:genus type:thm:quasi-invertible:many_results} and \ref{sec:genus type:thm:sum:many_results}, we can construct many examples of desired families, for example, by taking boundary sums of the 4-manifolds $X_i$ with codimension-0 submanifolds of $S^4$. See also Lemma~\ref{sec:genus:sum:lem:submanifold:generate}. 	In the (2) case, we can similarly construct many examples of desired families by taking connected sums instead of boundary sums. In the (3) case, we can give desired families of closed 4-manifolds by taking connected sums with $S^4(L)$ for an $S^2$-link $L$ in $S^4$. 
\end{proof}

As shown in Remark~\ref{sec:genus type:remark:characterization}, every exotic pair of closed 4-manifolds admits an oriented  diffeomorphism invariant of genus function type that can detect the pair.  By contrast, we show that there exist various infinite families of exotic 4-manifolds which can not be detected by any such invariant. 

\begin{theorem}\label{sec:genus type:thm:2-handlebody:undetectable}For any $4$-dimensional $2$-handlebody $X$, there exist infinitely many pairwise exotic $4$-manifolds HIHC-equivalent to  $X$ such that, for any oriented diffeomorphism invariant $G$ of genus function type, their $G$ invariants are pairwise equivalent. 
\end{theorem}

\begin{theorem}\label{sec:genus type:thm:symplectic:undetectable}
	Suppose  that a $4$-manifold $X$ with connected boundary admits an embedding into a symplectic $4$-manifold with weakly convex boundary. Then, there exist infinitely many pairwise exotic $4$-manifolds HIHC-equivalent to  $X$ such that, for any oriented diffeomorphism invariant $G$ of genus function type, their $G$ invariants are pairwise equivalent and are  algebraically equivalent to that of $X$. 
\end{theorem}

\begin{proof}[Proof of Theorems~$\ref{sec:genus type:thm:2-handlebody:undetectable}$ and $\ref{sec:genus type:thm:symplectic:undetectable}$]Due to Theorem~\ref{sec:genus type:thm:invertible}, the proofs are identical with those of Theorems~$\ref{sec:exotic:thm:2-handlebody:infinite}$ and $\ref{sec:exotic:thm:symplectic:infinitely}$.  
\end{proof}

We close this section by making a few remarks. 
 \begin{remark}
 	$(1)$ The topological genus function of a 4-manifold is defined similarly to the (smooth) genus function by using locally flat embedded surfaces instead of smoothly embedded surfaces, and this function also has long been studied, see  \cite{Ray} and the references therein, for example. Topological genus functions clearly form an (unoriented) diffeomorphism invariant of genus function type. Due to the transversality theorem for locally flat proper submanifolds (see \cite[Theorem~1.5]{BKKPR} and \cite{Ray}),  the proof of  Lemma~\ref{sec:genus:quasi-invertible:lem:genus:3-handle} shows that topological genus functions are preserved under attachment of a 3- and a 4-handle.   	
 	 Since two homeomorphic 4-manifolds have equivalent topological genus functions, they cannot be used for detecting exotic 4-manifolds.  However, results in this section for topological genus functions  might be useful for  determining their values  and also for detecting homeomorphism types by using stabilities of algebraic inequivalences. 
 	
 	$(2)$ Due to our convention, we defined the notion of genus function type only for (oriented) diffeomorphism invariants of compact 4-manifolds, but in the same manner, we can define the notion more generally for invariants of 4-manifolds which are not necessarily compact. For example, the genus function and the invariants $E^{SW}$ and  $g^{\mathrm{dp}}$ defined in Example~\ref{sec:genus type:example:invariant} are of genus function type in the latter sense. For such an invariant $G$, we can show that the $G$ invariants of a 4-manifold and its interior are algebraically equivalent similarly to the proof of Theorem~\ref{sec:genus:quasi-invertible:thm:genus preserving}. Thus, for any two 4-manifolds having algebraically inequivalent $G$ invariants,  their interiors are not orientedly diffeomorphic. On the other hand, we can also define the notion of genus function type for (oriented) diffeomorphism invariants of compact 4-manifolds having connected (possibly empty) boundary. For such invariants, our results in this section still hold for 4-manifolds with connected (possibly empty) boundary,  as easily seen from the proofs. 
 \end{remark}


\end{document}